\documentclass[11pt,a4paper]{article}

\usepackage[a4paper,top=1in,bottom=1in,left=1in,right=1in]{geometry} % For sufficient margins on A4
\usepackage{setspace}
\usepackage{natbib}
\usepackage{amsmath,amssymb,amsfonts,amsthm,bm}
\usepackage{graphicx,color,xcolor}
\usepackage{multirow,booktabs,array}
\usepackage{float,subfigure,rotating}
\usepackage{algorithm,algorithmicx,algpseudocode}
\usepackage{url,hyperref}
\usepackage{lipsum}
\usepackage{indentfirst}
\usepackage{authblk}
\usepackage{textcomp}
\usepackage{threeparttable}
\usepackage{pdflscape}
\usepackage{enumerate}
\hypersetup{
    colorlinks=true,
    linkcolor=blue,
    anchorcolor=blue,
    citecolor=blue,
    urlcolor=blue,
    CJKbookmarks=true
}

\newtheorem{theorem}{Theorem}
\newtheorem{assumption}{Assumption}
\newtheorem{lemma}{Lemma}
\newtheorem{remark}{Remark}
\newtheorem{corollary}{Corollary}
\newtheorem{proposition}{Proposition}

 \newcommand{\bP}{\mathbf{P}} \newcommand{\bA}{\mathbf{A}}
\newcommand{\bu}{\mathbf{u}} \newcommand{\ba}{\mathbf{a}} \newcommand{\bU}{\mathbf{U}}
\newcommand{\bV}{\mathbf{V}} \newcommand{\bD}{\mathbf{D}} \newcommand{\by}{\mathbf{y}}
\newcommand{\bX}{\mathbf{X}} \newcommand{\bI}{\mathbf{I}} 
 \newcommand{\bz}{\mathbf{z}} \newcommand{\bw}{\mathbf{w}}
 \newcommand{\be}{\mathbf{e}} \newcommand{\bmu}{\boldsymbol{\mu}}
\newcommand{\btheta}{\boldsymbol{\theta}} \newcommand{\bbeta}{\boldsymbol{\beta}}
 \newcommand{\bpsi}{\boldsymbol{\psi}}
\newcommand{\bPsi}{\boldsymbol{\Psi}} \newcommand{\beps}{\boldsymbol{\epsilon}}

\renewcommand{\hat}{\widehat} \renewcommand{\tilde}{\widetilde}
 \DeclareMathOperator*{\argmin}{argmin}
\def\Var{\operatorname{Var}} \def\tr{\operatorname{tr}}

\makeatletter \@addtoreset{equation}{section} \makeatother

\author{Jingfu Peng}
\affil{Yau Mathematical Sciences Center, Tsinghua University}

\date{\today}  % Optional: remove date

\title{Mallows-type model averaging: Non-asymptotic analysis\\ and all-subset combination}

\begin{document}

\maketitle

%Abstract
\begin{abstract}
%\begin{spacing}{1.5}
Model averaging (MA) and ensembling play a crucial role in statistical and machine learning practice. When multiple candidate models are considered, MA techniques can be used to weight and combine them, often resulting in improved predictive accuracy and better estimation stability compared to model selection (MS) methods. In this paper, we address two important problems in combining least squares estimators from both theoretical and practical perspectives. We first establish several oracle inequalities for least squares MA via minimizing a Mallows' $C_p$ criterion under an arbitrary candidate model set. Compared to existing studies, these oracle inequalities yield faster excess risk and directly imply the asymptotic optimality of the resulting MA estimators under milder conditions. Moreover, we consider candidate model construction and investigate the problem of optimal all-subset combination for least squares estimators, which remains an important yet not fully understood topic in the existing literature. We show that there exists a fundamental limit to achieving the optimal all-subset MA risk. To attain this limit, we propose a novel Mallows-type MA procedure based on a dimension-adaptive $C_p$ criterion. The implicit ensembling effects of several MS procedures are also revealed and discussed. We conduct several numerical experiments to support our theoretical findings and demonstrate the effectiveness of the proposed Mallows-type MA estimator.
%\end{spacing}
\end{abstract}


\tableofcontents

\addtocontents{toc}{\setcounter{tocdepth}{2}}
\section{Introduction}\label{sec:intro}

Model averaging (MA or ensemble learning) has been an active research topic in statistics, econometrics, and machine learning for over 30 years, with numerous approaches proposed for combining models to support decision-making. These include forecast combination \citep{Bates01121969}, Bayesian MA (BMA) \citep[see][and the references therein]{Draper1995, chatfield1995model, Hoeting1999}, bagging \citep{Breiman1996Bagging}, stacking \citep{WOLPERT1992241, Breiman1996Stacked}, random forests \citep{Breiman2001Forest}, AIC/BIC-based weighting \citep{Buckland1997, Hjort2003, Liang2011}, adaptive regression by mixing \citep{Yang2001, yang2004aggregating, Yuan2005, Wang2014}, and exponential weighting \citep{George1986, Leung2006}, among many other useful techniques. These classical MA methods have been successfully applied to a wide range of problems, such as mitigating model selection (MS) uncertainty (e.g., by BMA), constructing minimax adaptive estimators (e.g., by ARM), improving risk performance over MS \citep[see, e.g.,][]{Peng2022improvability, Xu2022model, le2022model, Chen2023error}, and conducting variable importance diagnostics in high-dimensional learning \citep[see, e.g.,][]{Ye2018}. For a comprehensive review of MA and ensemble learning, see \cite{claeskens2008model}, \cite{Wang2009}, \cite{MR3887687}, and \cite{Sagi2018}.

One of the most fundamental problems in MA is the combination of least squares estimators. In this setting, multiple candidate regression models are estimated using the least squares method, and a data-driven weighting scheme is then designed to aggregate these estimators based on the same dataset. To the best of our knowledge, an early but not very well-known study on least squares MA was conducted by \cite{Blaker1999}, where two nested models were combined by minimizing a Mallows' $C_p$ criterion \citep{Mallows1973}. This work is one of the earliest applications of what is now referred to as Mallows MA (MMA) methods. \cite{Leung2006} proposed an exponential weighting method to achieve the optimal MS risk over a collection of least squares estimators. In the context of multiple nested models, \cite{Hansen2007least} established that the MMA estimator achieves an asymptotic optimality (AOP), i.e., it is asymptotically equivalent to the optimal convex combination of candidate least squares estimators with discretized weights in terms of statistical loss. Later, the AOP property has become a predominant justification for the superiority of least squares MA approaches. Under certain restrictions on the candidate models, \cite{WAN2010277} established the MMA's AOP for general non-nested candidate models with continuous weights. A similar setting was adopted by \cite{Zhang_2021}, in which milder and more interpretable assumptions for the AOP were given. Building upon the least squares MA framework, various Mallows-type MA strategies have been developed to combine more general regression estimators \citep[see, e.g.,][]{HANSEN201238, ZHANG201382, Zhang2016jlm, Zhang2020Parsimonious, Ando2014, Ando2017, Cheng2015,  liu2015distribution, LIAO201935, FANG2022219, Li2022AdaBoost, SUN20231355, yu2024unified, Zhu2024Stability, Chen2024Forests, TU2025105999}.

Although Mallows-type MA approaches with AOP properties have been formulated within various general modeling frameworks, two important aspects of their theoretical foundation and practical implementation in the least squares MA setting continue to pose open challenges.

\emph{Is there any finite-sample performance guarantee for the MMA estimators?} In MA approaches with AOP properties, while asymptotic theory provides rigourous risk characterization as $n \to \infty$, it offers limited performance guarantees in the realistic settings where the sample size $n$ is finite. Indeed, in the MS context, \cite{Kabaila_2002} demonstrated that while AIC-based MS estimators can achieve AOP in terms of MS loss within a typical nested framework, they may perform inefficiently in finite sample settings; see also \cite{Yang2005, Yang_2007} for related discussions. As remarked in the first paragraph below Theorem~4 of \cite{Wang2009}, the footnote on page 278 of \cite{WAN2010277}, and Remark 6 in \cite{Liao_Tsay_2020}, such a decoupling between asymptotic theory and finite-sample performance may also occur for the MA estimators with AOP properties.

To the best of our knowledge, the only work on the finite-sample risk performance of MMA with general candidate models is given in Proposition 7.2 of \cite{Bellec2018}. It established an oracle inequality for MMA under Gaussian errors. However, \cite{Bellec2018}'s result shows the excess risk of MMA to the optimal MA risk converges at a rate no faster than $n^{-1/2}$, regardless of the number of candidate models. As remarked in Section~7 of \cite{Bellec2018}, it remains unclear whether his bound is tight, particularly when the number of candidate models is small. Therefore, in the setting where least squares estimators from general candidate models are combined, whether a sharper finite-sample risk bound of MMA exists remains an open question.

\emph{How to construct candidate model set for least squares MA estimators?} The asymptotic analysis in \cite{WAN2010277, Zhang_2021} and the oracle inequalities established in \cite{Bellec2018} provide valuable insights into MMA with general candidate models. However, these works do not investigate how the construction of the candidate model set influences the resulting MA estimators. Consider a typical setting of least squares MA, where the true regression function follows a linear model with $p$ regressors, and candidate models are constructed using different subsets of these regressors. In this setting, the ideal choice of candidate model set consists of all subsets of the $p$ regressors, resulting in $2^p$ least squares estimators. The optimal MA risk over these $2^p$ estimators should be regarded as the target for least squares MA.

In the existing literature, the achievability of the optimal all-subset combination remains largely an open problem. Some approaches, such as the two-stage least squares MA methods \citep[see, e.g.,][]{ELLIOTT2013357, Lee2020}, have been developed in an attempt to approximate this ideal risk of MA. However, their theoretical optimality has not been proven. \citet{Zhu2023Scalable} proposed a scalable MA method that aims to achieve the optimal all-subset MA risk under both orthogonal and general regression settings. Its theoretical guarantees depend on specific regularity conditions imposed on the optimal MA risk and the dimensionality. More recently, \cite{Peng2024Optimality} demonstrated that if the relative importance of regressors is largely known, then the optimal all-subset combination can be achieved by the nested MMA. In contrast, when the order of regressors is completely unknown, no method can attain the optimal all-subset MA risk. However, without prior ordering information of regressors, \cite{Peng2024Optimality} does not provide upper bounds on how closely an estimator can approach the optimal all-subset MA risk.

\subsection{Contributions}

In this paper, we address the aforementioned challenges in least squares MA. First, we establish several oracle inequalities for least-squares MMA estimators based on general candidate model set. These inequalities are derived under the finite fourth-moment condition on random error terms, as imposed in \cite{WAN2010277} and \cite{Zhang_2021}. Compared to the classical AOP theory, our risk bounds hold for any sample size, providing a finite-sample performance guarantee for the MMA estimators relative to the optimal convex combination of candidates. By letting $n \to \infty$, our oracle inequalities lead to milder and comparable conditions for AOP in risk compared to the loss-based AOP results in \cite{WAN2010277} and \cite{Zhang_2021}, respectively.

Second, from a technical perspective, we employ a \emph{shifted empirical process} method \citep[see, e.g.,][]{Baraud2000,Wegkamp2003} to obtain a non-exact oracle inequality, which yields faster convergence rate compared to that in \cite{Bellec2018}. As a byproduct, our established risk bounds also imply the achievability of the optimal MA risk with all-nested models under weaker conditions on the random error terms, relaxing the sub-Gaussian assumption in \cite{Peng2024Optimality}.

Third, we establish the fundamental limits of achieving the optimal all-subset MA risk. We show that even in the setting where the regressors are orthogonal and random error is Gaussian, the minimax risk ratio of any regression estimator relative to the optimal all-subset MA risk cannot converge to $1$ as $n \to \infty$. Specifically, when the dimension of true model $p$ is fixed, which corresponds to the typical parametric setting, the minimax risk ratio can be strictly larger than 1. Moreover, if $p$ diverges to infinity, the minimax risk ratio is lower bounded by a rate of $2\log p$.

Forth, under a similar setting as that in the lower bound, we propose a dimension adaptive Mallows-type MA to combine least squares estimators. We show that the resulting MA estimator attains the optimal convergence rate towards the risk of the optimal all-subset MA. To the best of our knowledge, this is the first MA estimator with a theoretically provable optimality in achieving the best all-subset combination, without imposing rate restrictions on the optimal MA risk. The connections between all-subset MA, soft/hard-thresholding estimators \citep{Donoho1994Ideal}, and the risk inflation MS criterion \citep{Foster1994} are also discussed. Simulation results further support our theoretical findings.

\subsection{Other related work}

This paper builds upon the line of insightful research by \cite{Hansen2007least}, \cite{WAN2010277}, and \cite{Zhang_2021}, which focuses on deriving the optimal convex combinations of estimators in a fixed design setting. Beyond this viewpoint, several other lines of research on MA have also been explored in the existing literature.

\emph{Aggregation of general estimation procedures.} Aggregation is a long-standing topic in statistical learning theory, aiming to combine general statistical procedures/estimators under various weight constraints \citep[see, e.g.,][]{Yang2000Mixing, Nemirovski2000, Catoni2004, Tsybakov2003, Wang2014}. The optimality of aggregation is measured by a minimax regret, i.e., the minimax gap between the aggregated estimator and the optimal aggregated risk over general candidate procedures and true models. When candidate estimators of the regression mean vector $\bmu$ have the affine form $\hat{\bmu}_m = \mathbf{A}_m \by + \mathbf{b}_m,m=1,\ldots,M_n$, some aggregation strategies have been proposed \citep{Dalalyan2012, Chernousova2013, Dai2014, Golubev2016, Bellec2018, Bellec2020}, and the minimax regret optimality has been studied \citep{Dalalyan2012, Bellec2018}. Although incorporating the deterministic intercepts $\mathbf{b}_1,\ldots,\mathbf{b}_{M_n}$ offers greater flexibility for candidate construction and also enables an application of the minimax lower bounds from \cite{Tsybakov2003, Rigollet2011}, this setup does not capture the fundamental difficulty of convex aggregation of $\mathbf{A}_1 \by, \ldots, \mathbf{A}_{M_n} \by$, which is more common in practice. For example, all estimators presented in Section~1.2 of \cite{Dalalyan2012} have the linear form without intercept terms. Our work focuses on a fundamental case in which each $\mathbf{A}_m$ is a projection matrix, and we establish a minimax lower bound for attaining the optimal all-subset MA risk in terms of risk ratio, along with several matching upper bounds.

%Particularly, the paper on aggregating affine estimator is closely related to the theme of this paper, where the goal of aggregation is mainly achievie the optimal single model in the minimax rate. However, these paper is has two different points compared to that in this paper. First, their candidate estimatr to be combined is generalize to affinie estimator $\mathbf{A} \by + \mathbf{b}$ \citep{Dalalyan2012, Bellec2018, Chernousova2013, Bellec2020, Golubev2016}. This brodness althouugh can be used to adopt the prevsious lower bound in \citep{Tsybakov2003}, but cannot capture the essential difficulty of combining more basic candidate estimator with only linear smoothers. As illustrated in \cite{Dalalyan2012}, the exaples presented in their paper is all the linear smoother. However, even in the most basic setting of combing the least squares estimator, the lower bound in the minimax sense is still remain open. This problem may be challeging, but the lower bound on achieving the all-subset least squares estimator partially address this problem. But about general aggregation of least squares estimator, the optimal aggregation rate in the sense of \cite{Tsybakov2003} still remains open.

%There some statement about the advantage of aggregation. \cite{Mourtada2023}

\emph{Ensemble learning under random design regression.} Recently, there has been growing interest in the asymptotic risk analysis of ensemble estimators in high-dimensional random design regression \citep[see, e.g.,][]{LeJeune2020, ando2023high, Bellec2024Corrected, Du2023Subsample, Du2024Extrapolated, patil2024asymptotically, wu2023ensemble}. The construction of candidate models and theoretical objectives in these studies differ from our work. For instance, \cite{ando2023high} combines minimum-norm least squares estimators from different subsets of regressors and samples. While an asymptotic expression for the out-of-sample prediction risk of the MA estimator is derived using random matrix theory, the study does not provide a theory for estimating the optimal weights or constructing the candidate model set, both of which are addressed in our work. Similarly, \cite{Bellec2024Corrected} considers a setting where penalized least squares estimators are constructed from different subsets of the sample drawn from the entire dataset, and these estimators are combined using equal weights. In contrast, our approach treats different subsets of regressors as candidate models and determines the weights in a data-driven manner.

\subsection{Organization}

We formally set up the regression problem and introduce the Mallows-type MA estimators in Section~\ref{sec:setup}. Section~\ref{sec:general} presents oracle inequalities for combining least squares estimators from general linear subspaces, with a brief discussion of their implications in the nested candidate model setup. In Section~\ref{sec:all-subset}, we establish both lower and upper bounds for achieving the optimal risk of all-subset MA. Section~\ref{sec:simu} provides numerical results, followed by a discussion in Section~\ref{sec:discuss}. The proofs of the main results are provided in the Appendix.

\section{Problem setup}\label{sec:setup}

\subsection{Setup and notation}\label{subsec:notation}

We study the problem of estimating an unknown mean vector $\bmu\triangleq(\mu_1,\ldots,\mu_n)^{\top} \in \mathbb{R}^n$ from noisy observations
\begin{equation}\label{eq:model}
  \by = \bmu +\beps,
\end{equation}
where $\by \triangleq (y_1,\ldots,y_n)^{\top} \in \mathbb{R}^n$, and $\beps \triangleq (\epsilon_1,\ldots,\epsilon_n)^{\top} \in \mathbb{R}^n$ consists of independent random errors with mean zero and variance $\sigma^2$. We assume that the random errors $\epsilon_i$ satisfy the following fourth-moment condition.
\begin{assumption}\label{ass:4_moment}
  The random error terms satisfy $\mathbb{E}\epsilon_i^4 \leq \nu < \infty$, where $\nu$ is a positive constant.
\end{assumption}
The objective is to construct an estimator $\hat{\bmu}$ of $\bmu$ based on the observation $\by$. For any estimator $\hat{\bmu}$, its performance is assessed by the normalized squared loss $L_n(\hat{\bmu}, \bmu) \triangleq n^{-1}\| \hat{\bmu} - \bmu \|^2$ and the corresponding squared risk $R_n(\hat{\bmu}, \bmu) \triangleq \mathbb{E} L_n(\hat{\bmu}, \bmu)$, where $\| \cdot \|$ denotes the Euclidean norm.

Since the true mean vector $\bmu$ may reside in an unknown subspace of $\mathbb{R}^n$, we consider a collection of $M_n$ candidate subspaces, $\mathbb{V}_1,\ldots,\mathbb{V}_{M_n}$, where each $\mathbb{V}_{m}$ is a linear subspace of $\mathbb{R}^n$ with dimension $k_m$. Given $\mathbb{V}_m$, we estimate $\bmu$ using the least squares estimator
\begin{equation*}\label{eq:mu_hat}
  \hat{\bmu}_{m} = \bP_{m}\by \triangleq \argmin_{\bmu \in \mathbb{V}_m}\| \by - \bmu \|^2,
\end{equation*}
where $\bP_{m} $ is the projection matrix on $\mathbb{V}_m$. This setup also encompasses an important case considered in the least squares MA literature \citep[see, e.g.,][]{Hansen2007least, WAN2010277, Zhang_2021}, where $\mu_i=\sum_{j=1}^{\infty}\beta_j x_{ij}$ and $\hat{\bmu}_m$ denotes the least-squares estimator that projects onto the span of a finite subset of the covariates. Let $\bw\triangleq\left(w_{1}, \ldots, w_{M_{n}}\right)^{\top}$ be a weight vector in $\mathcal{W} \triangleq \{\bw \in [0,1]^{M_n}: \sum_{m=1}^{M_n}w_m=1 \}$. The least squares MA estimator of $\bmu$ based on the candidate model set $\mathcal{M} \triangleq \{ \mathbb{V}_1,\ldots,\mathbb{V}_{M_n} \}$ is defined as
\begin{equation}\label{eq:ma}
\hat{\bmu}_{\bw|\mathcal{M}}\triangleq\sum_{m=1}^{M_{n}} w_{m} \hat{\bmu}_{m}=\bP(\bw)\by,
\end{equation}
where $\bP(\bw)\triangleq \sum_{m=1}^{M_{n}} w_{m}\bP_{m}$, and the subscript $\bw|\mathcal{M}$ emphasizes the dependence of the MA estimator on the candidate model set $\mathcal{M}$.

The performance of the MA estimator (\ref{eq:ma}) is measured by $L_n(\hat{\bmu}_{\bw|\mathcal{M}}, \bmu)$ and $R_{n}(\hat{\bmu}_{\bw|\mathcal{M}}, \bmu)$. From the perspective of risk minimization, the optimal MA risk is defined as
\begin{equation}\label{eq:optimal_MA_risk}
  R_{n}(\hat{\bmu}_{\bw^*|\mathcal{M}}, \bmu) \triangleq \min_{\bw \in \mathcal{W}}R_{n}(\hat{\bmu}_{\bw|\mathcal{M}}, \bmu).
\end{equation}
This represents the lowest possible MA risk given the candidate models $\mathcal{M}=\{\mathbb{V}_1,\ldots,\mathbb{V}_{M_n}\}$ at the true mean vector $\bmu$. The goal of constructing specific MA procedures can be divided into two parts: (\romannumeral1) estimating the weight vector $\tilde{\bw}$ based on the data and showing that its risk $\mathbb{E}L_{n}(\hat{\bmu}_{\tilde{\bw}|\mathcal{M}}, \bmu)$ approaches $R_{n}(\hat{\bmu}_{\bw^*|\mathcal{M}}, \bmu)$ as closely as possible; and (\romannumeral2) designing an appropriate set of candidate models such that $R_{n}(\hat{\bmu}_{\bw^*|\mathcal{M}}, \bmu)$ is both efficient and achievable.

In this paper, we investigate the aforementioned two goals from a theoretical perspective. We use the notation $\lesssim$ for comparison of two positive sequences, where $a_n \lesssim b_n$ denotes $a_n = O(b_n)$. Also, $a_n\asymp b_n$ denotes both $a_n\lesssim b_n$ and $b_n\lesssim a_n$. We use $a_n \sim b_n$ to denote $\lim_{n \to \infty}a_n/b_n= 1$. For any two real numbers $a$ and $b$, we use notation $a \wedge b = \min(a,b)$ and $a \vee b = \max(a,b)$. We use the notation $a_+ = a \vee 0$ to denote the nonnegative part of a real number $a$, and $\mathrm{sgn}(a)$ to denote its sign.

\subsection{MMA with general candidate models}\label{sec:MMA_general}

A widely used approach for estimating the weight vector is to minimize the Mallows-type MA criterion:
\begin{equation}\label{eq:c}
C_{n}(\bw|\mathcal{M},\lambda)\triangleq n^{-1}\left\|\by-\bP(\bw)\by\right\|^{2}+2\lambda^2\hat{\sigma}^{2}\tr\bP(\bw),
\end{equation}
where $\hat{\sigma}^{2}$ is an estimator for $\sigma^2$, and $\lambda$ is a penalty parameter. When $\lambda$ is set as $\lambda_1 \triangleq \sqrt{1/n}$ and $\mathbb{E}[\hat{\sigma}^{2}] = \sigma^2$, it can be verified that for any non-random $\bw$, $\mathbb{E}C_{n}(\bw|\mathcal{M},\lambda_1) =  R_n(\hat{\bmu}_{\bw|\mathcal{M}}, \bmu) + \sigma^2$. Thus, $C_n(\bw|\mathcal{M},\lambda_1)$ provides an unbiased estimate of the MA risk up to a constant \citep{Hansen2007least}. As $\lambda$ increases, the criterion $C_n(\bw|\mathcal{M},\lambda)$ imposes a larger penalty on weights associated with larger models, which may lead to a more parsimonious MA estimator; see \cite{Zhang2020Parsimonious} for the choice of $\lambda=\sqrt{(\log n)/n}$. 

The first goal of this paper is to investigate the Mallows-type criterion (\ref{eq:c}) under the classical choice $\lambda_1 \triangleq \sqrt{1/n}$. The estimated weight vector obtained by the classical MMA is then given by $\hat{\bw}_1\triangleq\argmin _{\bw \in \mathcal{W}} C_{n}(\bw|\mathcal{M}, \lambda_1)$. The resulting MMA estimator is
\begin{equation}\label{eq:mahat}
\hat{\bmu}_{\hat{\bw}_1|\mathcal{M}}=\sum_{m=1}^{M_{n}} \hat{w}_{1m} \hat{\bmu}_m,
\end{equation}
where $\hat{w}_{1m}$ is the $m$-th element of $\hat{\bw}_1$.

When no additional prior restrictions are imposed on the candidate models $\mathbb{V}_1,\ldots,\mathbb{V}_{M_n}$, the works of \cite{WAN2010277} and \cite{Zhang_2021} have deeply studied the asymptotic performance of (\ref{eq:mahat}) under Assumption~\ref{ass:4_moment}. It was proved by \cite{WAN2010277} that, if  
\begin{equation}\label{eq:wan_zhang}
  \frac{ [ M_n \sum_{m=1}^{M_n} (\| (\bI - \bP_m)\bmu \|^2 + \sigma^2 k_m) ]^{1/2} }{n R_{n}(\hat{\bmu}_{\bw^*|\mathcal{M}}, \bmu)} \to 0,
\end{equation}
or by \cite{Zhang_2021} that 
\begin{equation}\label{eq:wan_zhang_2}
  \frac{  M_n^2 }{n R_{n}(\hat{\bmu}_{\bw^*|\mathcal{M}}, \bmu)} \to 0,
\end{equation}
then 
\begin{equation}\label{eq:aop_wan_zhang}
  \frac{L_n(\hat{\bmu}_{\hat{\bw}_1|\mathcal{M}}, \bmu )}{\min_{\bw \in \mathcal{W}}L_n(\hat{\bmu}_{\bw|\mathcal{M}}, \bmu )} \to 1
\end{equation}
in probability. To the best of our knowledge, conditions (\ref{eq:wan_zhang})--(\ref{eq:wan_zhang_2}) collectively give the mildest known condition under which the MMA estimator can achieve (\ref{eq:aop_wan_zhang}) under Assumption~\ref{ass:4_moment} and for general candidate model set $\mathcal{M}$.

The asymptotic result in (\ref{eq:aop_wan_zhang}) focuses on the large-sample limit as $n \to \infty$. In this paper, we investigate the finite-sample risk behavior of the MMA estimator. Let $\mathbb{U}$ denote a subspace of $\mathbb{R}^n$ of interest (e.g., $\mathbb{R}^n$ or the $\ell^2$ bounded set $\mathbb{B}_2^L \triangleq \{\bmu: \|\bmu \|^2/n \leq L  \}$), and $\mathbf{M}(M_n) \triangleq \{ \mathcal{M}: \mathrm{Card}(\mathcal{M}) = M_n \}$ represents the collection of candidate model sets containing $M_n$ models. In this paper, our first goal is to answer the following question:
\begin{description}
  \item[Q1.] How can we construct a finite-sample upper bound on $\mathbb{E}L_n(\hat{\bmu}_{\hat{\bw}_1|\mathcal{M}}, \bmu ) - R_{n}(\hat{\bmu}_{\bw^*|\mathcal{M}}, \bmu)$ that holds for all $\bmu \in \mathbb{U}$ and $\mathcal{M} \in \mathbf{M}(M_n)$ under Assumption~\ref{ass:4_moment}?
\end{description}
The answer to Q1 can provide a finite-sample performance guarantee for the classical MMA estimator (\ref{eq:mahat}) over the general class of candidate model sets $\mathbf{M}(M_n)$.

\subsection{Construction of candidate models}\label{subsec:candidate}

Another critical factor that affects the performance of the MA estimator \eqref{eq:ma} is the choice of the candidate model set $\mathcal{M}$. In general, the subspaces in $\mathcal{M}$ may have arbitrary relationships. To facilitate theoretical analysis, we consider a structured setting in which all subspaces are spanned by vectors from a given \emph{complete orthogonal basis} $\{\bpsi_1, \ldots, \bpsi_p\}$, as specified in the following assumption.
\begin{assumption}\label{ass:complete}
  There exists a complete orthogonal basis $\{\bpsi_1, \ldots, \bpsi_p\}$ such that $\bpsi_j \in \mathbb{R}^n$, $n^{-1}\|\bpsi_j\|^2 = 1$, and $\bpsi_j^\top \bpsi_{j'} = 0$ for $j \neq j'$. Furthermore, the ture regression mean vector $\bmu$ has the representation
  \begin{equation}\label{eq:complete}
  \bmu = \sum_{j=1}^{p} \theta_j \bpsi_j,
  \end{equation}
  where $1 \leq p \leq n$, and $\theta_j = \bpsi_j^\top \bmu / n$.
\end{assumption}

A complete orthogonal basis satisfying (\ref{eq:complete}) with $p=n$ always exists, given that $\bmu \in \mathbb{R}^n$.  In practice, commonly used transformations such as the discrete cosine transform \citep[see, e.g.,][]{Rao1990} and the discrete wavelet transform \citep[see, e.g.,][]{Daubechies1988} can be adopted to construct $\{\bpsi_1, \ldots, \bpsi_n\}$. In the linear regression setting where $\bmu = \bX \bbeta$ with $\bX \in \mathbb{R}^{n \times d}$ and $\bbeta \in \mathbb{R}^d$, a complete basis with $p \leq \min(n, d)$ can be constructed using the singular value decomposition (SVD) of $\bX$ \citep[see, e.g.,][]{Jeffers1967, Zhu2023Scalable} . The theory and methods developed in this paper apply to any given complete orthogonal basis that satisfies condition (\ref{eq:complete}). In the numerical experiments in Section~\ref{sec:compare}, we discuss the use of SVD to construct the basis.

Given an index set $\mathcal{I} \subseteq \{1,\ldots,p \}$, let $\bPsi_{\mathcal{I}} \in \mathbb{R}^{n \times |\mathcal{I}|}$ denote the regressor matrix whose $j$-th column corresponds to $\bpsi_j$ for $j \in \mathcal{I}$. The estimator of $\bmu$ based on model ${\mathcal{I}}$ is then given by
\begin{equation}\label{eq:estimator}
    \hat{\bmu}_{\mathcal{I}} = \bP_{\mathcal{I}}\by \triangleq \bPsi_{\mathcal{I}}(\bPsi_{\mathcal{I}}^{\top}\bPsi_{\mathcal{I}})^{-1}\bPsi_{\mathcal{I}}^{\top}\by.
\end{equation}
In this paper, we consider two representative methods for constructing candidate model set $\mathcal{M} = \{ \mathcal{I}_1,\ldots, \mathcal{I}_{M_n} \}$.

The first approach considers nested candidate models \citep[see, e.g.,][]{Shibata1980, Breiman01031983, Li1987, Hansen2007least}. Specifically, we define the candidate model set as $\mathcal{M}_{AN} \triangleq \left\{ \{1 \}, \{1,2 \},\ldots,\{1,2,\ldots,p \}  \right\}$. Let $R_{n}(\hat{\bmu}_{\bw^*|\mathcal{M}_{AN}}, \bmu)$ denote the optimal MA risk based on all nested candidate models in $\mathcal{M}_{AN}$. The achievability of $R_{n}(\hat{\bmu}_{\bw^*|\mathcal{M}_{AN}}, \bmu)$ has been studied in \cite{Peng2024Optimality} under the assumption that $\epsilon_i$ follows a sub-Gaussian distribution. This raises the following question:
\begin{description}
  \item[Q2.] Can the optimal all-nested MA risk $R_{n}(\hat{\bmu}_{\bw^*|\mathcal{M}_{AN}}, \bmu)$ still be attainable under Assumption~\ref{ass:4_moment}?
\end{description}

The successful application of nested MA relies on the assumption that $|\theta_j|$ are ordered in descending magnitude. Define the candidate model set with all-subset models $\mathcal{M}_{AS}\triangleq \{ \mathcal{I} :  \mathcal{I} \subseteq \{1,\ldots,p \} \}$, and define the ideal MA risk based on all-subset models as $R_{n}(\hat{\bmu}_{\bw^*|\mathcal{M}_{AS}}, \bmu) \triangleq \min_{\bw }R_{n}(\hat{\bmu}_{\bw|\mathcal{M}_{AS}}, \bmu)$. Section~5 of \cite{Peng2024Optimality} shows that when $|\theta_j|$ are ordered, we have $R_{n}(\hat{\bmu}_{\bw^*|\mathcal{M}_{AN}}, \bmu) \sim R_{n}(\hat{\bmu}_{\bw^*|\mathcal{M}_{AS}}, \bmu)$, and the nested MMA estimator can attain $R_{n}(\hat{\bmu}_{\bw^*|\mathcal{M}_{AS}}, \bmu)$ asymptotically. However, if the ordering assumption is violated, the optimal all-nested MA risk $R_{n}(\hat{\bmu}_{\bw^*|\mathcal{M}_{AN}}, \bmu)$ may suffer a loss in efficiency \citep[see Section~3.3 of][]{Peng2023shrinkage}. In this setup, how to construct an estimator that approaches $R_{n}(\hat{\bmu}_{\bw^*|\mathcal{M}_{AS}}, \bmu)$ as closely as possible remains unknown.
\begin{description}
  \item[Q3.] What is the fundamental limit of achieving $R_{n}(\hat{\bmu}_{\bw^*|\mathcal{M}_{AS}}, \bmu)$ under the general assumption $\bmu \in \mathbb{R}^n$? Moreover, how can we construct an estimator to attain this limit?
\end{description}

Note that Q1 investigates the risk performance of the classical MMA estimator without imposing restrictions on the candidate models. Q2 and Q3 focus on constructing specific MA estimators with explicit consideration of candidate model construction. Addressing these questions will enhance both the theoretical understanding and practical application of MA.

\section{General candidate models}\label{sec:general}

\subsection{Oracle inequalities}\label{subsec:oracle}

In this subsection, we establish several oracle inequalities for the classical MMA estimator (\ref{eq:mahat}) based on general candidate model set $\mathcal{M}$.

\begin{proposition}\label{theo:sharp}
Suppose Assumption~\ref{ass:4_moment} holds. Then, for any candidate model set $\mathcal{M} \in \mathbf{M}(M_n)$ and any $\bmu \in \mathbb{R}^n$, there exists a constant $C>0$ such that
  \begin{equation}\label{eq:sharp}
  \begin{split}
     \mathbb{E}L_n(\hat{\bmu}_{\hat{\bw}_1|\mathcal{M}}, \bmu ) & \leq R_{n}(\hat{\bmu}_{\bw^*|\mathcal{M}}, \bmu) +  Cn^{-1}\left(\sum_{m=1}^{M_n}\| (\bI - \bP_m)\bmu \|^2\right)^{1/2} + C n^{-1}\left(\sum_{m=1}^{M_n} k_m\right)^{1/2} \\
       & \qquad\qquad\qquad\quad\;\,+ C n^{-1}\mathbb{E}\left|\hat{\sigma}^2 - \sigma^2\right|\max_{1 \leq m \leq M_n}k_m,
  \end{split}
\end{equation}
where $\hat{\bmu}_{\hat{\bw}_1|\mathcal{M}}$ is the MMA estimator defined in (\ref{eq:mahat}).
\end{proposition}

The inequality (\ref{eq:sharp}) is referred to as a sharp oracle inequality for the MA estimator, where the leading constant in the optimal MA risk term is exactly one \citep[see, e.g.,][]{Dalalyan2012}. Such an inequality provides a uniform performance guarantee for the classical MMA estimator $\hat{\bmu}_{\hat{\bw}_1|\mathcal{M}}$ relative to the optimal MA risk $R_n(\hat{\bmu}_{\bw^*|\mathcal{M}},\bmu)$ over all $\mathcal{M}\in\mathbf{M}(M_n)$ and $\bmu\in\mathbb{R}^n$, regardless of how fast or slowly $R_n(\hat{\bmu}_{\bw^*|\mathcal{M}},\bmu)$ converges. For the examination of the order of $R_n(\hat{\bmu}_{\bw^*|\mathcal{M}},\bmu)$, see Section~\ref{sec:nested} where a specific $\mathcal{M}$ is considered. Suppose that $\mathbb{E}|\hat{\sigma}^2 - \sigma^2| = O(\sqrt{M_n/n})$ and $\bmu \in \mathbb{B}_2^L$. Then, moving the optimal MA risk $R_{n}(\hat{\bmu}_{\bw^*|\mathcal{M}}, \bmu)$ to the left side of (\ref{eq:sharp}) and taking the maximum over $\bmu \in \mathbb{B}_2^L$ and $\mathcal{M} \in \mathbf{M}(M_n)$ yields the worst-case risk bound: 
\begin{equation}\label{eq:chuntian2}
  \max_{\bmu \in \mathbb{B}_2^L,\mathcal{M} \in \mathbf{M}(M_n)}\Bigl\{\mathbb{E}L_n(\hat{\bmu}_{\hat{\bw}_1|\mathcal{M}}, \bmu ) - R_{n}(\hat{\bmu}_{\bw^*|\mathcal{M}}, \bmu)\Bigl\} \leq C\left( \frac{M_n}{n} \right)^{1/2}.
\end{equation}
The upper bound in (\ref{eq:chuntian2}) provides a worst-case performance guarantee for the MMA estimator across a general class of candidate model sets. It means that for any $\bmu \in \mathbb{B}_2^L$ and $\mathcal{M} \in \mathbf{M}(M_n)$, $\mathbb{E}L_n(\hat{\bmu}_{\hat{\bw}_1|\mathcal{M}}, \bmu )$ cannot exceed $R_{n}(\hat{\bmu}_{\bw^*|\mathcal{M}}, \bmu)$ by more than a term of order $( M_n/n )^{1/2}$. It is also observed that the left-hand side of (\ref{eq:chuntian2}) must be positive when $M_n$ is larger than a sufficiently large constant. This is implied by reducing the general $\bmu$ and $\mathbf{M}(M_n)$ to the specific setting established in Theorem~\ref{theo:lower}, as shown later. However, even with a fixed number of candidate models, the right-hand side of (\ref{eq:chuntian2}) converges no faster than $n^{-1/2}$.

To achieve faster converging rate when $M_n$ is small, we combine the shifted empirical process technique \citep[see, e.g.,][]{Baraud2000, Wegkamp2003, cao2005oracle} with the results in \cite{Zhang_2021} to derive the following (non-exact) oracle inequality.

\begin{theorem}\label{tho:1}
Suppose that Assumption~\ref{ass:4_moment} holds. For an arbitrary quantity $0 < \delta < \infty$ that can depend on $n$, the risk of the MMA estimator (\ref{eq:mahat}) is upper bounded by
\begin{equation}\label{eq:riskbound1}
\begin{split}
   \mathbb{E}L_n(\hat{\bmu}_{\hat{\bw}_1|\mathcal{M}}, \bmu )&\leq (1+C\delta) R_{n}(\hat{\bmu}_{\bw^*|\mathcal{M}}, \bmu) + \frac{C(1+\delta)^3M_n}{\delta n}  + \frac{C(1+\delta)M_n^2}{n}+  \frac{C(1+\delta)^2\mathbb{E}(\hat{\sigma}^2 - \sigma^2)^2}{\delta},
\end{split}
\end{equation}
where $C$ is a positive constant independent of $n$ and $\delta$.
\end{theorem}

Comparing the sharp oracle inequality (\ref{eq:sharp}) with (\ref{eq:riskbound1}), we observe that the leading constant in (\ref{eq:riskbound1}) is greater than one. Suppose that $\mathbb{E}(\hat{\sigma}^2 - \sigma^2)^2 = O(M_n/n)$. Due to the arbitrariness of $\delta$, by choosing $\delta = \delta_n$ such that $\delta_n \to 0$ and $(1+\delta_n)^3/\delta_n = O(M_n)$ if $M_n \to \infty$, we obtain the following worst-case bound
\begin{equation}\label{eq:riskbound2}
  \max_{\bmu \in \mathbb{R}^n,\mathcal{M} \in \mathbf{M}(M_n)}\Bigl\{\mathbb{E}L_n(\hat{\bmu}_{\hat{\bw}_1|\mathcal{M}}, \bmu ) - [1+o(1)] R_{n}(\hat{\bmu}_{\bw^*|\mathcal{M}}, \bmu)\Bigl\} \leq \frac{C M_n^2}{n}.
\end{equation}
Note that (\ref{eq:riskbound2}) holds over the broader parameter space $\mathbb{R}^n$ than $\mathbb{B}_2^L$ in (\ref{eq:chuntian2}). By ``absorbing'' some higher-order terms into $[1+o(1)]R_{n}(\hat{\bmu}_{\bw^*|\mathcal{M}}, \bmu)$ using the shifted empirical process technique, (\ref{eq:riskbound2}) guarantees a faster worst-case convergence rate compared to (\ref{eq:chuntian2}) when $M_n \lesssim n^{1/3}$.

\begin{remark}
  To the best of our knowledge, the non-exact oracle inequality for the classical MMA with variance estimation presented in Theorem~\ref{tho:1} has not been established in the existing literature. There are two related works. The first is by \cite{Bellec2018}, where affine estimators are considered as candidates. When $\sigma^2$ is assumed to be known and $\epsilon_i$ follows a Gaussian distribution, Proposition 7.2 in \cite{Bellec2018} implies that
\begin{equation}\label{eq:A27}
  \max_{\bmu \in \mathbb{B}_2^L,\mathcal{M} \in \mathbf{M}(M_n)}\Bigl\{\mathbb{E}L_n(\hat{\bmu}_{\hat{\bw}_1|\mathcal{M}}, \bmu ) - R_{n}(\hat{\bmu}_{\bw^*|\mathcal{M}}, \bmu)\Bigl\} \leq C\left( \frac{\log M_n}{n} \right)^{1/2}.
\end{equation}
However, this bound still cannot guarantee a fast convergence rate when a small number of candidate models are combined. The second is \cite{Zhang_2021}, which establishes the asymptotic optimality of (\ref{eq:aop_wan_zhang}) under the condition $M_n^2/(n R_{n}(\hat{\bmu}_{\bw^*|\mathcal{M}}, \bmu)) \to 0$. 
\end{remark}

\subsection{Implications for AOP with general candidates}

Based on the oracle inequalities established in Proposition~\ref{theo:sharp} and Theorem~\ref{tho:1}, the AOP of the MMA estimator (\ref{eq:mahat}) is obtained; see Section~\ref{sec:proof_add_nested} of the Appendix for the proof.
\begin{corollary}\label{cor:aop}
  Suppose Assumption~\ref{ass:4_moment} holds. For any $\bmu \in \mathbb{R}^n$, if the candidate model set  $\mathcal{M}$ and the variance estimator $\hat{\sigma}^2$ satisfy the one of the following conditions:
  \begin{equation}\label{eq:condition_our}
  \frac{ [ \sum_{m=1}^{M_n} (\| (\bI - \bP_m)\bmu \|^2 + \sigma^2 k_m) ]^{1/2} }{n R_{n}(\hat{\bmu}_{\bw^*|\mathcal{M}}, \bmu)} \to 0 \quad \text{and}\quad \frac{\mathbb{E}|\hat{\sigma}^2 - \sigma^2|\max_{1 \leq m \leq M_n}k_m }{n R_{n}(\hat{\bmu}_{\bw^*|\mathcal{M}}, \bmu)} \to 0, 
\end{equation}
or 
\begin{equation}\label{eq:variance_cond}
  \frac{M_n^2 }{n R_{n}(\hat{\bmu}_{\bw^*|\mathcal{M}}, \bmu)} \to 0 \quad \text{and}\quad \frac{n \mathbb{E}(\hat{\sigma}^2 - \sigma^2)^2}{M_n} \to 0,
\end{equation}
then the MMA estimator achieves the AOP:
\begin{equation}\label{eq:aop_our}
  \frac{\mathbb{E}L_n(\hat{\bmu}_{\hat{\bw}_1|\mathcal{M}}, \bmu )}{R_{n}(\hat{\bmu}_{\bw^*|\mathcal{M}}, \bmu)} \leq 1 + o(1),\quad n \to \infty.
\end{equation}
\end{corollary}

The first equations in (\ref{eq:condition_our})--(\ref{eq:variance_cond}) are the key requirements for regulating the candidate model set to achieve AOP in MA risk. Comparing the first equation in (\ref{eq:condition_our}) with (\ref{eq:wan_zhang}), we observe that the numerator in (\ref{eq:condition_our}) removes an $M_n$ factor compared to (\ref{eq:wan_zhang}). Thus, Corollary~\ref{cor:aop} suggests that achieving AOP in terms of risk imposes milder conditions than those required for loss. Condition~(\ref{eq:variance_cond}) is same as (\ref{eq:wan_zhang_2}) but additionally incorporates the variance estimation. The variance estimation requirements in (\ref{eq:condition_our})--(\ref{eq:variance_cond}) are satisfied in several scenarios: (\romannumeral1) when $\hat{\sigma}^2$ is assumed to be known \citep[see. e.g.,][]{Bellec2018, Zhang_2021}, (\romannumeral2) when $\mathbb{E}(\hat{\sigma}^2 - \sigma^2)^2 = O(1/n)$ \citep[see, e.g., Section~4.2 of][]{Peng2024Optimality} and $M_n \to \infty$, which ensure (\ref{eq:variance_cond}), and (\romannumeral3) when using the estimator in Theorem 2 of \cite{WAN2010277} under some additional conditions on $\max_{1 \leq m \leq M_n}k_m$ and $R_{n}(\hat{\bmu}_{\bw^*|\mathcal{M}}, \bmu)$ that ensure (\ref{eq:condition_our}). For simplicity, we assume that $\sigma^2$ is known. Table~\ref{tab:addlabel} summarizes existing results on the AOP of MMA under Assumption~\ref{ass:4_moment} and a general candidate model set $\mathcal{M}$.

\begin{table}[htbp]
  \centering
  \caption{Sufficient conditions on $\mathcal{M}$ for achieving AOP under Assumption~\ref{ass:4_moment}.}
  \begin{tabular}{lcc|cc}
    \hline
    \multirow{2}{*}{Article} & \multirow{2}{*}{$\mathcal{M}$ Condition} & \multirow{2}{*}{$\bmu \in$} & \multicolumn{2}{c}{Asymptotic Optimality in} \\
    \cline{4-5}
    & & & Loss & Risk \\
    \hline
    \cite{WAN2010277} & $\frac{[ M_n \sum_{m=1}^{M_n} (\| (\bI - \bP_m)\bmu \|^2 + \sigma^2 k_m) ]^{1/2}  }{n R_{n}(\hat{\bmu}_{\bw^*|\mathcal{M}}, \bmu)} \to 0$ & $\mathbb{R}^n$ & \checkmark &  \\
    \cite{Zhang_2021} & $\frac{M_n^2 }{n R_{n}(\hat{\bmu}_{\bw^*|\mathcal{M}}, \bmu)} \to 0$ & $\mathbb{R}^n$ & \checkmark &  \\
    \hline
    \multirow{2}{*}{This paper} & $\frac{ [ \sum_{m=1}^{M_n} (\| (\bI - \bP_m)\bmu \|^2 + \sigma^2 k_m) ]^{1/2} \wedge M_n^2 }{n R_{n}(\hat{\bmu}_{\bw^*|\mathcal{M}}, \bmu)} \to 0$ & $\mathbb{R}^n$ &  & \checkmark \\
    & $\frac{ (nM_n)^{1/2} \wedge M_n^2 }{n R_{n}(\hat{\bmu}_{\bw^*|\mathcal{M}}, \bmu)} \to 0$ & $\mathbb{B}_2^L$ &  & \checkmark \\
    \hline
  \end{tabular}
  \label{tab:addlabel}
\end{table}

\begin{remark}
  Both the loss and risk versions of AOP are widely adopted in the literature \citep[see, e.g.,][]{Zhang2020Parsimonious, Peng2024Optimality, yu2024unified}. They have been established simultaneously under the comparable conditions; see, e.g., Theorem~3 of \cite{Zhang2020Parsimonious} and Theorem~1 and Corollary~A.1 of \cite{Peng2024Optimality}. It is worth noting that a recent study by \cite{xu2024asymptotic} reveals that a fundamental difference may exist between (\ref{eq:aop_wan_zhang}) and (\ref{eq:aop_our}) when the true model is included in $\mathcal{M}$. %In general setting, whether an intrinsic difference exists between (\ref{eq:aop_wan_zhang}) and (\ref{eq:aop_our}) remains unknown.
\end{remark}

%The oracle inequalities and asymptotic results in the this section provide a uniform performance guarantee for the MMA estimator across general candidate model constructions. However, they offer limited insight into how candidate model construction influences the resulting MA estimator. In the following sections, we investigate how to construct candidate models to obtain a well-performing MA estimator.

\subsection{Implications for all-nested MA}\label{sec:nested}

This subsection demonstrates that the oracle inequalities in Section~\ref{subsec:oracle} are important tools to answer the all-nested MA problem posed in Question 2. The nested MA plays a key role toward achieving the optimal all-subset MA risk when the regression coefficients are ordered \citep[see Section~5 of][]{Peng2024Optimality}. This problem has been extensively studied in \cite{Peng2024Optimality} and \cite{Peng2023shrinkage} under sub-Gaussian and Gaussian assumptions on the random error term, respectively. We show that the optimal all-nested MA risk $R_{n}(\hat{\bmu}_{\bw^*|\mathcal{M}_{AN}}, \bmu)$ remains attainable under the weaker Assumption~\ref{ass:4_moment}.

%In this section, we aim to anwer this quesion Q2
%If an specific candidate model set $\mathcal{M}$ is considered such that the optimal MA risk $R_{n}(\hat{\bmu}_{\bw^*|\mathcal{M}}, \bmu)$ converges slowly, although such an optimal MA risk can be attained, the such an MA estimator is dount.

%To address this question, in the rest of this paper, we discuss the candidate model construction problem of least squares MA raised in Section~\ref{}

The approach is to construct nested candidate models based on a system of weakly geometrically increasing blocks \citep{cavalier2001penalized} and then apply the general MMA bound from Theorem~\ref{tho:1}. Define $\rho_n= 1/\log p$, $j_1=\lceil \log p \rceil$, $j_t=j_{t-1}+\lfloor j_1(1+\rho_n)^{t-1}  \rfloor $ for $t=2,\ldots,T_n-1$, and $j_{T_n}=p$, where $
T_n\triangleq \argmin_{m \in \mathbb{N}} \{(j_1 + \sum_{t=2}^{m}\lfloor j_1(1+\rho_{n})^{t-1}\rfloor) \geq p\}$. We then construct the group-wise candidate model set
$$
\mathcal{M}_G \triangleq \{ \{1,\ldots,j_1 \}, \{1,\ldots,j_2 \},\ldots, \{1,\ldots,j_{T_n} \} \}.
$$
Let $\hat{\bmu}_{\hat{\bw}_1|\mathcal{M}_G}$ denote the MMA estimator (\ref{eq:mahat}) with $\mathcal{M} = \mathcal{M}_G$.
\begin{corollary}\label{tho:3}
  Suppose Assumption~\ref{ass:4_moment} holds, and $p \to \infty$ as $n \to \infty$. The nested MMA estimator $\hat{\bmu}_{\hat{\bw}_1|\mathcal{M}_G}$ satisfies the following upper bound for any $\bmu \in \mathbb{R}^n$:
  \begin{equation}\label{eq:tho:3}
\begin{split}
    \mathbb{E}L_n(\hat{\bmu}_{\hat{\bw}_1|\mathcal{M}_G}, \bmu )&\leq [1+o(1)]R_{n}(\hat{\bmu}_{\bw^*|\mathcal{M}_{AN}}, \bmu)+\frac{C(\log p)^4}{n} + C(\log p)^2 \mathbb{E}(\hat{\sigma}^2 - \sigma^2)^2, \\
\end{split}
\end{equation}
where $C>0$ is a constant independent of $n$.
\end{corollary}

Corollary~\ref{tho:3} establishes the achievability of the optimal MA risk for all nested candidate models. Consider the representative case where $p = n$, and suppose that $\mathbb{E}(\hat{\sigma}^2 - \sigma^2)^2 =O((\log n)^2/n)$. Dividing both sides of (\ref{eq:tho:3}) in Corollary~\ref{tho:3} by $R_{n}(\hat{\bmu}_{\bw^*|\mathcal{M}_{AN}}, \bmu)$, we observe that if 
\begin{equation}\label{eq:cond:nested}
  \frac{(\log n)^4}{nR_{n}(\hat{\bmu}_{\bw^*|\mathcal{M}_{AN}}, \bmu)} \to 0,
\end{equation}
then
\begin{equation*}
  \frac{\mathbb{E}L_n(\hat{\bmu}_{\hat{\bw}_1|\mathcal{M}_G}, \bmu )}{R_{n}(\hat{\bmu}_{\bw^*|\mathcal{M}_{AN}}, \bmu)} \leq 1+ o(1).
\end{equation*}
This result suggests that as long as $R_{n}(\hat{\bmu}_{\bw^*|\mathcal{M}_{AN}}, \bmu)$ does not converge too fast, the full potential of nested MA remains attainable under Assumption~\ref{ass:4_moment}. Condition~(\ref{eq:cond:nested}) is comparable to those imposed under the sub-Gaussian setting \citep[Theorem~3 of][]{Peng2024Optimality}, differing only in a logarithmic term in the numerator. Therefore, the advantage of model averaging over model selection in the nested model setting established in \cite{Peng2022improvability} can still be realized by MMA under the more general Assumption~\ref{ass:4_moment}. There are several representative cases under which (\ref{eq:cond:nested}) can be satisfied. The first case is when $|\theta_j|$ in (\ref{eq:complete}) decays at a polynomial rate, i.e., $|\theta_j| = j^{-\alpha_1} $, $\alpha_1 > 1/2$. In this case, we have $R_{n}(\hat{\bmu}_{\bw^*|\mathcal{M}_{AN}}, \bmu) \asymp n^{-1+1 /\left(2 \alpha_1\right)}$ \citep{Peng2022improvability, Peng2024Optimality}, and hence (\ref{eq:cond:nested}) is satisfied for any $\alpha_1 > 1/2$. The second case is when $|\theta_j|$ decays at an exponential rate, i.e., $|\theta_j| = \exp (- j^{\alpha_2})$, $\alpha_2>0$. In this case, $R_{n}(\hat{\bmu}_{\bw^*|\mathcal{M}_{AN}}, \bmu) \asymp(\log n)^{1 / \alpha_2}/n$ \citep{Peng2022improvability, Peng2024Optimality}, and (\ref{eq:cond:nested}) holds when $0<\alpha_2<1/4$. Beyond these cases, if there exists potential model misspecification such that $R_{n}(\hat{\bmu}_{\bw^*|\mathcal{M}_{AN}}, \bmu)$ converges at a very slow rate, then (\ref{eq:cond:nested}) may also be satisfied.

%The analysis in this section broadens the understanding of the optimal nested MA beyond the sub-Gaussian assumptions studied in \cite{Peng2024Optimality} and \cite{Peng2023shrinkage}. See Table~2 for a detailed comparison.
%\begin{table}[htbp]
%  \centering
%  \caption{Conditions on the optimal MA risk and random error terms for achieving the optimal nested MA risk with $p=n$.}
%  \renewcommand{\arraystretch}{1.3} % Increases row spacing for better readability
%  \begin{tabular}{lcc}
%    \hline
%    Article & Optimal MA risk  & Error distribution \\
%    \hline
%    \cite{Peng2024Optimality}  & $\frac{(\log n)^2}{nR_{n}(\hat{\bmu}_{\bw^*|\mathcal{M}_{AN}}, \bmu)} \to 0$ & Sub-Gaussian \\
%    \cite{Peng2023shrinkage}  & $\frac{\log \log n}{nR_{n}(\hat{\bmu}_{\bw^*|\mathcal{M}_{AN}}, \bmu)} \to 0$ & Gaussian \\
%    This paper & $\frac{(\log n)^4}{nR_{n}(\hat{\bmu}_{\bw^*|\mathcal{M}_{AN}}, \bmu)} \to 0$ & Assumption~\ref{ass:4_moment} \\
%    \hline
%  \end{tabular}
%  \label{tab:conditions_optimal_MA}
%\end{table}

\section{All-subset candidate models}\label{sec:all-subset}

In this section, we study the all-subset MA problem under the orthogonal basis that satisfies Assumption~\ref{ass:complete}. Following the classical AOP theory, we assess the performance of an estimator $\hat{\bmu}$ by the risk ratio $R_n( \hat{\bmu},  \bmu )/R_n(\hat{\bmu}_{\bw^*|\mathcal{M}_{AS}}, \bmu )$, which quantifies its risk relative to the optimal all-subset MA risk at $\bmu$. Recall that $R_n(\hat{\bmu}_{\bw^*|\mathcal{M}_{AS}}, \bmu ) = \min_{\bw }R_{n}(\hat{\bmu}_{\bw|\mathcal{M}_{AS}}, \bmu)$ denotes the optimal risk of linearly combining $2^p$ least squares estimators corresponding to the models in $\mathcal{M}_{AS}= \{ \mathcal{I} :  \mathcal{I} \subseteq \{1,\ldots,p \} \}$, thus representing the ideal performance of least squares MA.

\subsection{Fundamental limit}\label{subsec:fundamental}

In this subsection, we establish two minimax lower bounds for the risk ratio. Since the minimax lower bound is on the negative side (limit of achieving the optimal all-subset MA risk), we assume that the random errors follow a Gaussian distribution. When a more general error distribution class is considered, such as that in Assumption~\ref{ass:4_moment}, the problem of achieving $R_n(\hat{\bmu}_{\bw^*|\mathcal{M}_{AS}}, \bmu )$ certainly can not be easier.

Define the \emph{hardest cube} as
\begin{equation}\label{eq:hardest_cube}
  \Theta^* \triangleq \biggl\{ \btheta \in \mathbb{R}^p: 0 \leq |\theta_j| \leq \sqrt{\frac{2\sigma^2 \log p}{n}} \biggl\}.
\end{equation}
For any parameter space $\Theta \subseteq \mathbb{R}^p$, let $\mathcal{C}(\Theta) \triangleq \{\bmu = \sum_{j=1}^{p}\theta_j\bpsi_j: \btheta \in \Theta  \}$ denote the associated class of regression mean vectors. We have the following minimax lower bounds.

\begin{theorem}\label{theo:lower}
  Suppose that $\epsilon_1,\ldots,\epsilon_n$ are i.i.d. $N(0,\sigma^2)$. For any parameter space $\Theta$ satisfying $\Theta^*\subseteq\Theta$, there exists a sufficiently large integer $\underline{p}>0$ independent of $n$, such that if $p \geq \underline{p}$, then
    \begin{equation}\label{eq:lower_fixed}
    \min_{\hat{\bmu}}\max_{\bmu \in \mathcal{C}(\Theta) } \frac{R_n\left( \hat{\bmu},  \bmu \right)}{n^{-1}+R_n\left(\hat{\bmu}_{\bw^*|\mathcal{M}_{AS}}, \bmu \right)} > \underline{c}, 
  \end{equation}
  where $\underline{c}>1$ is a constant depending on $\underline{p}$. 
  If $p \to \infty$ as $n \to \infty$, then
  \begin{equation}\label{eq:lower_1}
    \min_{\hat{\bmu}}\max_{\bmu \in \mathcal{C}(\Theta) } \frac{R_n\left( \hat{\bmu},  \bmu \right)}{n^{-1}+R_n\left(\hat{\bmu}_{\bw^*|\mathcal{M}_{AS}}, \bmu \right)} \geq \left[ 1 + o(1) \right]2\log p,
  \end{equation}
  where the minimum is taken over all measurable estimators $\hat{\bmu}$ based on $\by$.

\end{theorem}

Theorem~\ref{theo:lower} suggests that there exist fundamental limits of achieving the optimal all-subset MA risk. For any parameter space $\Theta$ contains $\Theta^*$ (e.g., the whole space $\mathbb{R}^p$), even in the parametric case where there exists a fixed dimensional true model, the maximum risk ratio over $\Theta$ is strictly larger than a constant $\underline{c}>1$ for any estimator. Furthermore, in the diverging dimension scenario where $p \to \infty$, the minimax risk ratio diverges to $\infty$ at the asymptotic rate $2\log p$. Note that the minimax lower bounds in (\ref{eq:lower_fixed})--(\ref{eq:lower_1}) also serve as the minimax lower bounds for the risk ratio $R_n( \hat{\bmu},  \bmu )/R_n(\hat{\bmu}_{\bw^*|\mathcal{M}_{AS}}, \bmu )$ and for $nR_n( \hat{\bmu},  \bmu )$, which also imply that $\min_{\hat{\bmu}}\max_{\bmu \in \mathcal{C}(\Theta) }[R_n( \hat{\bmu},  \bmu ) - R_n(\hat{\bmu}_{\bw^*|\mathcal{M}_{AS}}, \bmu )] \geq 0$ and $\min_{\hat{\bmu}}\max_{\bmu \in \mathcal{C}(\Theta) }R_n( \hat{\bmu},  \bmu ) \gtrsim (\log p)/n $. 

The minimax lower bounds established in Theorem~\ref{theo:lower} have several important implications. First, they broaden the scope of the classical AOP theory, which justifies the optimality of MA by demonstrating that the risk ratio approaches one asymptotically \citep[see, e.g.,][]{Hansen2007least, WAN2010277}. Our results show that even in the setting where $p$ is fixed, achieving $R_n( \hat{\bmu},  \bmu )/R_n(\hat{\bmu}_{\bw^*|\mathcal{M}_{AS}}, \bmu ) \to 1$ is theoretically impossible for any estimators unless the parameter space is restricted to a more structured subset than $\Theta^*$; see, for example, the weakly ordered space in Theorem 5 of \cite{Peng2024Optimality}. Second, these lower bounds serve as fundamental benchmarks for the best achievable convergence rate of any estimator relative to the optimal all-subset MA risk $ R_n(\hat{\bmu}_{\bw^*|\mathcal{M}_{AS}}, \bmu )$ plus the parametric rate $1/n$. If an estimator attains these benchmarks, it can be concluded that this estimator is minimax optimal in terms of the risk ratio, and cannot be further improved without imposing additional data assumptions.

%This theorem also implies that the asymptotic optimality should not be the sole justification of MA. Even in this simple setting with fixed dimension $p$, the asymptotic optimality with constant one cannot be satisfied in the general conditions. %Even in the fixed dimension, the asymptotic optimality over all non-nested candidate models cannot be satisfied unless addition assumption is imposed on the parameter space $\mathcal{M}$.

\begin{remark}
The minimax lower bounds established in Theorem~\ref{theo:lower} extend Theorem 6 in \cite{Peng2024Optimality} in several directions. First, they are derived under more general parameter spaces and dimensionality compared to the permutation space and the specific setting $p=n$ considered in \cite{Peng2024Optimality}. In addition, the lower bound in (\ref{eq:lower_1}) is asymptotically exact, rather than only characterizing the minimax rate in order.
\end{remark}

\subsection{Attainability}\label{subsec:attainability}

In this subsection, we introduce an MA estimator based on a Mallows-type criterion \eqref{eq:c}, which attains the minimax lower bounds established in Theorem~\ref{theo:lower}. The proposed method has three key features: it considers all univariate models as candidate models, imposes a hypercube constraint on the weight vector, and sets the penalty parameter $\lambda$ to adapt to the dimension $p$. We refer to this strategy as \textbf{A}veraging via \textbf{d}imension \textbf{a}daptive \textbf{p}enalty (Adap), which is constructed in two steps.

\begin{description}
  \item[Step 1: ]~Define the univariate candidate model set as $\mathcal{M}_U \triangleq \{ \{1 \}, \{2 \},\ldots,\{ p \}  \}$. The $j$-th candidate model is estimated by
\begin{equation}\label{eq:esti_j}
  \hat{\bmu}_{j} = \bpsi_j(\bpsi_j^{\top}\bpsi_j)^{-1}\bpsi_j^{\top}\by =  \tilde{\theta}_j \bpsi_j,
\end{equation}
where $\tilde{\theta}_j \triangleq n^{-1}\bpsi_j^{\top}\by $.
  \item[Step 2: ]~Estimate the model weights by
\begin{equation}\label{eq:MMA_cor}
\begin{split}
     & \hat{\bw}_{2} \triangleq \argmin_{\bw \in \mathcal{H}} \biggl\{n^{-1}\Bigl\| \by - \sum_{j=1}^{p}w_j \hat{\bmu}_{j} \Bigl\|^2 + 2\lambda^2_{2} \sigma^2 \bw^{\top} \boldsymbol{1} \biggl\}, \\
\end{split}
\end{equation}
where $\mathcal{H} \triangleq [0,1]^p$, $\lambda_{2} \triangleq \sqrt{ (2\log p)/n }$, and $\boldsymbol{1} \triangleq (1,\ldots,1)^{\top}$. The resulting Adap estimator is then given by
\begin{equation}\label{eq:MMA_est} 
  \hat{\bmu}_{\hat{\bw}_{2}|\mathcal{M}_U}=\sum_{j=1}^{p} \hat{w}_{2j} \hat{\bmu}_{j} = \sum_{j=1}^{p} \hat{w}_{2j} \tilde{\theta}_j \bpsi_j,
\end{equation}
where $\hat{w}_{2j}$ denotes the $j$-th element of $\hat{\bw}_{2}$.
\end{description}

\begin{theorem}\label{theo:upper}
Suppose that for each $1 \leq j \leq p$, the term $n^{-1}\bpsi_j^{\top}\beps$ follows a Gaussian distribution $N(0, \sigma^2/n)$. If $p$ is fixed, there must exist a constant $\bar{C}>1$ which is independent of $n$ such that
    \begin{equation*}
    \max_{\bmu \in \mathbb{R}^n} \frac{\mathbb{E} L_n\left( \hat{\bmu}_{\hat{\bw}_{2}|\mathcal{M}_U},  \bmu \right)}{n^{-1}+R_{n}(\hat{\bmu}_{\bw^*|\mathcal{M}_{AS}}, \bmu)} \leq \bar{C}.
  \end{equation*}
  If $p \to \infty$ as $n \to \infty$, then
  \begin{equation}\label{eq:upper_2}
    \max_{\bmu \in \mathbb{R}^n} \frac{\mathbb{E} L_n\left( \hat{\bmu}_{\hat{\bw}_{2}|\mathcal{M}_U},  \bmu \right)}{n^{-1}+R_{n}(\hat{\bmu}_{\bw^*|\mathcal{M}_{AS}}, \bmu)} \leq[1+o(1)]2\log p.
  \end{equation}
\end{theorem}

The Gaussian condition on $n^{-1}\bpsi_j^{\top}\beps$ can be satisfied when $\epsilon_1,\ldots,\epsilon_n$ are i.i.d. $N(0,\sigma^2)$. For cases where $\epsilon_i$'s do not satisfy the Gaussian assumption, we provide a detailed analysis in Corollary~\ref{coro:non-Gaussian}. Theorem~\ref{theo:upper} establishes that the Adap estimator $\hat{\bmu}_{\hat{\bw}_{2}|\mathcal{M}_U}$ achieves the minimax lower bound in terms of risk ratio given in Theorem~\ref{theo:lower}. Specifically, when $p$ is fixed (i.e., a standard parametric setting), the maximum risk ratio of the proposed estimator over all $\bmu \in \mathbb{R}^n$ remains bounded. When $p \to \infty$, the maximum risk ratio of $\hat{\bmu}_{\hat{\bw}_{2}|\mathcal{M}_U}$ matches the lower bound in (\ref{eq:lower_1}), indicating that $\hat{\bmu}_{\hat{\bw}_{2}|\mathcal{M}_U}$ is an optimal MA estimator for the all-subset MA task up to a parametric-rate term $1/n$.

Note that the maximum risk-ratio bounds in Theorem~\ref{theo:upper} hold over all $\bmu \in \mathbb{R}^n$, whereas the matching lower bounds in Theorem~\ref{theo:lower} are valid for any subset $\mathcal{C}(\Theta)$ with $\Theta^* \subseteq \Theta$. This implies that the cube $\Theta^*$ indeed characterizes the most difficult parameter region for achieving the optimal all-subset MA risk. The optimal all-subset MA risk in Theorem~\ref{theo:upper} is conditioned on a given orthogonal basis $\{\bpsi_1,\ldots,\bpsi_p\}$. Ideally, to make $R_{n}(\hat{\bmu}_{\bw^*|\mathcal{M}_{AS}}, \bmu)$ efficient, the basis should provide an \emph{economical representation} of the unknown mean vector $\bmu$. That is, the coefficients $\theta_j$ in (\ref{eq:complete}) should exhibit certain sparse pattern \citep[see, e.g.,][]{Beran2000}. In practice, Adap can be implemented based on PCs \citep{Jeffers1967}. Our numerical results in Section~\ref{sec:simu} indicate that this choice often leads to satisfactory performance across a variety of settings.

\begin{remark}

The weight constraint $\mathcal{H}$ has also been adopted by \cite{Ando2014, Ando2017}, \cite{lin2023optimal}, and \cite{Peng2023shrinkage} to develop MA procedures. In addition, different penalty choices in (\ref{eq:c}) have been considered, such as $\lambda_1 =\sqrt{1/n}$ in \cite{Hansen2007least} and $\lambda_3=\sqrt{(\log n)/n}$ in \cite{Zhang2020Parsimonious}. However, none of these methods are designed to achieve the optimal MA risk of all-subset models. 

\end{remark}

\begin{remark}
  \cite{Zhu2023Scalable} considered a similar procedure to (\ref{eq:MMA_cor}), where the penalty is set to $\lambda_1 = \sqrt{1/n}$. Their theoretical analysis follows the classical AOP principle aiming to achieve an asymptotic loss-ratio of one under the assumption $p/[nR_{n}(\hat{\bmu}_{\bw^*|\mathcal{M}_{AS}}, \bmu)] \to 0$. As pointed out in the paragraph following Condition~C.2 in \cite{Zhu2023Scalable}, this assumption can be satisfied when all candidate models are misspecified. A specific example provided by \cite{Zhu2023Scalable} is the case where $nR_{n}(\hat{\bmu}_{\bw^*|\mathcal{M}_{AS}}, \bmu) \asymp n$, which is reasonable and may occur under relatively large model misspecification. In contrast, our setup considers (\ref{eq:complete}), where $p$ can be either fixed or diverge with $n$, without imposing the restriction on the convergence rate of $R_{n}(\hat{\bmu}_{\bw^*|\mathcal{M}_{AS}}, \bmu)$ due to the model misspecification. 
\end{remark}

\subsection{Extensions}

For the setting where the error variance $\sigma^2$ is unknown, we construct a specific estimator $\hat{\sigma}^2$ under the Gaussian assumption. If $p<n$, define the projection matrix $\bP \triangleq \sum_{j=1}^p n^{-1} \bpsi_j \bpsi_j^{\top}$ and estimate the variance by 
  \begin{equation}\label{eq:var_est}
    \hat{\sigma}^2=\frac{\|(\bI-\bP) \by\|^2}{n-p}. 
  \end{equation}
  Define the estimated vector $\tilde{\bw}_2$ by replacing $\sigma^2$ in (\ref{eq:MMA_cor}) with $\hat{\sigma}^2$. We establish the following upper bounds on the risk ratio of the Adap estimator $\hat{\bmu}_{\tilde{\bw}_2|\mathcal{M}_U}$. 
  
  \begin{corollary}\label{coro:unknown-variance}
  Suppose that $\epsilon_1,\ldots,\epsilon_n$ are i.i.d. $N(0,\sigma^2)$, and the basis $\{\bpsi_1,\ldots,\bpsi_p\}$ satisfies Assumption~\ref{ass:complete}.  Then, there exists a constant $\bar{C}>0$ such that  
    \begin{equation}\label{eq:universal_case_3}
    \max_{\bmu \in \mathbb{R}^n} \frac{\mathbb{E} L_n\left( \hat{\bmu}_{\tilde{\bw}_2|\mathcal{M}_U},  \bmu \right)}{n^{-1}+R_{n}(\hat{\bmu}_{\bw^*|\mathcal{M}_{AS}}, \bmu)} \leq \begin{cases} \bar{C}, & p \text{\: is fixed} \\ [1+o(1)]2\log p , &  p \to \infty, \, p(\log p )^2 = O(n), \text{ as } n\to \infty.\end{cases} 
  \end{equation}

\end{corollary}

Corollary~\ref{coro:unknown-variance} shows that the Adap estimator with the variance estimator in (\ref{eq:var_est}) also achieves the minimax lower bounds established in Theorem~\ref{theo:lower}. Compared with the known-variance result in Theorem~\ref{theo:upper}, the upper bound for the unknown-variance case requires a slightly stronger restriction on the dimension $p$. Specifically, it typically requires $p = O(n/(\log n)^2)$, whereas the upper bound in (\ref{eq:upper_2}) remains valid for the larger dimension $p = n$. 

Note that Theorem~\ref{theo:upper} and Corollary~\ref{coro:unknown-variance} rely on the Gaussian assumption either on $n^{-1}\bpsi_j^{\top}\beps$ or on $\epsilon_i$. In the following, we investigate the performance of the Adap estimator when the Gaussian assumption is relaxed. We introduce a general Adap estimator, which is defined similarly to the estimator in (\ref{eq:esti_j})--(\ref{eq:MMA_est}), except that the weight vector is estimated by 
\begin{equation}\label{eq:general_MMA_cor}
     \hat{\bw}_{c} \triangleq \argmin_{\bw \in \mathcal{H}} \biggl\{n^{-1}\Bigl\| \by - \sum_{j=1}^{p}w_j \hat{\bmu}_{j} \Bigl\|^2 + 2\lambda^2_{c} \sigma^2\bw^{\top} \boldsymbol{1} \biggl\}, \\
\end{equation}
where $\lambda_c \triangleq \sqrt{ (c\log p)/n }$, and $c>0$ is a constant. Accordingly, the Adap estimator in (\ref{eq:MMA_est}) is extended to $\hat{\bmu}_{\hat{\bw}_{c}|\mathcal{M}_U}$, which reduces to the original estimator (\ref{eq:MMA_est}) when $c = 2$. 

\begin{corollary}\label{coro:non-Gaussian}

Suppose that $\{\bpsi_1, \ldots, \bpsi_p\}$ satisfy Assumption~\ref{ass:complete}, and that there exists a constant $C_0>0$ such that $\max_{j,i}|\psi_{ji}| \leq C_0$, where $\psi_{ji}$ denotes the $i$-th element of $\bpsi_j$. If Assumption~\ref{ass:4_moment} holds, then for any constant $c\geq \bar{c}$, there exists a constant $\bar{C}>2$ such that 
  \begin{equation}\label{eq:universal_case}
    \max_{\bmu \in \mathbb{R}^n} \frac{\mathbb{E} L_n\left( \hat{\bmu}_{\hat{\bw}_{c}|\mathcal{M}_U},  \bmu \right)}{n^{-1}+R_{n}(\hat{\bmu}_{\bw^*|\mathcal{M}_{AS}}, \bmu)} \leq \begin{cases} \bar{C}, & p \text{\: is fixed} \\ [1+o(1)]\bar{C}\log p , &  p \to \infty, \, p = O(n^{1/4}) \text{ as } n\to \infty,\end{cases} 
  \end{equation}
where $\bar{c}>0$ is some constant. Furthermore, if there exists a constant $B>0$ such that for all $k\geq 2$, $\mathbb{E}|\epsilon_i|^k \leq \frac{1}{2} k!B^{k-2} \sigma^2$, then for any constant $c>\bar{c}$, there exists a constant $\bar{C}>2$ such that  
  \begin{equation}\label{eq:universal_case_2}
    \max_{\bmu \in \mathbb{R}^n} \frac{\mathbb{E} L_n\left( \hat{\bmu}_{\hat{\bw}_{c}|\mathcal{M}_U},  \bmu \right)}{n^{-1}+R_{n}(\hat{\bmu}_{\bw^*|\mathcal{M}_{AS}}, \bmu)} \leq \begin{cases} \bar{C}, & p \text{\: is fixed} \\ [1+o(1)]\bar{C}\log p , &  p \to \infty \text{ as } n\to \infty.\end{cases} 
  \end{equation}
\end{corollary}

The statement in \eqref{eq:universal_case} suggests that, when the random errors satisfy the general fourth-moment condition in Assumption~\ref{ass:4_moment}, the Adap estimator $\hat{\bmu}_{\hat{\bw}_c|\mathcal{M}_U}$ with a suitably enlarged penalty parameter $\lambda_c=\sqrt{(c\log p)/n}$ can achieve the minimax-optimal risk ratio up to a constant factor when $p=O(n^{1/4})$. Under the additional Bernstein condition on $\epsilon_i$, this rate optimality (up to a constant factor) can be extended to the wider range of $p$, as shown in (\ref{eq:universal_case_2}). The dependence of $\hat{\bw}_{c}$ on the unknown variance $\sigma^2$ in (\ref{eq:general_MMA_cor}) does not pose a significant difficulty in this case if one does not require the exact control of the constant $\bar{C}$. For example, if an upper bound $\bar{\sigma}^2$ on $\sigma^2$ is available or can be estimated to certain accuracy, then the unknown variance $\sigma^2$ in (\ref{eq:general_MMA_cor}) can be absorbed into the penalty parameter $\lambda_c$. Specifically, we replace the factor $\lambda^2_{c} \sigma^2$ in (\ref{eq:general_MMA_cor}) by $\lambda^2_{c}$. Consequently, we may choose $c$ sufficiently large such that $c>\bar{c}\bar{\sigma}^2$. Under this choice, the results in Corollary~\ref{coro:non-Gaussian} continue to hold for the generalized Adap estimator $\hat{\bmu}_{\hat{\bw}_{c}|\mathcal{M}_U}$ for some $\bar{C}$ large enough. In practice, if the upper bound $\bar{\sigma}^2$ is also unavailable, one can instead treat $c\sigma^2$ as an unknown tuning parameter and select it in a data-driven manner. Note that the minimax lower bounds in Theorem~\ref{theo:lower} are actually established for the Gaussian errors. It remains unknown whether the constant inflation in \eqref{eq:universal_case} and \eqref{eq:universal_case_2} is an intrinsic phenomenon under these general error assumptions.

\subsection{The ensemble effect of several MS procedures}\label{sec:MS_dis} 

The proposed Adap estimator (\ref{eq:MMA_est}) is closely related to several classical MS procedures in the existing literature. From the proof in Section~\ref{subsubsec:equivalent} of the Appendix, we see that the estimated coefficients in (\ref{eq:MMA_est}) have the closed form
\begin{equation*}
  \hat{w}_{2j} \tilde{\theta}_j= \left( 1 - \frac{\lambda^2_2 \sigma^2}{\tilde{\theta}_j^2} \right)_+\tilde{\theta}_j,\quad j =1,\ldots,p,
\end{equation*}
which coincides with the non-negative garrotte estimator proposed by \citet{Breiman1995garrote}. The MS consistency of such estimator has been established in \citet{Zou2006} and \citet{Yuan2007}, and its minimax risk-ratio optimality with respect to the optimal all-subset MS risk was demonstrated in \citet{Gao1998}. However, to the best of our knowledge, it was previously unknown that the non-negative garrotte estimator also has certain ensemble effect, as established in Theorem~\ref{theo:upper} through its attainment of the minimax optimal rate to $R_{n}(\hat{\bmu}_{\bw^*|\mathcal{M}_{AS}}, \bmu)$.

The risk inflation criterion (RIC) \citep{Foster1994} and the Lasso \citep{Tibshirani1996lasso} are two well-known MS strategies. Under the orthogonal design setting in Assumption~\ref{ass:complete}, both reduce to the soft-thresholding estimator \citep{Donoho1994Ideal}:
\begin{equation}\label{eq:soft}
  \hat{\bmu}_{\mathrm{ST}}= \sum_{j=1}^{p} \mathrm{sgn}(\tilde{\theta}_j) \bigl( |\tilde{\theta}_j| - \lambda_2\sigma \bigl)_+ \bpsi_j.
\end{equation}
In addition to~\eqref{eq:soft}, a closely related method is the hard-thresholding estimator:
\begin{equation}\label{eq:hard}
  \hat{\bmu}_{\mathrm{HT}}= \sum_{j=1}^{p} 1_{\{ |\tilde{\theta}_j| >\lambda_2\sigma\}}\tilde{\theta}_j \bpsi_j.
\end{equation}
By connecting the results in Section~4 of \cite{Donoho1994Ideal} to our MA framework, we find that both $\hat{\bmu}_{\mathrm{ST}}$ and $\hat{\bmu}_{\mathrm{HT}}$ achieve the optimal all-subset MA in terms of the minimax risk ratio, as stated in the following corollary.
\begin{corollary}\label{coro:upper_soft_hard}
  Let $\hat{\bmu}_{\cdot\mathrm{T}}$ denote either $\hat{\bmu}_{\mathrm{ST}}$ or $\hat{\bmu}_{\mathrm{HT}}$. Under the same assumptions as in Theorem~\ref{theo:upper}, if $p \to \infty$ as $n \to \infty$, then
  \begin{equation}\label{eq:upper_4}
    \max_{\bmu \in \mathbb{R}^n} \frac{R_n\left( \hat{\bmu}_{\cdot\mathrm{T}},  \bmu \right)}{n^{-1}+R_{n}(\hat{\bmu}_{\bw^*|\mathcal{M}_{AS}}, \bmu)} \leq[1+o(1)]2\log p.
  \end{equation}
\end{corollary}

Due to the regularization effects that induce sparse solutions, techniques such as Lasso and RIC have been regarded as MS methods in some of the literature, and therefore may suffer from MS uncertainty that could potentially be alleviated through MA schemes. However, our analysis in this subsection demonstrates that under orthogonal settings, certain properly tuned MS procedures can in fact attain the fastest possible convergence rate to the optimal all-subset MA risk. This finding may address the open question raised at the end of Section~6 of \cite{Wang2009} regarding the relationship between MA and penalized MS approaches. The unveiled ensemble effect underlying these MS methods suggests that they can exhibit competitive predictive performance compared to MA when the variables are not highly correlated. The numerical results presented in Section~\ref{sec:simu} support this theoretical understanding.

\section{Numerical studies}\label{sec:simu}

In this section, we conduct several numerical experiments to illustrate the theoretical results developed in Sections~\ref{sec:general}--\ref{sec:all-subset} and to compare the performance of several MA and MS procedures.

\subsection{Assessing the achievability of the optimal all-nested MA risk}\label{sec:simu_1}

\begin{figure}[t!]
\centering
%\vspace{-0.35cm}	
%\subfigtopskip=2pt
%\subfigbottomskip=2pt
%\subfigcapskip=-5pt
%\setlength{\abovecaptionskip}{2pt}

\subfigure[t-distribution]{
\includegraphics[width=5in]{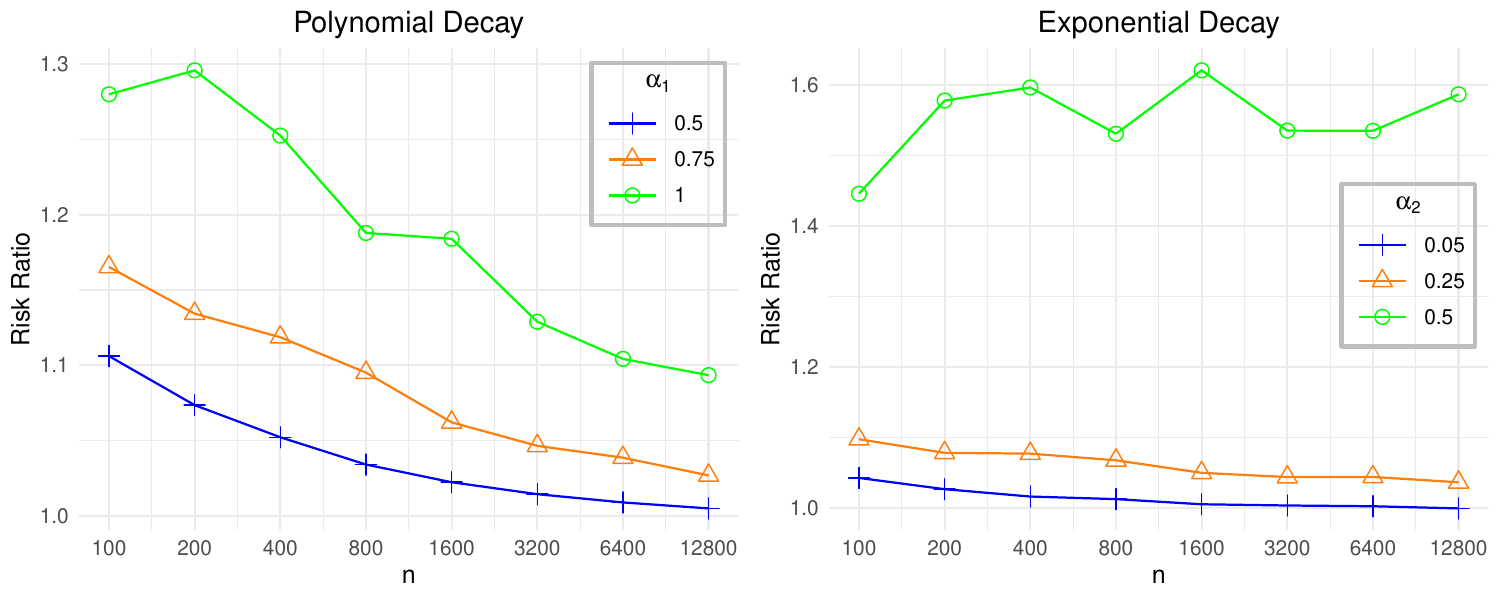}}%加大宽度

\subfigure[Pareto distribution]{
\includegraphics[width=5in]{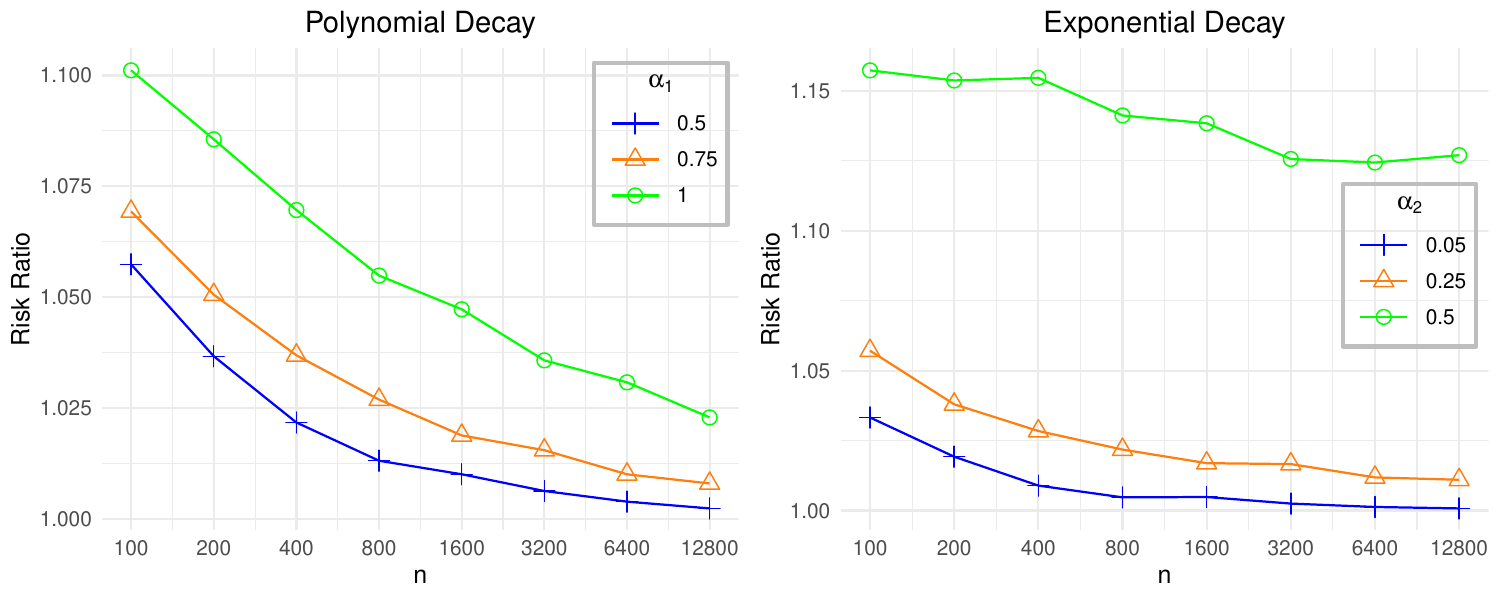}}%加大宽度

\caption{Risk ratio of the MMA estimator $\hat{\boldsymbol{\mu}}_{\hat{\mathbf{w}}_1|\mathcal{M}_G}$ under the polynomially and exponentially decaying coefficients. Results for $t$-distributed errors are shown in row (a), and those for Pareto-distributed errors are shown in row (b).}
\label{fig:risk_simu}
\end{figure}

The data are generated from (\ref{eq:model}) and (\ref{eq:complete}) with the canonical basis $\{\bpsi_j = \sqrt{n}\be_j,j=1,\ldots,n\}$ and $p=n$, where $\be_j \in \mathbb{R}^n$ is the vector with 1 in its $j$-th element and 0 elsewhere. The coefficients $\theta_j, j =1,\ldots, p$ in (\ref{eq:complete}) are set as the ordered sequence $\theta_{(j)}, j =1,\ldots, p$ under two settings:
\begin{description}
  \item[Polynomial decay:]$\theta_{(j)} = j^{-\alpha_1}$, with $0.5 < \alpha_1 < \infty$.
  \item[Exponential decay:]$\theta_{(j)} = \exp(-j^{\alpha_2})$, with $0 < \alpha_2 < \infty$.
\end{description}
The random error terms $\epsilon_1,\ldots,\epsilon_n$ are i.i.d. from two heavy-tailed distributions. The first is a $t$-distribution with $\mathrm{df} = 5$. The second is a Pareto distribution, where $|\epsilon_i|$ follows a Pareto Type I distribution with shape parameter 5 and scale parameter 1. For each distribution, the variance $\sigma^2$ is adjusted such that the signal-to-noise ratio (SNR) $ \sum_{j=1}^{n} \theta_j^2 / \sigma^2 $ equals 5. The sample size $n$ increases from 100 to 12800 on a logarithmic scale. The risk ratio is computed as the averaged loss of the nested MMA estimator $\hat{\bmu}_{\hat{\bw}_1|\mathcal{M}_G}$ over 1000 replications, divided by the optimal all-nested MA risk. The results are presented in Figure~\ref{fig:risk_simu}.

From the left panels of Figure~\ref{fig:risk_simu}, we observe that the risk ratios in the polynomial decay case gradually decrease toward 1 as the sample size increases. The exponential case with $\alpha_2 = 0.05$ also exhibits an obvious downward trend, which supports the AOP result in Section~\ref{sec:nested} that the optimal nested MA risk can be attained when $R_{n}(\hat{\bmu}_{\bw^*|\mathcal{M}_{AN}}, \bmu)$ converges slower than $(\log n)^4/n$. In contrast, for the exponential case with $\alpha_2 = 0.5$, a substantial gap between the risk ratio and 1 exists even when the sample sizes are sufficiently large, suggesting that it is difficult to achieve $R_{n}(\hat{\bmu}_{\bw^*|\mathcal{M}_{AN}}, \bmu)$ when the coefficients decay fast.

\subsection{Assessing the achievability of the optimal all-subset MA risk}\label{sec:simu_2}

\begin{figure}[!htbp]
    \centering
    \subfigure[Polynomially decaying signal with $\alpha_1 = 1$]{
    \begin{minipage}[t]{1\linewidth}
    \centering
       \includegraphics[width=5in]{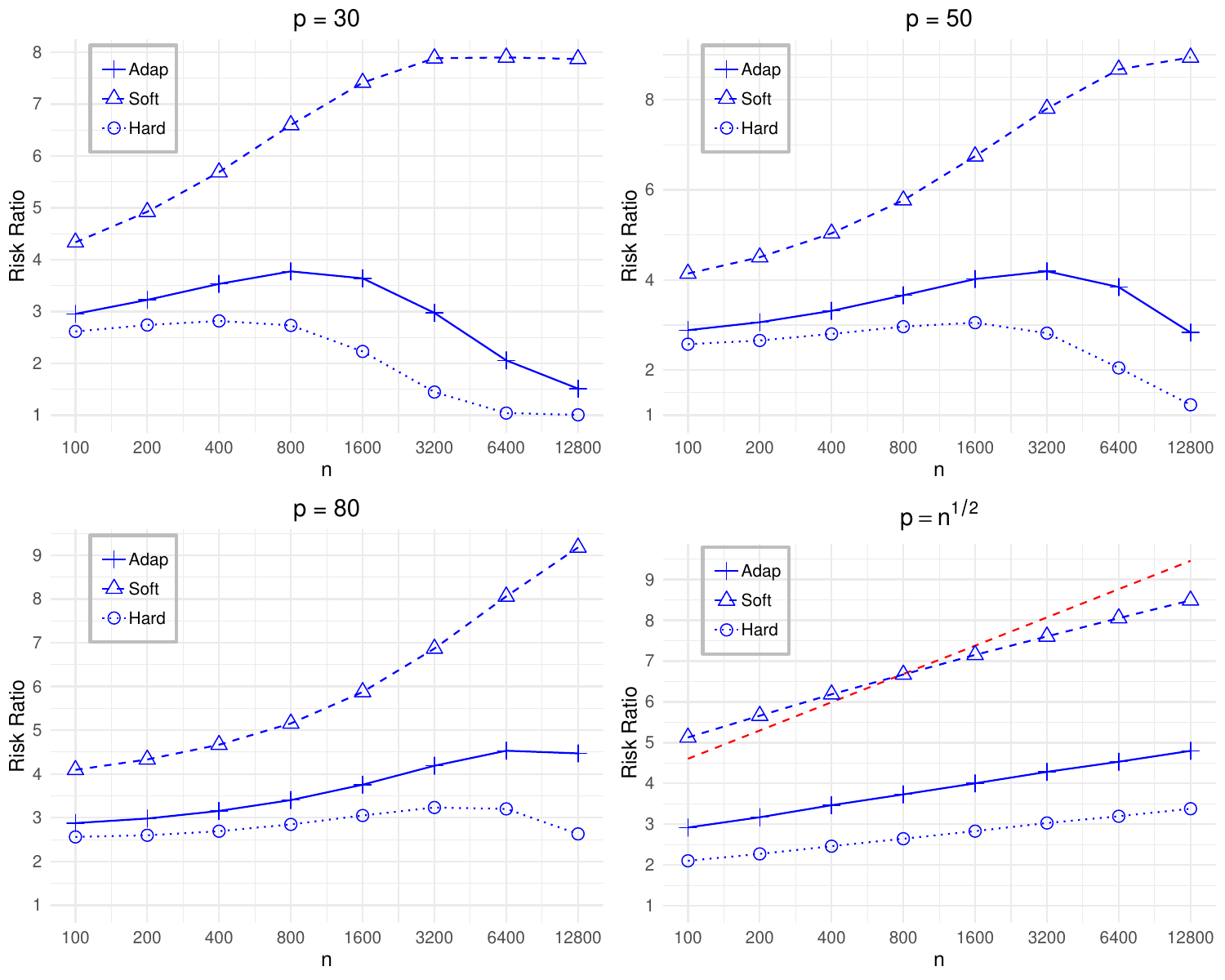}
       % \hspace{2cm}
    \end{minipage}
    }

    \subfigure[Exponentially decaying signal with $\alpha_2 = 0.5$]{
    \begin{minipage}[t]{1\linewidth}
    \centering
       \includegraphics[width=5in]{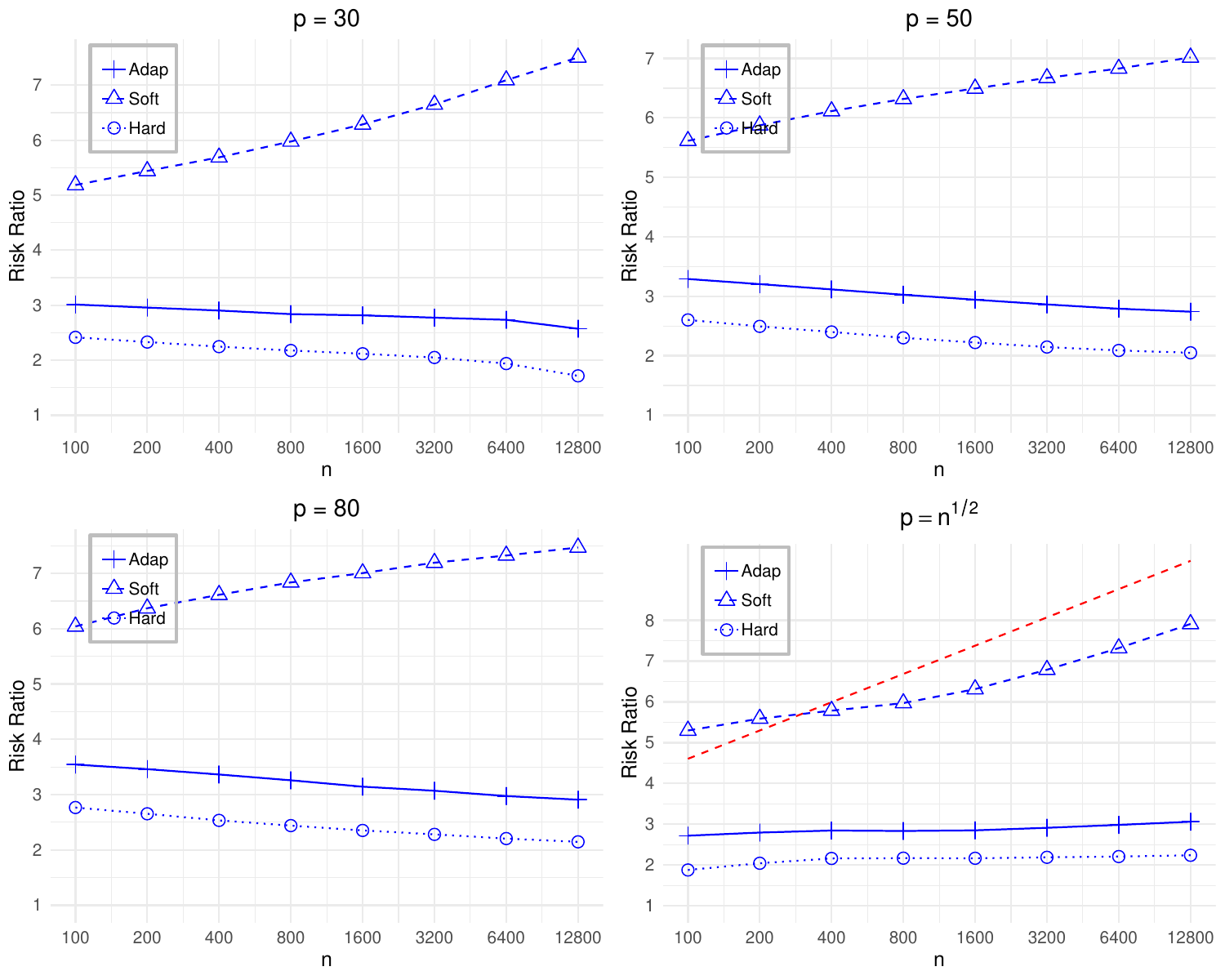}
       % \hspace{2cm}
    \end{minipage}
    }

    \caption{The risk ratios of the three competing methods under different signal decay scenarios. In each subfigure, the red dashed line in the bottom-right panel represents the curve of $2\log p$.}
    \label{fig:risk_non}
\end{figure}

The data are generated from the same model as that in Section~\ref{sec:simu_1} with $p=30,50, 80$, and $\lfloor n^{1/2}\rfloor$. In each simulation replication, the coefficients $\theta_1,\ldots,\theta_p$ in (\ref{eq:complete}) are generated as a random permutation of the ordered sequence $\theta_{(1)},\ldots,\theta_{(p)}$. This setup is designed to mimic scenarios where the importance of variables is unknown to statisticians, under which the nested MA strategy is not favorable. The random error terms $\epsilon_1,\ldots,\epsilon_n$ are i.i.d. from $N(0,\sigma^2)$. We plot the risk ratios of the Adap estimator (\ref{eq:MMA_est}) and the soft/hard-thresholding estimators (\ref{eq:soft})--(\ref{eq:hard}) relative to the optimal all-subset MA risk. The results are presented in Figure~\ref{fig:risk_non}.

From Figure~\ref{fig:risk_non}, we observe that although all three methods have been theoretically shown to be minimax optimal, their empirical performance differs under the two specific simulation settings considered. Specifically, the hard-thresholding and Adap estimators perform better than the soft-thresholding estimator. This observation is due to the soft-thresholding estimator tends to overshrink large signals and thus incurs greater bias. For a more detailed theoretical comparison of the thresholding estimators, see \cite{Guo2024Signal}. In the fixed-dimensional setting, the risk ratios of both the Adap and hard-thresholding estimators remain bounded. In the diverging-dimension regime, the risk ratios of all three methods lie below the curve $2\log p$, which support the minimax upper bounds in Theorem~\ref{theo:upper} and Corollary~\ref{coro:upper_soft_hard}.

\subsection{Comparing several different procedures}\label{sec:compare}

In this subsection, we compare several methods through simulation studies and a real data example, where the orthogonal basis in Assumption~\ref{ass:complete} is constructed from the PCs \citep[see, e.g.,][]{Jeffers1967}. For the simulation, the data are generated from a PC regression model $\by = \bU\btheta + \beps$, where $\bU = [\bu_1, \ldots, \bu_p]$ is obtained from the SVD $\bX = \bU \bD \bV^{\top}$, the diagonal matrix $\bD$ contains singular values $\lambda_1 \geq \lambda_2 \geq \cdots \geq \lambda_p > 0$, the noise term $\beps \sim N(\boldsymbol{0}, \sigma^2 \bI_n)$, and $p$ denotes the rank of $\bX$. The matrix $\bX$ follows a multivariate normal distribution $N(\boldsymbol{0}, \Sigma)$, where $\Sigma = (0.5^{|i-j|})_{1 \leq i,j \leq d}$, $n = 500$, and $d = 1000$. In the simulation, we consider both the ordered and unordered coefficient $\theta_j$, which are the same as those considered in Sections~\ref{sec:simu_1}--\ref{sec:simu_2}. The ordered cases are designed to mimic scenarios in which the signal strength projected onto the PCs decays in alignment with the order of singular values. This phenomenon has been observed in some classical statistical problems \citep{Hocking1976}; see also the real data example below. However, such alignment does not always occur \citep[see, e.g.,][]{Bair01032006}. The unordered cases are thus used to model more general data structure.

\begin{figure}[t!]
\centering
%\vspace{-0.35cm}	
%\subfigtopskip=2pt
%\subfigbottomskip=2pt
%\subfigcapskip=-5pt
%\setlength{\abovecaptionskip}{2pt}

\subfigure[Ordered coefficient]{
\includegraphics[width=5.5in]{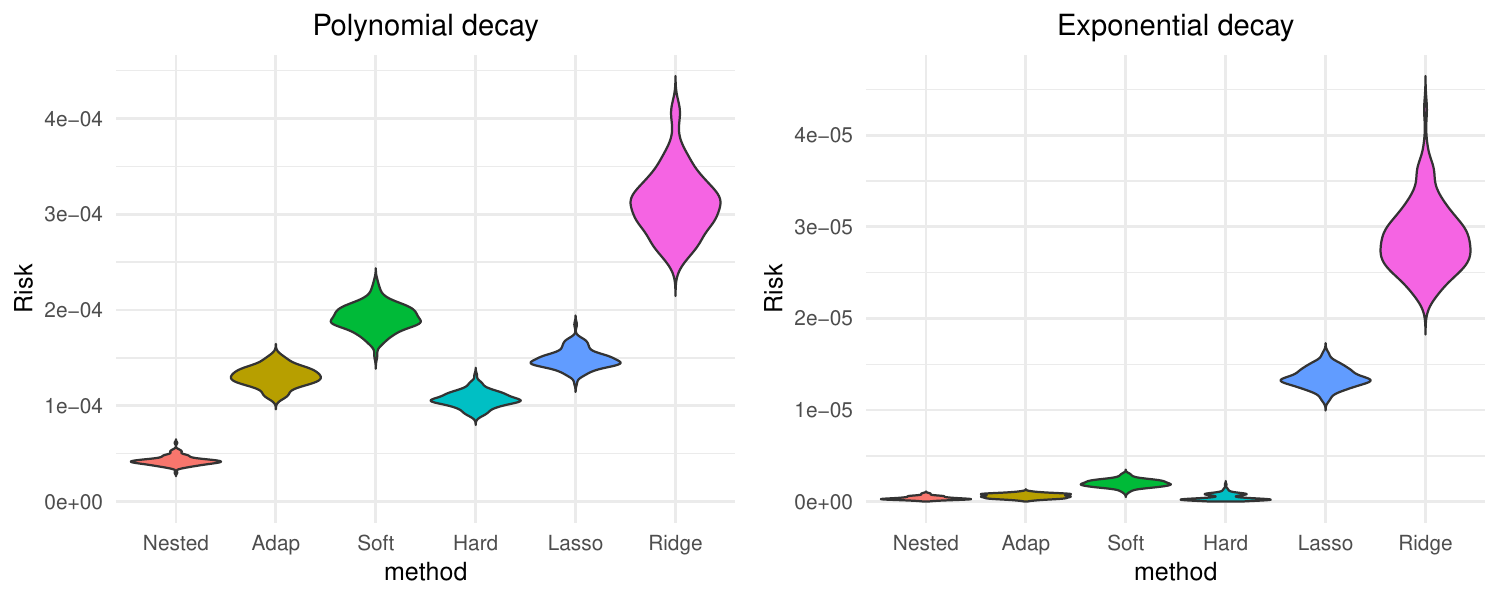}}%加大宽度

\subfigure[Unordered coefficient]{
\includegraphics[width=5.5in]{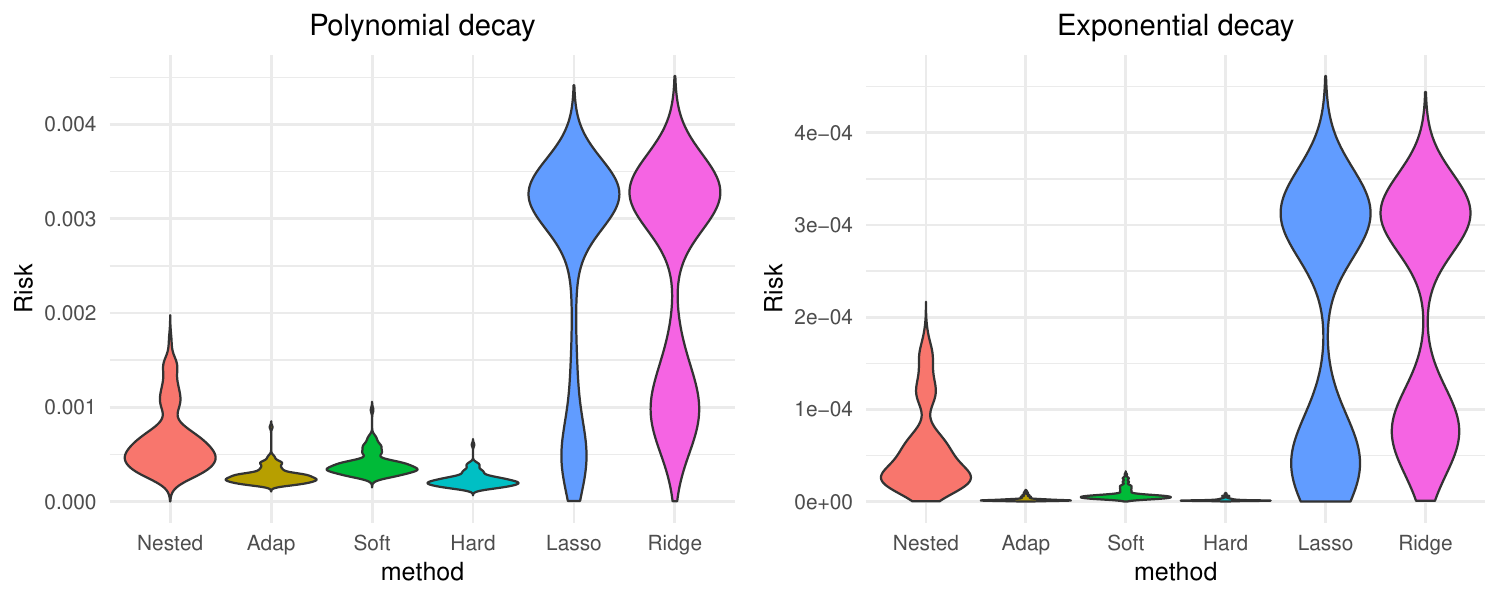}}%加大宽度

\caption{Risk comparison of six competing methods. Results for the ordered cases are presented in row (a), and those for the unordered cases are shown in row (b).}
\label{fig:risk_compare}
\end{figure}

Since each $\bu_j$ has unit norm, we define an orthogonal basis $\{\bpsi_1, \ldots, \bpsi_p\}$ by setting $\bpsi_j = \sqrt{n} \bu_j$ for $j = 1, \ldots, p$. Based on this basis, we construct the nested MMA estimator $\hat{\bmu}_{\hat{\bw}_1|\mathcal{M}_G}$ described in Section~\ref{sec:nested}, the Adap estimator (\ref{eq:MMA_est}), the soft-thresholding estimator (\ref{eq:soft}), and the hard-thresholding estimator (\ref{eq:hard}) as competing methods. In addition, we include the Lasso method \citep{Tibshirani1996lasso} and ridge regression \citep{hoerl1970ridge} as representative modeling procedures based on the design matrix $\bX$. The regularization parameters in these two methods are selected via 5-fold cross-validation. From the simulation results in Figure~\ref{fig:risk_compare}, we observe that when $\theta_j, j = 1, \ldots, p$ are ordered, the nested MMA estimator performs quite well. In contrast, when the ordering structure is violated, as illustrated in Figure~\ref{fig:risk_compare}(b), the Adap estimator and soft/hard-thresholding estimators appear to be more efficient. It is also worth noting that the Lasso, when applied to the original design matrix $\bX$, performs poorly in our simulation. This is not surprising, as the Lasso is suboptimal when the regressors are correlated, regardless of the choice of tuning parameters \citep[see, e.g.,][]{pathak2024design}.

PC regression is a popular method in chemometric data analysis. We further compare the above six methods using the gas chromatographic (GC) retention dataset \citep{Varmuza2009}. The dataset contains 209 observations of the GC retention indices of polycyclic aromatic compounds as the response variable $\by$, along with 467 molecular descriptors as the predictors $\bX$. The dataset is available at https://github.com/cran/chemometrics. The implementation of the competing methods is the same as in the simulation studies, except for the treatment of the unknown $\sigma^2$. For $\hat{\bmu}_{\hat{\bw}_1|\mathcal{M}_G}$, the variance is estimated from the residuals of the regression estimator based on the first 100 PCs. For the Adap estimator, we consider the general form in (\ref{eq:general_MMA_cor}) and treat $c\sigma^2$ as the tuning parameter, which is selected via five-fold cross-validation. The similar approach is applied to the soft/hard-thresholding estimators. 

After 100 repetitions of random training-testing splits, the average prediction errors (with standard errors in parentheses) for the methods Nested, Adap, Soft, Hard, Lasso, and Ridge are 89.90 (8.33), 98.42 (9.55), 102.10 (9.22), 121.68 (11.17), 125.01 (37.37), and 265.24 (42.56), respectively. Indeed, in the GC retention data set, the correlations between $\by$ and the PCs exhibit a decreasing pattern as the PC index increases. For example, the absolute correlations between $\by$ and the PCs are 0.94, 0.16, and 0.09 for the 1st, 10th, and 20th PCs, respectively. Therefore, it is expected that the nested MMA estimator $\hat{\bmu}_{\hat{\bw}_1|\mathcal{M}_G}$ still performs well in this setting, which is also consistent with the observations in Figure~\ref{fig:risk_compare}(a). However, if we exchange the order of the 1st--20th and 131st--150th PCs and re-examine the six methods, the prediction errors of Nested, Adap, Soft, Hard, Lasso, and Ridge are 195.22 (21.20), 97.22 (8.56), 110.66 (9.98), 120.60 (12.54), 103.35 (21.53), and  290.80 (57.41), respectively. This result shows that the nested MMA estimator is sensitive to the specific ordering of the PCs, whereas the Adap and the soft/hard-thresholding estimators are more robust to the unordered settings. This observation is also consistent with the results shown in Figure~\ref{fig:risk_compare}(b).

\section{Conclusion and discussion}\label{sec:discuss}

This paper addresses two important problems in the theory and application of Mallows-type MA. First, we establish a finite-sample risk guarantee for the MMA estimator. The results are derived under general candidate model constructions, without imposing assumptions on the model structure or regressor design.

The second part of this paper focuses on specific candidate constructions. In our setup, the candidate models for MA are formed using different subsets of a given orthogonal basis. This assumption is natural and mild in the case of nested model spaces, as the nesting inherently induces a basis with orthogonal properties \citep{Xu2022model, Peng2024Optimality}. Moreover, in establishing the minimax lower bound for the optimal all-subset MA risk, the orthogonality constraint is a reasonable simplification, as it represents the most fundamental setting for analyzing the statistical limit of combining all-subset least squares estimators. %Notably, the lower bound derived under the orthogonal setup also serves as a lower bound for the general case beyond the orthogonal scenario, since the orthogonal design is a special case in the broader regressor designs.

The Adap estimator (\ref{eq:MMA_est}), which achieves the minimax optimal rate for all-subset MA, is constructed based on a given orthogonal basis. From a practical standpoint, such a basis can be obtained through an orthogonalization algorithm. Our numerical results suggest that the SVD of $\bX$ provides a viable method for constructing this basis. From a theoretical standpoint, the orthogonal setting serves as a starting point for understanding MS procedures \citep[see, e.g.,][]{Barron1999, birge2001gaussian, massart2007concentration}. In the context of MA, however, such a foundational understanding remained limited even in this basic setting prior to our work. Our paper takes a step toward filling this gap.

%when considering all-subset candidate models, the situation becomes significantly more complex, especially under general regression designs where orthogonality does not hold. the derivation of the upper bound—and, in particular, the successful implementation of the proposed Adaptive (Adap) estimator—relies crucially on the orthogonal basis construction. This dependence highlights the technical barrier in extending our approach to non-orthogonal settings.

Extending the all-subset MA theory developed in this paper to more general regressor design settings remains a challenging problem. In the context of all-subset MS, various relaxed forms of orthogonality have been proposed to establish the optimality of penalized MS methods \citep[see, e.g.,][]{Candes2006, Bickel2009, Meinshausen2009, Raskutti2011, Bellec2018Slope}. However, it is still an open question how to formulate analogous and suitable assumptions for all-subset MA, where the focus lies in achieving optimal model combination. Moreover, without any restrictions on the correlations among regressors, an all-subset comparison approach becomes essential for achieving the optimal rate \citep[see, e.g.,][]{yang1999model, Wang2014}, and the associated MS problem escalates to NP-hard complexity \citep[see, e.g.,][]{Natarajan1995, Zhang2014Lower}. A theoretical framework that addresses both the methodological and computational complexities of all-subset MA under the general correlation structures is still lacking. We leave these problems as directions for future research.

\section*{Acknowledgments}

The author would like to thank Professor Xinyu Zhang for several insightful discussions since the fall of 2021, which have broadened my knowledge and deepened my understanding about model averaging. The author also thanks Professor Yuhong Yang for helpful comments on an earlier version of this paper.

\newpage
\appendix

\section*{Appendix}\label{appendix}
\addcontentsline{toc}{section}{Appendix}
%\addtocontents{toc}{\setcounter{tocdepth}{1}}
%\renewcommand{\theequation}{A.\arabic{equation}}
%\renewcommand{\thesection}{A.\arabic{section}}
\renewcommand{\thesection}{\Alph{section}}
\numberwithin{equation}{section}

\section{Proof of the results in Section~\ref{sec:general}}\label{sec:a:proof_general}

% Note: in this sample, the section number is hard-coded in. Following
% proper LaTeX conventions, it should properly be coded as a reference:

%In this appendix we prove the following theorem from
%Section~\ref{sec:textree-generalization}:

\subsection{Proof of Proposition~\ref{theo:sharp}}

Given an arbitrary candidate model set $\mathcal{M}$, recall that $\bP_m$ denotes the projection matrix associated with the $m$-th candidate least squares estimator, and $\bP(\bw) = \sum_{m=1}^{M_n}w_m\bP_m$. For notational simplicity, we define $\bA(\bw) \triangleq \bI - \bP(\bw)$ in the proof.

By the definition of $\hat{\bw}_1$ in Section~\ref{sec:MMA_general}, we have
\begin{equation}\label{eq:A1}
\begin{split}
   n\mathbb{E}C_n(\hat{\bw}_1|\mathcal{M},\lambda_1) &\leq n\mathbb{E}C_n(\bw^*|\mathcal{M},\lambda_1) = \mathbb{E}\left[ \| \by - \bP(\bw^*) \by \|^2 + 2 \hat{\sigma}^2 \tr \bP(\bw^*) \right]\\
   &= \| \bA(\bw^*) \bmu \|^2 + \sigma^2 \tr \bA^2(\bw^*)  + 2\mathbb{E}(\hat{\sigma}^2)\tr \bP(\bw^*)  \\
   & = \| \bA(\bw^*) \bmu \|^2 + \sigma^2\tr \bP^2(\bw^*)+2(\mathbb{E}\hat{\sigma}^2 - \sigma^2)\tr \bP(\bw^*) + n\sigma^2\\
     &   = n R_{n}(\hat{\bmu}_{\bw^*|\mathcal{M}}, \bmu)+2(\mathbb{E}\hat{\sigma}^2 - \sigma^2)\tr \bP(\bw^*)+ n\sigma^2, 
\end{split}
\end{equation}
where the last equality follows from the fact that for any non-random $\bw$, 
\begin{equation}\label{eq:leavv1}
  n R_{n}(\hat{\bmu}_{\bw|\mathcal{M}}, \bmu) = \mathbb{E} \| \hat{\bmu}_{\bw|\mathcal{M}} - \bmu \|^2 = \mathbb{E} \| \bP(\bw)\beps - \bA(\bw)\bmu \|^2  = \| \bA(\bw) \bmu \|^2 + \sigma^2\tr \bP^2(\bw). 
\end{equation}
The loss function of the MMA estimator can be decomposed as
\begin{equation}\label{eq:A2}
\begin{split}
   n L_{n}(\hat{\bmu}_{\hat{\bw}_1|\mathcal{M}}, \bmu) & =nC_n(\hat{\bw}_1|\mathcal{M},\lambda_1) - \| \beps \|^2 - 2 \left\langle \bA(\hat{\bw}_1)\bmu, \beps \right\rangle \\
     & + 2 \left[\beps^{\top}\bP(\hat{\bw}_1)\beps - \sigma^2\tr\bP(\hat{\bw}_1)\right] + 2(\sigma^2 - \hat{\sigma}^2)\tr\bP(\hat{\bw}_1).
\end{split}
\end{equation}
Combining inequalities \eqref{eq:A1}--\eqref{eq:A2}, we obtain
\begin{equation}\label{eq:sa_lu}
  \begin{split}
     n\mathbb{E}L_{n}(\hat{\bmu}_{\hat{\bw}_1|\mathcal{M}}, \bmu) & \leq nR_{n}(\hat{\bmu}_{\bw^*|\mathcal{M}}, \bmu)- 2 \mathbb{E}  \left\langle \bA(\hat{\bw}_1)\bmu, \beps \right\rangle  + 2\mathbb{E} \left[\beps^{\top}\bP(\hat{\bw}_1)\beps - \hat{\sigma}^2\tr\bP(\hat{\bw}_1)\right]\\
     &+2(\mathbb{E}\hat{\sigma}^2 - \sigma^2)\tr \bP(\bw^*) - 2\mathbb{E}\left[(\hat{\sigma}^2 - \sigma^2)\tr \bP(\hat{\bw}_1)\right]\\
     & \leq nR_{n}(\hat{\bmu}_{\bw^*|\mathcal{M}}, \bmu) + 2\mathbb{E} \sup_{\bw \in \mathcal{W}} \left| \left\langle \bA(\bw)\bmu, \beps \right\rangle \right|\\
     &+ 2\mathbb{E} \sup_{\bw \in \mathcal{W}} \left|\beps^{\top}\bP(\bw)\beps - \sigma^2\tr\bP(\bw)\right|+4\mathbb{E}\left|\hat{\sigma}^2 - \sigma^2\right|\sup_{\bw \in \mathcal{W}}\tr\bP(\bw), 
  \end{split}
\end{equation}
where the last inequality follows from the absolute value inequality and upper bounding the terms involving the random $\hat{\bw}_1$ by the uniform upper bound over all $\bw \in \mathcal{W}$.

We first bound the term $\mathbb{E} \sup_{\bw \in \mathcal{W}} | \langle \bA(\bw)\bmu, \beps \rangle |$ in (\ref{eq:sa_lu}). Since $\langle \bA(\bw)\bmu, \beps \rangle$ is a linear function in $\bw$, its supremum and infimum are attained at a vertex of $\mathcal{W}$. Thus,
\begin{equation*}
  \mathbb{E} \sup_{\bw \in \mathcal{W}} \left| \left\langle \bA(\bw)\bmu, \beps \right\rangle \right| \leq \mathbb{E} \max_{1 \leq m \leq M_n} \left| \left\langle (\bI - \bP_m)\bmu, \beps \right\rangle \right|.
\end{equation*}
Applying standard tail probability bounds, we derive
\begin{equation}\label{eq:shabi_lifan}
\begin{split}
   \mathbb{P}\left( \max_{1 \leq m \leq M_n} \left| \left\langle (\bI - \bP_m)\bmu, \beps \right\rangle \right| > t \right)  & \leq \sum_{m=1}^{M_n} \mathbb{P}\left(  \left| \left\langle (\bI - \bP_m)\bmu, \beps \right\rangle \right| > t \right) \\
     & \leq \sum_{m=1}^{M_n} \frac{\mathbb{E}\left\langle (\bI - \bP_m)\bmu, \beps \right\rangle^2}{t^2}\\
     & \leq \frac{\sigma^2 \sum_{m=1}^{M_n}\| (\bI - \bP_m)\bmu \|^2}{t^2},
\end{split}
\end{equation}
where the first inequality follows from the union bound, the second from Markov's inequality, and the third from
\begin{equation*}\label{eq:sa_lu2}
\begin{split}
   \mathbb{E}\left\langle (\bI - \bP_m)\bmu, \beps \right\rangle^2  & = \mathbb{E} [\beps^{\top}(\bI - \bP_m)\bmu \bmu^{\top}(\bI - \bP_m)\beps] = \sigma^2 \| (\bI - \bP_m)\bmu \|^2. \\
\end{split}
\end{equation*}
Integrating the tail probability bound in (\ref{eq:shabi_lifan}) yields
\begin{equation*}
  \begin{split}
     &\mathbb{E} \max_{1 \leq m \leq M_n} \left| \left\langle (\bI - \bP_m)\bmu, \beps \right\rangle \right| = \int_{0}^{\infty}\mathbb{P}\left( \max_{1 \leq m \leq M_n} \left| \left\langle (\bI - \bP_m)\bmu, \beps \right\rangle \right| > t \right) dt \\
       & \leq \int_{0}^{\infty} \min\left(1, \frac{\sigma^2 \sum_{m=1}^{M_n}\| (\bI - \bP_m)\bmu \|^2}{t^2}\right)dt\\
       & = \int_{0}^{\sigma \sqrt{\sum_{m=1}^{M_n}\| (\bI - \bP_m)\bmu \|^2}} 1 dt + \int_{\sigma \sqrt{\sum_{m=1}^{M_n}\| (\bI - \bP_m)\bmu \|^2}}^{\infty} \frac{\sigma^2 \sum_{m=1}^{M_n}\| (\bI - \bP_m)\bmu \|^2}{t^2} dt\\
       & = 2 \sigma \sqrt{\sum_{m=1}^{M_n}\| (\bI - \bP_m)\bmu \|^2}.
  \end{split}
\end{equation*}
Thus, we establish the bound
\begin{equation}\label{eq:part11}
  \mathbb{E} \sup_{\bw \in \mathcal{W}} \left| \left\langle \bA(\bw)\bmu, \beps \right\rangle \right| \leq 2 \sigma \sqrt{\sum_{m=1}^{M_n}\| (\bI - \bP_m)\bmu \|^2}.
\end{equation}

We then establish an upper bound for $\mathbb{E} \sup_{\bw \in \mathcal{W}} |\beps^{\top}\bP(\bw)\beps - \sigma^2\tr\bP(\bw)|$. Since the term $\beps^{\top}\bP(\bw)\beps - \sigma^2\tr\bP(\bw)$ is also linear in $\bw$, the supremum and infimum over $\mathcal{W}$ occur at the vertices of the simplex. Consequently, we obtain
\begin{equation*}
  \mathbb{E} \sup_{\bw \in \mathcal{W}} \left|\beps^{\top}\bP(\bw)\beps - \sigma^2\tr\bP(\bw)\right| \leq \mathbb{E} \max_{1 \leq m \leq M_n} \left|\beps^{\top}\bP_m\beps - \sigma^2\tr\bP_m\right|.
\end{equation*}
Define $\kappa = \mathbb{E} \epsilon_i^4 - 3\sigma^4$. For each $m$, the variance of $\beps^{\top}\bP_m\beps$ can be upper bounded by
\begin{equation*}
  \begin{split}
     \mathbb{E}(\beps^{\top}\bP_m\beps - \sigma^2\tr\bP_m )^2 & = \mathbb{E} (\beps^{\top}\bP_m\beps)^2 - (\sigma^2\tr\bP_m )^2 \\
       & = \sigma^4 [(\tr\bP_m)^2 + 2 \tr\bP_m]+\kappa \tr(\bP_m*\bP_m)- (\sigma^2\tr\bP_m )^2 \\
       & \leq \sigma^4 (k_m^2 + 2 k_m) + \kappa k_m - \sigma^4 k_m^2 \\
       & = (2 \sigma^4 + \kappa) k_m \leq C \sigma^4 k_m,
  \end{split}
\end{equation*}
where the second step follows from Lemma A.2 in \cite{Zhang_2021}, $*$ there denote the Hadamard product of two matrix, and the inequality follows from for the projection matrix $\bP_m$, $\tr(\bP_m*\bP_m) \leq \tr \bP_m$. Next, applying the union bound and Markov's inequality, we obtain the tail probability bound
\begin{equation*}
  \begin{split}
     \mathbb{P}\left( \max_{1 \leq m \leq M_n} |\beps^{\top}\bP_m\beps - \sigma^2\tr\bP_m | >t  \right) & \leq \sum_{m=1}^{M_n}\mathbb{P}\left(  |\beps^{\top}\bP_m\beps - \sigma^2\tr\bP_m | >t  \right) \\
       & \leq \sum_{m=1}^{M_n} \frac{\mathbb{E}(\beps^{\top}\bP_m\beps - \sigma^2\tr\bP_m )^2}{t^2}\\
       & \leq  \frac{C\sigma^4 \sum_{m=1}^{M_n}k_m}{t^2}.
  \end{split}
\end{equation*}
Integrating the tail probability yields
\begin{equation*}
  \begin{split}
     \mathbb{E} \max_{1 \leq m \leq M_n} \left|\beps^{\top}\bP_m\beps - \sigma^2\tr\bP_m\right| & = \int_{0}^{\infty} \mathbb{P}\left( \max_{1 \leq m \leq M_n} |\beps^{\top}\bP_m\beps - \sigma^2\tr\bP_m | >t  \right)dt \\
       & \leq \int_{0}^{\infty} \min \left(1, \frac{C \sigma^4 \sum_{m=1}^{M_n} k_m}{t^2} \right) dt \\
       & \leq \int_{0}^{\sqrt{C \sigma^4 \sum_{m=1}^{M_n} k_m}} 1 dt + \int_{\sqrt{C \sigma^4 \sum_{m=1}^{M_n} k_m}}^{\infty}\frac{C \sigma^4 \sum_{m=1}^{M_n} k_m}{t^2} dt \\
       & \leq C \sigma^2 \sqrt{\sum_{m=1}^{M_n} k_m}.
  \end{split}
\end{equation*}
Therefore, we establish the upper bound
\begin{equation}\label{eq:part12}
  \mathbb{E} \sup_{\bw \in \mathcal{W}} \left|\beps^{\top}\bP(\bw)\beps - \sigma^2\tr\bP(\bw)\right| \leq C \sigma^2 \sqrt{\sum_{m=1}^{M_n} k_m}.
\end{equation}

At last, we upper bound the term $\mathbb{E}|\hat{\sigma}^2 - \sigma^2|\sup_{\bw \in \mathcal{W}}\tr\bP(\bw)$ in (\ref{eq:sa_lu}). It is also observed that $\tr\bP(\bw)$ is a linear function respect to $\bw$. Thus, 
\begin{equation}\label{eq:variance_est_inequality}
  \mathbb{E}\left|\hat{\sigma}^2 - \sigma^2\right|\sup_{\bw \in \mathcal{W}}\tr\bP(\bw) \leq \mathbb{E}\left|\hat{\sigma}^2 - \sigma^2\right| \sup_{1 \leq m \leq M_n}\tr\bP_m \leq \mathbb{E}\left|\hat{\sigma}^2 - \sigma^2\right| \sup_{1 \leq m \leq M_n} k_m. 
\end{equation}

Finally, combining equations (\ref{eq:sa_lu}), (\ref{eq:part11})–(\ref{eq:variance_est_inequality}), we conclude the proof of Proposition~\ref{theo:sharp}.

\subsection{Proof of Theorem~\ref{tho:1}}

Before proceeding with the proof of Theorem~\ref{tho:1}, we state a useful lemma, which has already been established in Section A.2 of \cite{Zhang_2021}.
\begin{lemma}\label{lemma:1}
Let $\mathcal{M} \in \mathbf{M}(M_n)$ be a general candidate model set. Then, there exists a positive constant $C$ such that
  \begin{equation}\label{eq:lem1}
    \mathbb{E}\sup_{\bw \in \mathcal{W}}  \frac{\left\langle \bA(\bw)\bmu, \beps\right\rangle^2}{\left\| \bA(\bw) \bmu\right\|^2}  \leq CM_n,
  \end{equation}
  \begin{equation}\label{eq:lem2}
    \mathbb{E}\sup_{\bw \in \mathcal{W}}\frac{\left\langle \bA^2(\bw)\bmu, \beps\right\rangle^2}{\left\| \bA(\bw) \bmu \right\|^2}\leq CM_n^2,
  \end{equation}
    \begin{equation}\label{eq:lem3}
    \mathbb{E}\sup_{\bw \in \mathcal{W}}\frac{\left[\beps^{\top}\bP(\bw)\beps - \sigma^2\tr\bP(\bw)\right]^2}{\sigma^2\tr\bP^2(\bw)}\leq C M_n,
  \end{equation}
  and
  \begin{equation}\label{eq:lem4}
  \mathbb{E}\sup_{\bw\in\mathcal{W}} \frac{\left[\beps^{\top}\bP^2(\bw)\beps - \sigma^2\tr\bP^2(\bw)\right]^2}{\sigma^2\tr\bP^2(\bw)}\leq C M_n^2,
  \end{equation}
  where $\beps = (\epsilon_1,\ldots, \epsilon_n)^{\top}$ is the random error vector, and $\epsilon_i$ satisfies Assumption~\ref{ass:4_moment}.
\end{lemma}

We now proceed with the proof of the oracle inequality stated in Theorem~\ref{tho:1}. The main theoretical tool adopted in the proof is inspired by the technique developed in \cite{cao2005oracle, cao2006splines}, which is also known as the shifted empirical process method in some model selection literature \citep[see, e.g.,][]{Baraud2000, Wegkamp2003, Lecue2012}. The proof is divided into three steps. In the first step, we establish an intermediate oracle inequality for $\mathbb{E}L_{n}(\hat{\bmu}_{\hat{\bw}_1|\mathcal{M}}, \bmu)$. The remainder term in this inequality involves the plug-in risk $\mathbb{E}R_{n}(\hat{\bmu}_{\hat{\bw}_1|\mathcal{M}}, \bmu)$, whose precise definition is given below. In the second step, we derive an upper bound for $\mathbb{E}R_{n}(\hat{\bmu}_{\hat{\bw}_1|\mathcal{M}}, \bmu)$. Combining these two results yields the desired oracle inequality.

\subsubsection{An intermediate oracle inequality}

For any $0<\gamma<1$, the loss function of the MMA estimator can be decomposed as
\begin{equation}\label{eq:A3}
\begin{split}
   (1-\gamma)n L_{n}(\hat{\bmu}_{\hat{\bw}_1|\mathcal{M}}, \bmu) & =n L_{n}(\hat{\bmu}_{\hat{\bw}_1|\mathcal{M}}, \bmu)-\gamma n L_{n}(\hat{\bmu}_{\hat{\bw}_1|\mathcal{M}}, \bmu)\\
    &= n L_{n}(\hat{\bmu}_{\hat{\bw}_1|\mathcal{M}}, \bmu) - \gamma \left\| \bA(\hat{\bw}_1) \bmu \right\|^2 - \gamma \sigma^2 \operatorname{tr}\bP^2(\hat{\bw}_1)\\
      &+ 2\gamma\left\langle \bA(\hat{\bw}_1)\bmu, \beps \right\rangle - 2\gamma\left\langle \bA^2(\hat{\bw}_1)\bmu, \beps \right\rangle -\gamma \left[\beps^{\top}\bP^2(\hat{\bw}_1)\beps - \sigma^2 \operatorname{tr}\bP^2(\hat{\bw}_1) \right].
\end{split}
\end{equation}
Combining (\ref{eq:A3}) with the first inequality in (\ref{eq:sa_lu}), we get
\begin{equation}\label{eq:A4}
  \begin{split}
     (1-\gamma)n\mathbb{E} L_{n}(\hat{\bmu}_{\hat{\bw}_1|\mathcal{M}}, \bmu)& \leq nR_{n}(\hat{\bmu}_{\bw^*|\mathcal{M}}, \bmu)\\
      &+\mathbb{E}\left\{ (2\gamma-2)\langle \bA(\hat{\bw}_1)\bmu, \beps\rangle - \gamma(1-\gamma) \left\| \bA(\hat{\bw}_1) \bmu \right\|^2\right\}\\
       & +\mathbb{E}\left\{ -2\gamma \langle \bA^2(\hat{\bw}_1)\bmu, \beps\rangle - \gamma^2 \left\| \bA(\hat{\bw}_1) \bmu \right\|^2\right\}\\
       &+\mathbb{E}\left\{2 \left[\beps^{\top}\bP(\hat{\bw}_1)\beps - \sigma^2\operatorname{tr}\bP(\hat{\bw}_1)\right] - \gamma(1-\gamma)\sigma^2\operatorname{tr}\bP^2(\hat{\bw}_1)\right\}\\
       &+\mathbb{E}\left\{-\gamma\left[\beps^{\top}\bP^2(\hat{\bw}_1)\beps - \sigma^2\operatorname{tr}\bP^2(\hat{\bw}_1)\right] - \gamma^2 \sigma^2\operatorname{tr}\bP^2(\hat{\bw}_1)\right\}\\
       & +2\mathbb{E}\left[(\hat{\sigma}^2 - \sigma^2)\tr \bP(\bw^*)\right] - 2\mathbb{E}\left[(\hat{\sigma}^2 - \sigma^2)\tr \bP(\hat{\bw}_1)\right].
  \end{split}
\end{equation}
The task is now to construct the upper bounds for the remainder terms on right side of (\ref{eq:A4}), respectively. 

Note that the first remainder term in (\ref{eq:A4}) is upper bounded by
\begin{equation}\label{eq:A5}
\begin{split}
     &\mathbb{E}\left\{ (2\gamma-2)\langle \bA(\hat{\bw}_1)\bmu, \beps\rangle - \gamma(1-\gamma) \left\| \bA(\hat{\bw}_1) \bmu \right\|^2\right\}\\
     & \leq (2-2\gamma)\mathbb{E}\sup_{\bw \in \mathcal{W}}\left\{ - \langle \bA(\bw)\bmu, \beps\rangle - \frac{\gamma}{2} \left\| \bA(\bw) \bmu \right\|^2\right\}\\
     &\leq (2-2\gamma)\mathbb{E}\sup_{\bw \in \mathcal{W}}\left\{- \langle \bA(\bw)\bmu, \beps\rangle - \frac{\gamma}{2} \left\| \bA(\bw) \bmu \right\|^2 \right\}_{+} \\
     & \leq (2-2\gamma)\mathbb{E}\sup_{\bw \in \mathcal{W}}\left\{- \langle \bA(\bw)\bmu, \beps\rangle 1_{\left\{ -\langle \bA(\bw)\bmu, \beps\rangle \geq \frac{\gamma}{2} \left\| \bA(\bw) \bmu \right\|^2\right\}} \right\}\\
     & \leq (2-2\gamma)\mathbb{E}\sup_{\bw \in \mathcal{W}}\left[ \left| \langle \bA(\bw)\bmu, \beps\rangle \right| \frac{2\left|\langle \bA(\bw)\bmu, \beps\rangle \right|}{\gamma \left\| \bA(\bw) \bmu\right\|^2}  \right]\\
      &= \frac{4-4\gamma}{\gamma} \mathbb{E}\sup_{\bw \in \mathcal{W}}  \frac{\langle \bA(\bw)\bmu, \beps\rangle^2}{\left\| \bA(\bw) \bmu\right\|^2}  \leq\frac{C(4-4\gamma)M_n}{\gamma},
\end{split}
\end{equation}
where the forth step is due to $\eta 1\left\{\eta\geq x\right\}\leq |\eta| |\eta / x|$, and the last step follows from Lemma~\ref{lemma:1}.

Similarly, the second remainder term in (\ref{eq:A4}) is upper bounded by
\begin{equation}\label{eq:A6}
  \begin{split}
      & \mathbb{E}\left\{ -2\gamma \langle \bA^2(\hat{\bw}_1)\bmu, \beps\rangle - \gamma^2 \left\| \bA(\hat{\bw}_1) \bmu \right\|^2\right\} \\
      & \leq 2\gamma \mathbb{E}\sup_{\bw \in \mathcal{W}}\left\{ - \langle \bA^2(\bw)\bmu, \beps\rangle - \frac{\gamma}{2} \left\| \bA(\bw) \bmu \right\|^2\right\}\\
      & \leq 2\gamma \mathbb{E}\sup_{\bw \in \mathcal{W}}\left\{ - \langle \bA^2(\bw)\bmu, \beps\rangle - \frac{\gamma}{2} \left\| \bA(\bw) \bmu \right\|^2\right\}_+ \\
      & \leq 2\gamma \mathbb{E}\sup_{\bw \in \mathcal{W}}\left\{ - \langle \bA^2(\bw)\bmu, \beps\rangle 1_{\left\{- \langle \bA^2(\bw)\bmu, \beps\rangle\geq  \frac{\gamma}{2} \left\| \bA(\bw) \bmu \right\|^2\right\}}\right\} \\
      & \leq 2\gamma \mathbb{E}\sup_{\bw \in \mathcal{W}}\left[ \left| \langle \bA^2(\bw)\bmu, \beps\rangle \right| \frac{2\left|\langle \bA^2(\bw)\bmu, \beps\rangle\right|}{\gamma\left\| \bA(\bw) \bmu \right\|^2}\right]\\
      & \leq 4\mathbb{E}\sup_{\bw \in \mathcal{W}}\frac{\langle \bA^2(\bw)\bmu, \beps\rangle^2}{\left\| \bA(\bw) \bmu \right\|^2}\leq4CM_n^2.\\
  \end{split}
\end{equation}
The third remainder term is upper bounded by
\begin{equation}\label{eq:A66}
\begin{split}
    &   \mathbb{E}\left\{2 \left[\beps^{\top}\bP(\hat{\bw}_1)\beps - \sigma^2\operatorname{tr}\bP(\hat{\bw}_1)\right] - \gamma(1-\gamma)\sigma^2\operatorname{tr}\bP^2(\hat{\bw}_1)\right\}\\
    & \leq 2 \mathbb{E}\sup_{\bw \in \mathcal{W}}\left\{ \left[\beps^{\top}\bP(\bw)\beps - \sigma^2\operatorname{tr}\bP(\bw)\right] - \frac{\gamma(1-\gamma)}{2}\sigma^2\operatorname{tr}\bP^2(\bw) \right\}\\
        & \leq 2 \mathbb{E}\sup_{\bw \in \mathcal{W}}\left\{ \left[\beps^{\top}\bP(\bw)\beps - \sigma^2\operatorname{tr}\bP(\bw)\right] - \frac{\gamma(1-\gamma)}{2}\sigma^2\operatorname{tr}\bP^2(\bw) \right\}_+\\
        & \leq 2 \mathbb{E}\sup_{\bw \in \mathcal{W}}\left\{ \left[\beps^{\top}\bP(\bw)\beps - \sigma^2\operatorname{tr}\bP(\bw)\right]1_{\{ [\beps^{\top}\bP(\bw)\beps - \sigma^2\operatorname{tr}\bP(\bw)] \geq\frac{\gamma(1-\gamma)}{2}\sigma^2\operatorname{tr}\bP^2(\bw) \}}  \right\}\\
     &\leq \frac{4}{\gamma(1-\gamma)}\mathbb{E}\sup_{\bw \in \mathcal{W}}\frac{\left[\beps^{\top}\bP(\bw)\beps - \sigma^2\tr\bP(\bw)\right]^2}{\sigma^2\tr\bP^2(\bw)}\leq \frac{4C M_n}{\gamma(1-\gamma)}.
\end{split}
\end{equation}
The forth remainder term in (\ref{eq:A4}) can be upper bounded by
\begin{equation}\label{eq:A661}
\begin{split}
    &   \mathbb{E}\left\{-\gamma\left[\beps^{\top}\bP^2(\hat{\bw}_1)\beps - \sigma^2\operatorname{tr}\bP^2(\hat{\bw}_1)\right] - \gamma^2 \sigma^2\operatorname{tr}\bP^2(\hat{\bw}_1)\right\} \\
    & \leq \gamma \mathbb{E}\sup_{\bw \in \mathcal{W}}\left\{-\left[\beps^{\top}\bP^2(\bw)\beps - \sigma^2\operatorname{tr}\bP^2(\bw)\right] - \gamma \sigma^2\operatorname{tr}\bP^2(\bw)\right\} \\
    & \leq \gamma \mathbb{E}\sup_{\bw \in \mathcal{W}}\left\{-\left[\beps^{\top}\bP^2(\bw)\beps - \sigma^2\operatorname{tr}\bP^2(\bw)\right] - \gamma \sigma^2\operatorname{tr}\bP^2(\bw)\right\}_+\\
    & \leq \gamma \mathbb{E}\sup_{\bw \in \mathcal{W}}\left\{-\left[\beps^{\top}\bP^2(\bw)\beps - \sigma^2\operatorname{tr}\bP^2(\bw)\right]1_{\{ -[\beps^{\top}\bP^2(\bw)\beps - \sigma^2\operatorname{tr}\bP^2(\bw)]\geq\gamma \sigma^2\operatorname{tr}\bP^2(\bw)\}} \right\}\\
    & \leq \mathbb{E}\sup_{\bw\in\mathcal{W}} \frac{\left[\beps^{\top}\bP^2(\bw)\beps - \sigma^2\operatorname{tr}\bP^2(\bw)\right]^2}{\sigma^2\operatorname{tr}\bP^2(\bw)}\leq CM_n^2.
\end{split}
\end{equation}

The last line of (\ref{eq:A4}) involves two remainder terms related to $\hat{\sigma}^2$. The first term $\mathbb{E}[(\hat{\sigma}^2 - \sigma^2)\tr \bP(\bw^*)]$ can be upper bounded by 
\begin{equation}\label{eq:chuntian3}
  \begin{split}
     &\mathbb{E}\left[(\hat{\sigma}^2 - \sigma^2)\tr \bP(\bw^*)\right] \leq \mathbb{E} \left| \hat{\sigma}^2 - \sigma^2 \right| \tr \bP(\bw^*) \leq \sqrt{\mathbb{E}\left( \hat{\sigma}^2 - \sigma^2 \right)^2}\tr \bP(\bw^*)\\
     &\leq \sqrt{\mathbb{E}\left( \hat{\sigma}^2 - \sigma^2 \right)^2} \sqrt{n \tr \bP^2(\bw^*)} \leq \sqrt{\mathbb{E}\left( \hat{\sigma}^2 - \sigma^2 \right)^2} \sqrt{\frac{n^2  R_{n}(\hat{\bmu}_{\bw^*|\mathcal{M}}, \bmu)}{\sigma^2}}, 
  \end{split}
\end{equation}
where the second inequality follows from Jensen's inequality. For the third inequality, note that for any symmetric matrix $\bP$, $\tr \bP = \sum_{i=1}^{n} \lambda_i$, where $\lambda_i$ are the eigenvalues of $\bP$. By the Cauchy--Schwarz inequality, 
\begin{equation}\label{eq:eigen_cs}
  \sum_{i=1}^{n} \lambda_i \leq \sqrt{n\sum_{i=1}^{n}\lambda_i^2} = \sqrt{n\tr \bP^2}. 
\end{equation}
Finally, the last inequality in (\ref{eq:chuntian3}) follows from (\ref{eq:leavv1}), which gives that $nR_{n}(\hat{\bmu}_{\bw^*|\mathcal{M}}, \bmu) = \| \bA(\bw^*) \bmu \|^2 + \sigma^2\tr \bP^2(\bw^*)$. Similarly, the second term in the last line of (\ref{eq:A4}) can be upper bounded by  
\begin{equation}\label{eq:xiatian2026}
  \begin{split}
     &-\mathbb{E}\left[(\hat{\sigma}^2 - \sigma^2)\tr \bP(\hat{\bw}_1)\right] \leq \sqrt{\mathbb{E}\left(\hat{\sigma}^2-\sigma^2\right)^2} \sqrt{\mathbb{E}\left[\operatorname{tr} \mathbf{P}(\hat{\mathbf{w}}_1)\right]^2} .\\
        & \leq \sqrt{\mathbb{E}\left(\hat{\sigma}^2-\sigma^2\right)^2} \sqrt{n\mathbb{E} \tr \bP^2(\hat{\bw}_1)} \leq \sqrt{\mathbb{E}\left( \hat{\sigma}^2 - \sigma^2 \right)^2} \sqrt{\frac{n^2  \mathbb{E} R_{n}(\hat{\bmu}_{\hat{\bw}_1|\mathcal{M}}, \bmu)}{\sigma^2}}, 
  \end{split}
\end{equation}
where the first inequality follows from the Cauchy-Schwarz inequality, and the second inequality follows by the same argument as in (\ref{eq:eigen_cs}). The last inequality follows from the identity $nR_{n}(\hat{\bmu}_{\hat{\bw}_1|\mathcal{M}}, \bmu) \triangleq \| \bA(\hat{\bw}_1) \bmu \|^2 + \sigma^2\tr \bP^2(\hat{\bw}_1)$, which is obtained by plugging the MMA weight $\hat{\bw}_1$ into the right-hand side of (\ref{eq:leavv1}). We refer to $R_{n}(\hat{\bmu}_{\hat{\bw}_1|\mathcal{M}}, \bmu)$ as the plug-in risk thereafter. 

Substituting (\ref{eq:A5})--(\ref{eq:chuntian3}) and (\ref{eq:xiatian2026}) into (\ref{eq:A4}), we obtain that for any $0< \gamma < 1$,
\begin{equation}\label{eq:A15}
\begin{split}
   (1-\gamma)\mathbb{E} L_{n}(\hat{\bmu}_{\hat{\bw}_1|\mathcal{M}}, \bmu)& \leq R_{n}(\hat{\bmu}_{\bw^*|\mathcal{M}}, \bmu) +\frac{C(1-\gamma)M_n}{\gamma n}+\frac{CM_n}{\gamma(1-\gamma)n}+\frac{CM_n^2}{n} \\
   & + C\sqrt{\mathbb{E}(\hat{\sigma}^2 - \sigma^2)^2}\left[  \sqrt{R_{n}(\hat{\bmu}_{\bw^*|\mathcal{M}}, \bmu)} + \sqrt{\mathbb{E} R_{n}(\hat{\bmu}_{\hat{\bw}_1|\mathcal{M}}, \bmu)} \right].
\end{split}
\end{equation}
Note that the right-hand side involves the term $\sqrt{\mathbb{E} R_{n}(\hat{\bmu}_{\hat{\bw}_1|\mathcal{M}}, \bmu)}$, which arises from estimating $\sigma^2$. This term will be bounded in the next step.

\subsubsection{Upper bounding the plug-in risk} 

Similar to (\ref{eq:A3}), for any $0<\gamma<1$, the plug-in risk function $R_{n}(\hat{\bmu}_{\hat{\bw}_1|\mathcal{M}}, \bmu)$ of the MMA weights $\hat{\bw}_1$ can be decomposed as 
\begin{equation}\label{eq:pi_risk}
\begin{split}
   (1-\gamma)n R_{n}(\hat{\bmu}_{\hat{\bw}_1|\mathcal{M}}, \bmu) & =n R_{n}(\hat{\bmu}_{\hat{\bw}_1|\mathcal{M}}, \bmu)-\gamma n R_{n}(\hat{\bmu}_{\hat{\bw}_1|\mathcal{M}}, \bmu)\\
    &= n R_{n}(\hat{\bmu}_{\hat{\bw}_1|\mathcal{M}}, \bmu) - \gamma \left\| \bA(\hat{\bw}_1) \bmu \right\|^2 - \gamma \sigma^2 \operatorname{tr}\bP^2(\hat{\bw}_1),\\
\end{split}
\end{equation}
where the second equality follows from the definition $n R_{n}(\hat{\bmu}_{\hat{\bw}_1|\mathcal{M}}, \bmu) \triangleq \| \bA(\hat{\bw}_1) \bmu \|^2 + \sigma^2\tr \bP^2(\hat{\bw}_1)$, which is defined by plugging $\hat{\bw}_1$ in the right hand side of (\ref{eq:leavv1}) directly. Note that $R_{n}(\hat{\bmu}_{\hat{\bw}_1|\mathcal{M}}, \bmu)$ is still a random quantity. Similar to (\ref{eq:A2}), the plug-in risk function and the MMA criterion have the following relationship 
\begin{equation}\label{eq:leavv111}
\begin{split}
n R_n\left(\hat{\bmu}_{\hat{\bw}_1 | \mathcal{M}}, \bmu\right)& = n C_n\left(\hat{\bw}_1 | \mathcal{M}, \lambda_1\right)-\|\beps\|^2  -2\left\langle\bA^2(\hat{\bw}_1) \bmu, \beps\right\rangle \\
& +2\left[\beps^{\top} \bP(\hat{\bw}_1) \beps-\sigma^2 \tr \bP(\hat{\bw}_1)\right]  -\left[\beps^{\top} \bP^2(\hat{\bw}_1) \beps-\sigma^2 \tr \mathbf{P}^2(\hat{\bw}_1)\right] \\
& +2(\sigma^2-\hat{\sigma}^2) \tr \mathbf{P}(\hat{\mathbf{w}}_1) .
\end{split}
\end{equation}
Similar to (\ref{eq:A4}), we combine (\ref{eq:pi_risk})--(\ref{eq:leavv111}) with (\ref{eq:A1}) and get 
\begin{equation}\label{eq:summer2026_1}
  \begin{split}
     (1-\gamma) n \mathbb{E} R_n\left(\hat{\bmu}_{\hat{\bw}_1 | \mathcal{M}}, \bmu\right) & \leq n R_n\left(\hat{\bmu}_{\bw^* | \mathcal{M}}, \bmu\right) \\
       &+ \mathbb{E}\left\{-2\left\langle\mathbf{A}^2(\hat{\bw}_1) \boldsymbol{\mu}, \beps\right\rangle-\gamma\left\|\mathbf{A}(\hat{\bw}_1) \boldsymbol{\mu}\right\|^2\right\}\\
       & + \mathbb{E}\left\{2\left[\beps^{\top} \mathbf{P}(\hat{\bw}_1) \beps-\sigma^2 \operatorname{tr} \mathbf{P}(\hat{\bw}_1)\right]-\frac{\gamma}{2} \sigma^2 \operatorname{tr} \mathbf{P}^2(\hat{\bw}_1)\right\}\\
       & +\mathbb{E}\left\{-\left[\beps^{\top} \mathbf{P}^2(\hat{\bw}_1) \beps-\sigma^2 \operatorname{tr} \mathbf{P}^2(\hat{\bw}_1)\right]-\frac{\gamma}{2} \sigma^2 \operatorname{tr} \mathbf{P}^2(\hat{\bw}_1)\right\}\\
       & +2\mathbb{E}\left[\left( \hat{\sigma}^2-\sigma^2\right) \operatorname{tr} \mathbf{P}(\bw^*)\right]-2 \mathbb{E}\left[\left(\hat{\sigma}^2-\sigma^2\right) \operatorname{tr} \mathbf{P}(\hat{\bw}_1)\right]. 
  \end{split}
\end{equation}
It remains to upper bound the remainder terms in (\ref{eq:summer2026_1}), respectively. 

For the first remainder term in (\ref{eq:summer2026_1}), we have 
\begin{equation}\label{eq:summer2026_2}
  \begin{split}
       & \mathbb{E}\left\{-2\left\langle\mathbf{A}^2(\hat{\bw}_1) \boldsymbol{\mu}, \beps\right\rangle-\gamma\left\|\mathbf{A}(\hat{\bw}_1) \boldsymbol{\mu}\right\|^2\right\} \\
       & \leq \mathbb{E}\sup_{\bw \in \mathcal{W}}\left\{-2\left\langle\mathbf{A}^2(\bw) \boldsymbol{\mu}, \beps\right\rangle-\gamma\left\|\mathbf{A}(\bw) \boldsymbol{\mu}\right\|^2\right\}\\
       & \leq \mathbb{E}\sup_{\bw \in \mathcal{W}}\left\{-2\left\langle\mathbf{A}^2(\bw) \boldsymbol{\mu}, \beps\right\rangle-\gamma\left\|\mathbf{A}(\bw) \boldsymbol{\mu}\right\|^2\right\}_{+}\\
       & \leq \mathbb{E}\sup_{\bw \in \mathcal{W}} \left\{ -2\left\langle\mathbf{A}^2(\bw) \boldsymbol{\mu}, \beps\right\rangle 1_{\{ -2\left\langle\mathbf{A}^2(\bw) \boldsymbol{\mu}, \beps\right\rangle \geq \gamma\left\|\mathbf{A}(\bw) \boldsymbol{\mu}\right\|^2\}} \right\}\\
       & \leq \mathbb{E}\sup_{\bw \in \mathcal{W}} \left\{ \left|2\left\langle\mathbf{A}^2(\bw) \boldsymbol{\mu}, \beps\right\rangle \right| \frac{\left|2\left\langle\mathbf{A}^2(\bw) \boldsymbol{\mu}, \beps\right\rangle \right|}{\gamma\left\|\mathbf{A}(\bw) \boldsymbol{\mu}\right\|^2} \right\}\\
       & \leq \frac{4}{\gamma} \mathbb{E}\sup_{\bw \in \mathcal{W}}  \frac{\left\langle\mathbf{A}^2(\bw) \boldsymbol{\mu}, \beps\right\rangle^2}{\left\|\mathbf{A}(\bw) \boldsymbol{\mu}\right\|^2} \leq\frac{4CM_n^2}{\gamma}, 
  \end{split}
\end{equation}
where the last inequality follows from (\ref{eq:lem2}) in Lemma~\ref{lemma:1}. 

The second remainder term in (\ref{eq:summer2026_1}) is upper bounded by
\begin{equation}\label{eq:summer2026_3}
\begin{split}
    &   \mathbb{E}\left\{2\left[\beps^{\top} \mathbf{P}(\hat{\bw}_1) \beps-\sigma^2 \operatorname{tr} \mathbf{P}(\hat{\bw}_1)\right]-\frac{\gamma}{2} \sigma^2 \operatorname{tr} \mathbf{P}^2(\hat{\bw}_1)\right\}\\
    & \leq 2 \mathbb{E}\sup_{\bw \in \mathcal{W}}\left\{ \left[\beps^{\top}\bP(\bw)\beps - \sigma^2\operatorname{tr}\bP(\bw)\right] - \frac{\gamma}{4}\sigma^2\operatorname{tr}\bP^2(\bw) \right\}\\
        & \leq 2 \mathbb{E}\sup_{\bw \in \mathcal{W}}\left\{ \left[\beps^{\top}\bP(\bw)\beps - \sigma^2\operatorname{tr}\bP(\bw)\right] - \frac{\gamma}{4}\sigma^2\operatorname{tr}\bP^2(\bw) \right\}_+\\
        & \leq 2 \mathbb{E}\sup_{\bw \in \mathcal{W}}\left\{ \left[\beps^{\top}\bP(\bw)\beps - \sigma^2\operatorname{tr}\bP(\bw)\right]1_{\{ [\beps^{\top}\bP(\bw)\beps - \sigma^2\operatorname{tr}\bP(\bw)] \geq\frac{\gamma}{4}\sigma^2\operatorname{tr}\bP^2(\bw) \}}  \right\}\\
     &\leq \frac{8}{\gamma}\mathbb{E}\sup_{\bw \in \mathcal{W}}\frac{\left[\beps^{\top}\bP(\bw)\beps - \sigma^2\tr\bP(\bw)\right]^2}{\sigma^2\tr\bP^2(\bw)}\leq \frac{8C M_n}{\gamma}, 
\end{split}
\end{equation}
where the last inequality follows from (\ref{eq:lem3}) in Lemma~\ref{lemma:1}. 

The third remainder term in (\ref{eq:summer2026_1}) is upper bounded by 
\begin{equation}\label{eq:summer2026_4}
\begin{split}
    &   \mathbb{E}\left\{-\left[\beps^{\top} \mathbf{P}^2(\hat{\bw}_1) \beps-\sigma^2 \operatorname{tr} \mathbf{P}^2(\hat{\bw}_1)\right]-\frac{\gamma}{2} \sigma^2 \operatorname{tr} \mathbf{P}^2(\hat{\bw}_1)\right\}\\
    & \leq \mathbb{E}\sup_{\bw \in \mathcal{W}}\left\{ -\left[\beps^{\top}\bP^2(\bw)\beps - \sigma^2\operatorname{tr}\bP^2(\bw)\right] - \frac{\gamma}{2}\sigma^2\operatorname{tr}\bP^2(\bw) \right\}\\
        & \leq \mathbb{E}\sup_{\bw \in \mathcal{W}}\left\{ -\left[\beps^{\top}\bP^2(\bw)\beps - \sigma^2\operatorname{tr}\bP^2(\bw)\right] - \frac{\gamma}{2}\sigma^2\operatorname{tr}\bP^2(\bw) \right\}_+\\
        & \leq  \mathbb{E}\sup_{\bw \in \mathcal{W}}\left\{ -\left[\beps^{\top}\bP^2(\bw)\beps - \sigma^2\operatorname{tr}\bP^2(\bw)\right]1_{\{ -[\beps^{\top}\bP^2(\bw)\beps - \sigma^2\operatorname{tr}\bP^2(\bw)] \geq\frac{\gamma}{2}\sigma^2\operatorname{tr}\bP^2(\bw) \}}  \right\}\\
     &\leq \frac{2}{\gamma}\mathbb{E}\sup_{\bw \in \mathcal{W}}\frac{\left[\beps^{\top}\bP^2(\bw)\beps - \sigma^2\tr\bP^2(\bw)\right]^2}{\sigma^2\tr\bP^2(\bw)}\leq \frac{2C M_n^2}{\gamma}, 
\end{split}
\end{equation}
where the last inequality follows from (\ref{eq:lem4}) in Lemma~\ref{lemma:1}. 

Using the results in (\ref{eq:chuntian3}) and (\ref{eq:xiatian2026}), we see that the term in the last line of (\ref{eq:summer2026_1}) is upper bounded by 
\begin{equation}\label{eq:summer2026_5}
\begin{split}
     & 2\mathbb{E}\left[\left( \hat{\sigma}^2-\sigma^2\right) \operatorname{tr} \mathbf{P}(\bw^*)\right]-2 \mathbb{E}\left[\left(\hat{\sigma}^2-\sigma^2\right) \operatorname{tr} \mathbf{P}(\hat{\bw}_1)\right] \\
     & \leq C\sqrt{\mathbb{E}(\hat{\sigma}^2 - \sigma^2)^2} \left[ \sqrt{n^2 R_{n}(\hat{\bmu}_{\bw^*|\mathcal{M}}, \bmu)} +  \sqrt{n^2 \mathbb{E} R_{n}(\hat{\bmu}_{\hat{\bw}_1|\mathcal{M}}, \bmu)} \right]. 
\end{split}
\end{equation}

Inserting the upper bounds in (\ref{eq:summer2026_2})--(\ref{eq:summer2026_5}) into the right-hand side of (\ref{eq:summer2026_1}), we obtain 
\begin{equation*}
\begin{split}
    (1-\gamma) \mathbb{E} R_n\left(\hat{\bmu}_{\hat{\bw}_1 | \mathcal{M}}, \bmu\right) & \leq R_n\left(\hat{\bmu}_{\bw^* | \mathcal{M}}, \bmu\right) + \frac{CM_n}{n\gamma} + \frac{CM_n^2}{n\gamma} \\
     & +  C\sqrt{\mathbb{E}(\hat{\sigma}^2 - \sigma^2)^2}\left[  \sqrt{R_{n}(\hat{\bmu}_{\bw^*|\mathcal{M}}, \bmu)} + \sqrt{\mathbb{E} R_{n}(\hat{\bmu}_{\hat{\bw}_1|\mathcal{M}}, \bmu)} \right]. 
\end{split}
\end{equation*}
By setting $\gamma = 1/2$, the above bound implies that 
\begin{equation}\label{eq:summer2026_6}
  \begin{split}
     \mathbb{E} R_n\left(\hat{\bmu}_{\hat{\bw}_1 | \mathcal{M}}, \bmu\right) & \leq 2R_n\left(\hat{\bmu}_{\bw^* | \mathcal{M}}, \bmu\right) + \frac{CM_n^2}{n} \\
       & +  C\sqrt{\mathbb{E}(\hat{\sigma}^2 - \sigma^2)^2}\left[  \sqrt{R_{n}(\hat{\bmu}_{\bw^*|\mathcal{M}}, \bmu)} + \sqrt{\mathbb{E} R_{n}(\hat{\bmu}_{\hat{\bw}_1|\mathcal{M}}, \bmu)} \right].  
  \end{split}
\end{equation}
Define $X \triangleq \sqrt{\mathbb{E} R_{n}(\hat{\bmu}_{\hat{\bw}_1|\mathcal{M}}, \bmu)}$, $A \triangleq 2R_n(\hat{\bmu}_{\bw^* | \mathcal{M}}, \bmu) + \frac{CM_n^2}{n} + C\sqrt{\mathbb{E}(\hat{\sigma}^2 - \sigma^2)^2} \sqrt{R_{n}(\hat{\bmu}_{\bw^*|\mathcal{M}}, \bmu)}$, $B \triangleq C\sqrt{\mathbb{E}(\hat{\sigma}^2 - \sigma^2)^2}$. Then the inequality (\ref{eq:summer2026_6}) can be written as $X^2 - BX - A \leq 0$, which also implies that $X \leq \frac{B+\sqrt{B^2+4 A}}{2}$. And consequently, 
\begin{equation}\label{eq:summer2026_7}
  \sqrt{\mathbb{E} R_{n}(\hat{\bmu}_{\hat{\bw}_1|\mathcal{M}}, \bmu)} \leq C(B + \sqrt{A}) \leq C \left[ \sqrt{ R_n(\hat{\bmu}_{\bw^* | \mathcal{M}}, \bmu)} + \sqrt{\frac{M_n^2}{n}} +  \sqrt{\mathbb{E}(\hat{\sigma}^2 - \sigma^2)^2} \right]. 
\end{equation}

\subsubsection{Completing the proof}

Inserting the risk bound (\ref{eq:summer2026_7}) into the right side of (\ref{eq:A15}) yields 
\begin{equation*}
  \begin{split}
   (1-\gamma)\mathbb{E} L_{n}(\hat{\bmu}_{\hat{\bw}_1|\mathcal{M}}, \bmu)& \leq R_{n}(\hat{\bmu}_{\bw^*|\mathcal{M}}, \bmu) +\frac{C(1-\gamma)M_n}{\gamma n}+\frac{CM_n}{\gamma(1-\gamma)n}+\frac{CM_n^2}{n} \\
   & + C\sqrt{\mathbb{E}(\hat{\sigma}^2 - \sigma^2)^2}\left[  \sqrt{R_{n}(\hat{\bmu}_{\bw^*|\mathcal{M}}, \bmu)} + \sqrt{\frac{M_n^2}{n}} +  \sqrt{\mathbb{E}(\hat{\sigma}^2 - \sigma^2)^2} \right]\\
   & \leq R_{n}(\hat{\bmu}_{\bw^*|\mathcal{M}}, \bmu) +\frac{C(1-\gamma)M_n}{\gamma n}+\frac{CM_n}{\gamma(1-\gamma)n}+\frac{CM_n^2}{n} \\ 
   & + C \,\mathbb{E}(\hat{\sigma}^2 - \sigma^2)^2 + C\sqrt{\mathbb{E}(\hat{\sigma}^2 - \sigma^2)^2} \sqrt{R_{n}(\hat{\bmu}_{\bw^*|\mathcal{M}}, \bmu)}. 
\end{split}
\end{equation*}

Using the change of variable $\delta = \frac{\gamma}{1-\gamma}$, we have
\begin{equation*}
  \begin{split}
     \mathbb{E} L_{n}(\hat{\bmu}_{\hat{\bw}_1|\mathcal{M}}, \bmu) & \leq (1+\delta) R_{n}(\hat{\bmu}_{\bw^*|\mathcal{M}}, \bmu) + \frac{C(1+\delta)M_n}{\delta n}+ \frac{C(1+\delta)^3M_n}{\delta n} + \frac{C(1+\delta)M_n^2}{n}\\
       &+ C (1+\delta)\mathbb{E}(\hat{\sigma}^2 - \sigma^2)^2 + C\sqrt{\frac{(1+\delta)^2\mathbb{E}(\hat{\sigma}^2 - \sigma^2)^2}{\delta}} \sqrt{\delta R_{n}(\hat{\bmu}_{\bw^*|\mathcal{M}}, \bmu)} \\
       & \leq (1+\delta) R_{n}(\hat{\bmu}_{\bw^*|\mathcal{M}}, \bmu) + \frac{C(1+\delta)^3M_n}{\delta n}  + \frac{C(1+\delta)M_n^2}{n}\\
        &+ C (1+\delta)\mathbb{E}(\hat{\sigma}^2 - \sigma^2)^2+ \frac{C(1+\delta)^2\mathbb{E}(\hat{\sigma}^2 - \sigma^2)^2}{\delta} + C\delta R_{n}(\hat{\bmu}_{\bw^*|\mathcal{M}}, \bmu),
  \end{split}
\end{equation*}
which completes the proof of Theorem~\ref{tho:1}.

\subsection{Proof of Corollary~\ref{cor:aop}}\label{sec:proof_add_nested}

Dividing both sides of the oracle inequality in Proposition~\ref{theo:sharp} yields 
\begin{equation*}
  \begin{split}
     \frac{\mathbb{E}L_n(\hat{\bmu}_{\hat{\bw}_1|\mathcal{M}}, \bmu )}{R_{n}(\hat{\bmu}_{\bw^*|\mathcal{M}}, \bmu)} & \leq 1 +\frac{C[\sum_{m=1}^{M_n}\| (\bI - \bP_m)\bmu \|^2]^{1/2}+[\sum_{m=1}^{M_n} k_m]^{1/2}}{nR_{n}(\hat{\bmu}_{\bw^*|\mathcal{M}}, \bmu)}  + \frac{C\mathbb{E}|\hat{\sigma}^2 - \sigma^2|\max_{1 \leq m \leq M_n}k_m}{nR_{n}(\hat{\bmu}_{\bw^*|\mathcal{M}}, \bmu)}\\
     & \leq 1 + \frac{C[\sum_{m=1}^{M_n}\| (\bI - \bP_m)\bmu \|^2+ \sum_{m=1}^{M_n} k_m]^{1/2}}{nR_{n}(\hat{\bmu}_{\bw^*|\mathcal{M}}, \bmu)} + \frac{C\mathbb{E}|\hat{\sigma}^2 - \sigma^2|\max_{1 \leq m \leq M_n}k_m}{nR_{n}(\hat{\bmu}_{\bw^*|\mathcal{M}}, \bmu)}\\
     & \leq 1 + \frac{C[ \sum_{m=1}^{M_n} (\| (\bI - \bP_m)\bmu \|^2 + \sigma^2 k_m) ]^{1/2}}{nR_{n}(\hat{\bmu}_{\bw^*|\mathcal{M}}, \bmu)} + \frac{C\mathbb{E}|\hat{\sigma}^2 - \sigma^2|\max_{1 \leq m \leq M_n}k_m}{nR_{n}(\hat{\bmu}_{\bw^*|\mathcal{M}}, \bmu)}, 
  \end{split}
\end{equation*}
where the second inequality follows from $A^{1/2}+B^{1/2}\leq \sqrt{2}(A+B)^{1/2}$ for $A\geq 0$ and $B\geq 0$. Then, condition (\ref{eq:condition_our}) implies (\ref{eq:aop_our}). 

Dividing both sides of the oracle inequality in Theorem~\ref{tho:1} gives 
\begin{equation*}
\begin{split}
   \frac{\mathbb{E}L_n(\hat{\bmu}_{\hat{\bw}_1|\mathcal{M}}, \bmu )}{R_{n}(\hat{\bmu}_{\bw^*|\mathcal{M}}, \bmu)}&\leq 1+C\delta + \frac{C(1+\delta)^3M_n}{\delta nR_{n}(\hat{\bmu}_{\bw^*|\mathcal{M}}, \bmu)}  + \frac{C(1+\delta)M_n^2}{nR_{n}(\hat{\bmu}_{\bw^*|\mathcal{M}}, \bmu)}+  \frac{C(1+\delta)^2\mathbb{E}(\hat{\sigma}^2 - \sigma^2)^2}{\delta R_{n}(\hat{\bmu}_{\bw^*|\mathcal{M}}, \bmu)}. \\
\end{split}
\end{equation*}
Due to the arbitrariness of $\delta>0$, we can choose $\delta = \delta_n \to 0$ as $n \to \infty$ such that $ (1+ \delta_n)^3/ \delta_n = o(nR_{n}(\hat{\bmu}_{\bw^*|\mathcal{M}}, \bmu)/M_n)$. Such a sequence $\delta_n$ always exists, since $nR_{n}(\hat{\bmu}_{\bw^*|\mathcal{M}}, \bmu)/M_n \geq nR_{n}(\hat{\bmu}_{\bw^*|\mathcal{M}}, \bmu)/M_n^2 \to \infty$ under condition (\ref{eq:variance_cond}). With this choice of $\delta_n$, we have 
\begin{equation*}
  C\delta + \frac{C(1+\delta)^3M_n}{\delta nR_{n}(\hat{\bmu}_{\bw^*|\mathcal{M}}, \bmu)}  + \frac{C(1+\delta)M_n^2}{nR_{n}(\hat{\bmu}_{\bw^*|\mathcal{M}}, \bmu)} \to 0, 
\end{equation*}
and 
\begin{equation*}
  \frac{C(1+\delta)^2\mathbb{E}(\hat{\sigma}^2 - \sigma^2)^2}{\delta R_{n}(\hat{\bmu}_{\bw^*|\mathcal{M}}, \bmu)} =o\left[ \frac{nR_{n}(\hat{\bmu}_{\bw^*|\mathcal{M}}, \bmu)}{M_n}  \frac{\mathbb{E}(\hat{\sigma}^2 - \sigma^2)^2}{ R_{n}(\hat{\bmu}_{\bw^*|\mathcal{M}}, \bmu)} \right] = O\left[ \frac{n \mathbb{E}(\hat{\sigma}^2 - \sigma^2)^2}{M_n} \right] \to 0, 
\end{equation*}
where the last approximation follows from (\ref{eq:variance_cond}).

\subsection{Proof of the results in Section~\ref{sec:nested}}\label{sec:a:proof_nested}

\subsubsection{Preliminaries}\label{sec:pre_nested}

Given the complete basis $\{\bpsi_1, \ldots, \bpsi_p\}$ satisfying (\ref{eq:complete}), the candidate least squares estimator (\ref{eq:estimator}) admits the following spectral representation. The coefficient vector $\btheta \triangleq (\theta_1, \ldots, \theta_p)^{\top}$ is called the transform of $\bmu$ and is an isometry of $\bmu$ in $\mathbb{R}^p$. Define the empirical coefficients $\tilde{\theta}_j \triangleq  \by^{\top}\bpsi_j/n $ and the empirical random error terms $e_j \triangleq  \beps^{\top} \bpsi_j /n$. Accordingly, the vectors $\tilde{\btheta} \triangleq (\tilde{\theta}_1,\ldots,\tilde{\theta}_p)^{\top}$ and $\be \triangleq (e_1,\ldots,e_p)^{\top}$ are the transforms of $\by$ and $\beps$, respectively. The estimator (\ref{eq:estimator}) takes the form
\begin{equation}\label{eq:estimator_2}
 \hat{\bmu}_{\mathcal{I}} = \sum_{j \in {\mathcal{I}}} n^{-1}\by^{\top}\bpsi_j \bpsi_j = \sum_{j \in {\mathcal{I}}} \tilde{\theta}_j \bpsi_j.
\end{equation}

In the nested setup, the MA estimators based on $\mathcal{M}_{AN}$ and $\mathcal{M}_{G}$ can be expressed as
\begin{equation}\label{eq:est_all_nested}
  \begin{split}
     \hat{\bmu}_{\bw|\mathcal{M}_{AN}} & = \sum_{k=1}^{p}w_k \sum_{j =1}^{k} \tilde{\theta}_j \bpsi_j = \sum_{j=1}^{p}\lambda_j\tilde{\theta}_j \bpsi_j,\\
  \end{split}
\end{equation}
where $\lambda_j = \sum_{k=j}^{p}w_k$, and
\begin{equation}\label{eq:est_group}
  \begin{split}
     \hat{\bmu}_{\bw'|\mathcal{M}_{G}} & = \sum_{t=1}^{T_n}w_t' \sum_{l =1}^{j_{t}} \tilde{\theta}_l \bpsi_l = \sum_{j=1}^{p}\lambda_j'\tilde{\theta}_j \bpsi_j,\\
  \end{split}
\end{equation}
where $\lambda_j' = \sum_{k=t}^{T_n}w_k'$ for $ j_{t-1}+1 \leq j \leq j_t$. The risks of the MA estimators in (\ref{eq:est_all_nested})--(\ref{eq:est_group}) take the following forms:
\begin{equation*}\label{eq:risk_all_nested}
  R_n(\hat{\bmu}_{\bw|\mathcal{M}_{AN}}, \bmu) = \sum_{j=1}^{p}\left[ (1-\lambda_j)^2 \theta_j^2 + \lambda_j^2\sigma^2/n \right],
\end{equation*}
and
\begin{equation*}\label{eq:risk_all_nested}
  R_n(\hat{\bmu}_{\bw|\mathcal{M}_{G}}, \bmu) = \sum_{j=1}^{p}\left[ (1-\lambda_j')^2 \theta_j^2 + \lambda_j'^2\sigma^2/n \right].
\end{equation*}

\subsubsection{Proof of Corollary~\ref{tho:3}}

The proof of Corollary~\ref{tho:3} follows from Theorem~\ref{tho:1} and uses some proof techniques from Chapter 3.6 of \cite{Tsybakov2009}. Based on the oracle inequality in Theorem~\ref{tho:1}, there exists a constant $C>0$ such that for any $0< \delta <\infty$, we have
\begin{equation}\label{eq:jiaoleng1}
\begin{split}
   \mathbb{E}L_n(\hat{\bmu}_{\hat{\bw}_1|\mathcal{M}_G}, \bmu ) & \leq (1+C\delta) R_{n}(\hat{\bmu}_{\bw^*|\mathcal{M}_G}, \bmu) + \frac{C(1+\delta)^3T_n}{\delta n}  + \frac{C(1+\delta)T_n^2}{n}+  \frac{C(1+\delta)^2\mathbb{E}(\hat{\sigma}^2 - \sigma^2)^2}{\delta} \\
     & \leq [1+o(1)]R_{n}(\hat{\bmu}_{\bw^*|\mathcal{M}_G}, \bmu) + \frac{CT_n^2}{n} + CT_n \mathbb{E}(\hat{\sigma}^2 - \sigma^2)^2\\
     & \leq [1+o(1)]R_{n}(\hat{\bmu}_{\bw^*|\mathcal{M}_G}, \bmu) + \frac{C(\log p)^4}{n} + C(\log p)^2 \mathbb{E}(\hat{\sigma}^2 - \sigma^2)^2, 
\end{split}
\end{equation}
where the second inequality follows from the choice $\delta = \delta_n >0$, $\delta_n \to 0$, and $(1+\delta_n)^3/\delta_n = O(T_n)$, and the last inequality follows from the bound $T_n \leq C (\log p)^2$ given in Lemma 3.12 of \cite{Tsybakov2009}.

What remains is to establish a connection between $R_{n}(\hat{\bmu}_{\bw^*|\mathcal{M}_G}, \bmu)$ and the optimal MA risk over all nested models, $R_{n}(\hat{\bmu}_{\bw^*|\mathcal{M}_{AN}}, \bmu)$. This follows directly from Lemma 3.11 and Lemma 3.12 of \cite{Tsybakov2009}. Specifically, we have
\begin{equation}\label{eq:jiaoleng2}
  \begin{split}
       R_n(\hat{\bmu}_{\bw^*|\mathcal{M}_{G}}, \bmu)& \leq (1+3\rho_n)R_n(\hat{\bmu}_{\bw|\mathcal{M}_{AN}}, \bmu) + \frac{\sigma^2j_1}{n}\\
       & = (1+3\rho_n)R_n(\hat{\bmu}_{\bw|\mathcal{M}_{AN}}, \bmu) + \frac{C\log p}{n}.
  \end{split}
\end{equation}
By combining (\ref{eq:jiaoleng1}) with (\ref{eq:jiaoleng2}), we establish Corollary~\ref{tho:3}.

\section{Proof of the results in Section~\ref{sec:all-subset}}

\subsection{Preliminaries}

Let $\bw = (w_{{\mathcal{I}}})_{{\mathcal{I}} \subseteq \{1,\ldots,p \}}$ be a weight vector in  $\mathbb{R}^{2^{p}}$. The all-subset MA estimator based on $\bw$ is defined as
\begin{equation}\label{eq:MA_estimator}
  \hat{\bmu}_{\bw|\mathcal{M}_{AS}} = \sum_{{\mathcal{I}} \subseteq \{1,\ldots,p \}}w_{{\mathcal{I}}} \hat{\bmu}_{{\mathcal{I}}}.
\end{equation}
Using the spectral representation in Section~\ref{sec:pre_nested} again, we can write (\ref{eq:MA_estimator}) in an equivalent form
\begin{equation}\label{eq:MA_estimator_2}
  \hat{\bmu}_{\bw|\mathcal{M}_{AS}} = \sum_{{\mathcal{I}} \subseteq \{1,\ldots,p \}}w_{{\mathcal{I}}} \sum_{j \in {\mathcal{I}}} \tilde{\theta}_j \bpsi_j = \sum_{j=1}^{p}\left( \sum_{{\mathcal{I}}:j \in {\mathcal{I}}}w_{\mathcal{I}} \right) \tilde{\theta}_j \bpsi_j = \sum_{j=1}^{p} \gamma_j \tilde{\theta}_j \bpsi_j,
\end{equation}
where $\gamma_j \triangleq \sum_{{\mathcal{I}}:j \in {\mathcal{I}}}w_{\mathcal{I}}$. The performance of $\hat{\bmu}_{\bw|\mathcal{M}_{AS}}$ is measured by
\begin{equation}\label{eq:MA_risk}
  \begin{split}
     R_n\left(\hat{\bmu}_{\bw|\mathcal{M}_{AS}}, \bmu \right) & = n^{-1}\mathbb{E}\left\| \hat{\bmu}_{\bw|\mathcal{M}_{AS}}  - \bmu \right\|^2 = n^{-1}\mathbb{E}\left\|  \sum_{j=1}^{p} \gamma_j \tilde{\theta}_j \bpsi_j - \sum_{j=1}^{p} \theta_j \bpsi_j \right\|^2 \\
     &= \sum_{j=1}^{p} \mathbb{E}( \gamma_j \tilde{\theta}_j - \theta_j )^2 = \sum_{j=1}^{p}\left[\left(1-\gamma_j\right)^2\theta_j^2+\sigma^2\gamma_j^2/n\right].
  \end{split}
\end{equation}
The optimal all-subset MA risk is given by
\begin{equation}\label{eq:ideal_MA_risk}
\begin{split}
   R_n\left(\hat{\bmu}_{\bw^*|\mathcal{M}_{AS}}, \bmu \right) = \min_{\bw} R_n\left(\hat{\bmu}_{\bw^*|\mathcal{M}_{AS}}, \bmu \right) & = \sum_{j=1}^{p}\min_{\gamma_j}\left[\left(1-\gamma_j\right)^2\theta_j^2+\sigma^2\gamma_j^2/n\right] \\
     &     = \sum_{j=1}^{p}\frac{\theta_j^2 \sigma^2/n}{\theta_j^2+ \sigma^2/n}.
\end{split}
\end{equation}

\subsection{Proof of Theorem~\ref{theo:lower}}\label{sec:a:proof_lower}

The proof of the minimax lower bounds combines the Bayes risk analysis from \cite{Donoho1994Ideal, Averkamp2003} with the minimax problem reduction scheme in Chapter~3.3.2 of \cite{Tsybakov2009}. 

\subsubsection{Reduction to a minimax problem in a Gaussian sequence model}

For any measurable estimator $\hat{\bmu}$ based on $\by$, define its transformation coefficients as $\theta_j(\hat{\bmu}) \triangleq n^{-1} \hat{\bmu}^{\top} \bpsi_j , j=1,\ldots,p$. Note that $\theta_j(\hat{\bmu})$ is a statistic depending on $\by$. The risk of $\hat{\bmu}$ can be lower bounded as follows: 
  \begin{equation}\label{eq:tuijian5}
    \begin{split}
       R_n\left( \hat{\bmu},  \bmu \right) & = n^{-1}\mathbb{E}\left\| \hat{\bmu}  - \bmu \right\|^2 = n^{-1}\mathbb{E}\left\|  \sum_{j=1}^{p} \theta_j(\hat{\bmu}) \bpsi_j +\mathbf{b} - \sum_{j=1}^{p} \theta_j \bpsi_j \right\|^2 \geq \sum_{j=1}^{p} \mathbb{E}_{\by \sim P_{\btheta}} \left[ \theta_j(\hat{\bmu}) -  \theta_j \right]^2,
    \end{split}
  \end{equation}
  where $\mathbf{b}$ is the component of $\hat{\bmu}$ that is orthogonal to $\bpsi_j$, $j = 1, \ldots, p$. Here, $P_{\btheta}$ denotes the distribution of $\by$ under the model $\by = \sum_{j=1}^{p} \theta_j \bpsi_j + \beps$. The subscript $\by \sim P_{\btheta}$ in $\mathbb{E}_{\by \sim P_{\btheta}}$ is used to emphasize that the expectation is taken with respect to the random observation $\by$ drawn from $P_{\btheta}$. 
  
  Define the Gaussian sequence model as 
  \begin{equation}\label{eq:tanhua2}
    z_j = \theta_j + e_j, \quad j =1,\ldots,p,
  \end{equation}
  where $z_j \triangleq \by^{\top}\bpsi_j/n$, $\theta_j = \bmu^{\top}\bpsi_j/n$, and $e_j =  \beps^{\top} \bpsi_j /n$ are i.i.d. $N(0, \sigma^2/n)$. Let $\Bar{P}_{\btheta}$ denote the distribution of $\bz \triangleq (z_1, \ldots, z_n)^{\top}$ under the model (\ref{eq:tanhua2}), and let $\mathbb{E}_{\bz \sim \Bar{P}_{\btheta}}$ denote the corresponding expectation. The main idea of the following analysis is to reduce the expectation $\mathbb{E}_{\by \sim P_{\btheta}}$ in (\ref{eq:tuijian5}) to the expectation $\mathbb{E}_{\bz \sim \Bar{P}_{\btheta}}$. To this end, we follow the technique introduced in Chapter 3.3.2 of \cite{Tsybakov2009}. 
  
  When $\btheta = \boldsymbol{0}$, we have $\by = \beps$, where $\beps\sim N(\boldsymbol{0}, \sigma^2 \bI)$, and the density function of $\by$ is
  \begin{equation*}
    p_{\boldsymbol{0}}(\by) = \left( 2\pi \sigma^2 \right)^{-\frac{n}{2}} \exp\left(  - \frac{\sum_{i=1}^{n}y_i^2}{2\sigma^2} \right).
  \end{equation*}
  For general $\btheta$, we have
  \begin{equation*}
  \begin{split}
     p_{\btheta}(\by) & = \left( 2\pi \sigma^2 \right)^{-\frac{n}{2}} \exp\left(  - \frac{ \| \by - \sum_{j=1}^{p}\theta_j \bpsi_j \|^2}{2\sigma^2} \right) \\
       & = \left( 2\pi \sigma^2 \right)^{-\frac{n}{2}} \exp \left(  -\frac{\sum_{i=1}^{n}y_i^2 - 2 n \sum_{j=1}^{p}\theta_j z_j + n\sum_{j=1}^{p}\theta_j^2}{2\sigma^2} \right).
  \end{split}
  \end{equation*}
  Thus, the likelihood ratio between $p_{\btheta}$ and $p_{\boldsymbol{0}}$ is
  \begin{equation*}
    \frac{p_{\btheta}(\by)}{p_{\boldsymbol{0}}(\by)} = \exp \left(  \frac{  n \sum_{j=1}^{p}\theta_j z_j}{\sigma^2} - \frac{n\sum_{j=1}^{p}\theta_j^2}{2\sigma^2} \right) \triangleq S( \bz; \btheta ), 
  \end{equation*}
  where $\btheta \triangleq (\theta_1,\ldots,\theta_p)^{\top}$. Therefore, the last term in (\ref{eq:tuijian5}) can be written as
  \begin{equation}\label{eq:tuijian7}
    \begin{split}
       \mathbb{E}_{\by \sim P_{\btheta}} \left[ \theta_j(\hat{\bmu}) -  \theta_j \right]^2 & = \mathbb{E}_{\by \sim P_{\boldsymbol{0}}} \left[ \frac{p_{\btheta}(\by)}{p_{\boldsymbol{0}}(\by)} ( \theta_j(\hat{\bmu}) -  \theta_j )^2\right]\\
       & = \mathbb{E}_{\by \sim P_{\boldsymbol{0}}} \left[  ( \theta_j(\hat{\bmu}) -  \theta_j )^2S( \bz;\btheta )\right]\\
         & = \mathbb{E}_{\by \sim P_{\boldsymbol{0}}} \left\{\mathbb{E}_{\by \sim P_{\boldsymbol{0}}} \left[( \theta_j(\hat{\bmu}) -  \theta_j )^2 | \bz \right] S( \bz; \btheta ) \right\}, 
    \end{split}
  \end{equation}
  where the last equality follows from the formula of conditional expectation. By Jensen's inequality, we have
  \begin{equation}\label{eq:tuijian6}
  \begin{split}
       & \mathbb{E}_{\by \sim P_{\boldsymbol{0}}} \left[( \theta_j(\hat{\bmu}) -  \theta_j )^2 | \bz \right] \geq \left\{ \mathbb{E}_{\by \sim P_{\boldsymbol{0}}} \left[ \theta_j(\hat{\bmu}) | \bz \right] - \theta_j \right\}^2  = \left[\bar{\theta}_j(\bz) - \theta_j \right]^2,
  \end{split}
  \end{equation}
  where $\bar{\theta}_j(\bz) \triangleq \mathbb{E}_{\by \sim P_{\boldsymbol{0}}} [ \theta_j(\hat{\bmu}) | \bz ]$ depends on $\by$ only through $\bz$, which is the observation in the Gaussian sequence model (\ref{eq:tanhua2}). 
  
  Combining (\ref{eq:tuijian5}), (\ref{eq:tuijian7}), and (\ref{eq:tuijian6}), we obtain that for any estimator $\hat{\bmu}$ based on $\by$,
  \begin{equation}\label{eq:tanhua1}
    \begin{split}
         R_n\left( \hat{\bmu},  \bmu \right) &\geq \sum_{j=1}^{p} \mathbb{E}_{\by \sim P_{\btheta}} \left[ \theta_j(\hat{\bmu}) -  \theta_j \right]^2 \geq \sum_{j=1}^{p} \mathbb{E}_{\by \sim P_{\boldsymbol{0}}} \left\{ \left[\bar{\theta}_j(\bz) - \theta_j \right]^2  S( \bz; \btheta ) \right\}\\
          & = \sum_{j=1}^{p} \mathbb{E}_{\by \sim P_{\btheta}} \left[ \bar{\theta}_j(\bz) -  \theta_j \right]^2 = \sum_{j=1}^{p} \mathbb{E}_{\bz \sim \Bar{P}_{\btheta}} \left[ \bar{\theta}_j(\bz) -  \theta_j \right]^2.
    \end{split}
  \end{equation}
 Thus, the minimax risk ratio is lower bounded by 
  \begin{equation}\label{eq:tanhua3}
    \begin{split}
      & \min_{\hat{\bmu}}\max_{\bmu \in \mathcal{C}(\Theta)} \frac{R_n\left( \hat{\bmu},  \bmu \right)}{n^{-1} + R_n\left(\hat{\bmu}_{\bw^*|\mathcal{M}_{AS}}, \bmu \right)} \geq \min_{\hat{\bmu}}\max_{\bmu \in \mathcal{C}(\Theta)} \frac{\mathbb{E}_{\bz \sim \Bar{P}_{\btheta}}\sum_{j=1}^{p} \left[ \bar{\theta}_j(\bz) -  \theta_j \right]^2}{\frac{1}{n} + \sum_{j=1}^{p}\frac{\theta_j^2 \sigma^2/n}{\theta_j^2+ \sigma^2/n}}\\
        &\geq \min_{\hat{\boldsymbol{\vartheta}}}\max_{\btheta \in \Theta} \frac{\mathbb{E}_{\bz \sim \Bar{P}_{\btheta}}\sum_{j=1}^{p} ( \hat{\vartheta}_j -  \theta_j )^2}{\frac{1}{n}+\sum_{j=1}^{p}\frac{\theta_j^2 \sigma^2/n}{\theta_j^2+ \sigma^2/n}}\geq \min_{\hat{\boldsymbol{\vartheta}}}\max_{\btheta \in \Theta^*} \frac{\mathbb{E}_{\bz \sim \Bar{P}_{\btheta}}\sum_{j=1}^{p} ( \hat{\vartheta}_j -  \theta_j )^2}{\frac{1}{n} + \sum_{j=1}^{p}\frac{\theta_j^2 \sigma^2/n}{\theta_j^2+ \sigma^2/n}},
    \end{split}
  \end{equation}
  where the first inequality follows from (\ref{eq:tanhua1}) and (\ref{eq:ideal_MA_risk}). The second inequality follows from the fact that $\bar{\theta}_j(\bz)$ depends on $\by$ only through $\bz=(z_1,\ldots,z_p)^\top$, and hence the minimization can be extended to the class of all measurable estimators $\hat{\boldsymbol{\vartheta}}$ based solely on $\bz$. The last inequality follows from the inclusion $\Theta^*\subseteq\Theta$. Therefore, the last expression in (\ref{eq:tanhua3}) is exactly the minimax risk ratio problem in the Gaussian sequence model (\ref{eq:tanhua2}).

\subsubsection{A Bayes problem in one-dimensional case}

The main tool for lower bounding the last term in (\ref{eq:tanhua3}) is to evaluate the Bayes risk. We first consider a Bayesian problem in the one-dimensional case. Recall that the hardest cube is defined as $\Theta^* = \{ \btheta: 0 \leq |\theta_j| \leq \sqrt{2\sigma^2 \log p/n} \}$. For $0 < \kappa < 1$ and $0 < a \leq \sqrt{2\sigma^2 \log p/n}$, let
$$
F_{\kappa, a} \triangleq \kappa \delta_a + (1 - \kappa) \delta_0,
$$
where $\delta_c$ denotes the Dirac measure with unit mass at $c$, and the precise values of $\kappa$ and $a$ will be specified later. We are interested in the Bayes risk for estimating $\theta_1 \in \mathbb{R}$ based on the observation $z_1 = \theta_1 + e_1$, where the prior distribution of $\theta_1$ is $F_{\kappa, a}$, and $e_1$ is distributed as $N(0, \sigma^2/n)$. For an estimator $\hat{\vartheta}_1$, the Bayes risk is defined as $\mathbb{E}_{\theta_1 \sim F_{\kappa, a}}\mathbb{E}_{z_1 \sim \bar{P}_{\theta_1}}(\hat{\vartheta}_1 - \theta_1)^2$.

Let $f$ denote the density function of $e_1$, which is given by $f(x) = \frac{1}{\sigma / \sqrt{n}} \phi(\frac{x}{\sigma / \sqrt{n}})$, where $\phi$ denotes the density function of the standard normal distribution. In the univariate setting above, the Bayes estimator of $\theta_1$ is the posterior mean. Given $z_1 = x$, it is given by
\begin{equation}\label{eq:tuijian1}
  \begin{split}
     \vartheta_{\kappa, a}(x)& \triangleq \mathbb{E}( \theta_1 \mid z_1= x )  = 0 \times \mathbb{P}( \theta_1 = 0 \mid z_1= x ) + a \times \mathbb{P}( \theta_1 = a \mid z_1= x )  =\frac{\kappa f(x-a)}{\kappa f(x-a)+(1-\kappa) f(x)} a.
  \end{split}
\end{equation}
Thus, the Bayes risk of $\vartheta_{\kappa, a}$ is lower bounded by
\begin{equation}\label{eq:xiluo1}
  \begin{split}
      \mathbb{E}_{\theta_1 \sim F_{\kappa, a}} \mathbb{E}_{z_1 \sim \bar{P}_{\theta_1}}\left(\vartheta_{\kappa, a}-\theta_1\right)^2 &=\kappa \int_{-\infty}^{+\infty}\left[\vartheta_{\kappa, a}(x)-a\right]^2 f(x-a) d x+(1-\kappa) \int_{-\infty}^{+\infty} \vartheta^2_{\kappa, a}(x) f(x) d x\\
       & \geq \kappa \int_{-\infty}^{+\infty}\left[\vartheta_{\kappa, a}(x)-a\right]^2 f(x-a) d x\\
       & =\kappa a^2 \int_{-\infty}^{+\infty}\left[\frac{(1-\kappa) f(x)}{\kappa f(x-a)+(1-\kappa) f(x)}\right]^2 f(x-a) d x \\
       & =(1-\kappa)^2 \kappa a^2 \int_{-\infty}^{+\infty} \frac{f^2(x)}{[\kappa f(x-a)+(1-\kappa) f(x)]^2} f(x-a) d x,
  \end{split}
\end{equation}
where the second equality follows from (\ref{eq:tuijian1}).

Let us now lower bound the integrand in the last term of (\ref{eq:xiluo1}). Recall that $f(x)$ is the density function of $N(0,\sigma^2/n)$. For any $\alpha \in (0,1)$, there exists a positive constant $c = -\sigma \Phi^{-1}(\frac{1-\alpha}{2})$ such that
\begin{equation}\label{eq:tuijian2}
  \int_{-c/\sqrt{n}}^{c/\sqrt{n}} f(x) d x = \alpha,
\end{equation}
where $\Phi$ denotes the cumulative distribution function of the standard normal distribution. Furthermore, for any $\beta>0$, suppose that $\kappa$ and $a$ are chosen such that 
\begin{equation}\label{eq:xiluo2}
  \beta f\left(a+ \frac{c}{\sqrt{n}} \right) \geq \frac{\kappa}{1-\kappa} f(0). 
\end{equation}
Then, for any $x\in[a-c/\sqrt n,a+c/\sqrt n]$, the unimodality and symmetry of the Gaussian density imply that
\begin{equation*}
  \beta f(x) \geq \beta f\left( a+\frac{c}{\sqrt{n}} \right) \geq \frac{\kappa}{1-\kappa} f(0) \geq \frac{\kappa}{1-\kappa} f\left(x - a\right).
\end{equation*}
Consequently, for $a-c/\sqrt n\leq x\leq a+c/\sqrt n$, the integrand in the last term of (\ref{eq:xiluo1}) satisfies 
\begin{equation}\label{eq:tuijian3}
\begin{split}
     \frac{f^2(x)}{[\kappa f(x-a)+(1-\kappa) f(x)]^2} f(x-a) &\geq \frac{f^2(x)}{[(1-\kappa)\beta f(x)+(1-\kappa) f(x)]^2} f(x-a)\\
     & = \frac{f(x-a)}{(1-\kappa)^2(1+\beta)^2}. 
\end{split}
\end{equation}

Thus, for $0<\kappa<1$ and $0<a\leq \sqrt{2\sigma^2\log p/n}$, if $\kappa$ and $a$ are chosen such that (\ref{eq:xiluo2}) holds, then 
\begin{equation}\label{eq:ssss1}
  \begin{split}
        \mathbb{E}_{\theta_1 \sim F_{\kappa, a}} \mathbb{E}_{z_1 \sim \bar{P}_{\theta_1}}\left(\vartheta_{\kappa, a}-\theta_1\right)^2 &\geq (1-\kappa)^2 \kappa a^2 \int_{a-\frac{c}{\sqrt{n}}}^{a+\frac{c}{\sqrt{n}}} \frac{f^2(x)}{[\kappa f(x-a)+(1-\kappa) f(x)]^2} f(x-a) d x \\
       & \geq (1-\kappa)^2 \kappa a^2 \int_{a-\frac{c}{\sqrt{n}}}^{a+\frac{c}{\sqrt{n}}} \frac{f(x-a)}{(1-\kappa)^2(1+\beta)^2} d x\\
       & = \frac{\kappa a^2}{(1+\beta)^2} \int_{a-\frac{c}{\sqrt{n}}}^{a+\frac{c}{\sqrt{n}}} f(x-a) dx = \frac{\kappa a^2}{(1+\beta)^2} \int_{-\frac{c}{\sqrt{n}}}^{\frac{c}{\sqrt{n}}} f(x) dx\\
        &= \frac{\alpha}{(1+\beta)^2} \kappa a^2,
  \end{split}
\end{equation}
where the first inequality follows from (\ref{eq:xiluo1}), the second inequality follows from (\ref{eq:tuijian3}), and the last equality follows from (\ref{eq:tuijian2}).

\subsubsection{From one-dimensional case to multivariate case}

We now evaluate the Bayes risk for estimating the multivariate vector $\btheta$ in the Gaussian sequence model (\ref{eq:tanhua2}). Suppose that $\btheta$ follows the product prior distribution $Q_p \triangleq \otimes_{j=1}^pF_{\kappa, a}$, where the parameters $\kappa$ and $a$ satisfy $0<\kappa<1$, $0<a\leq \sqrt{2\sigma^2\log p/n}$, and the condition in (\ref{eq:xiluo2}). The Bayes risk is defined as $\mathbb{E}_{\btheta \sim Q_p} \mathbb{E}_{\bz \sim \bar{P}_{\btheta}} \| \hat{\btheta} - \btheta \|^2 = \mathbb{E}_{\btheta \sim Q_p} \mathbb{E}_{\bz \sim \bar{P}_{\btheta}} \sum_{j=1}^{p}(\hat{\theta}_j - \theta_j)^2$. Under this product prior, the components of $\btheta$ are independent and identically distributed according to $F_{\kappa,a}$. Therefore, the Bayes estimator of $\btheta$ given the observation $\bz=(z_1,\ldots,z_p)^\top$ in (\ref{eq:tanhua2}) is obtained by applying the univariate Bayes rule coordinate-wise, yielding
\begin{equation*}
    \boldsymbol{\vartheta}_{\kappa, a}(\bz) \triangleq \left[\vartheta_{\kappa, a}(z_1),\ldots,\vartheta_{\kappa, a}(z_p) \right]^{\top},
\end{equation*}
where $\vartheta_{\kappa, a}(\cdot)$ is the univariate Bayes rule defined in (\ref{eq:tuijian1}). Thus, for any estimator $\hat{\btheta}$, by the optimality of the Bayes estimator under squared error loss and (\ref{eq:ssss1}), we have 
\begin{equation}\label{eq:bayes_riskss}
  \mathbb{E}_{\btheta \sim Q_p} \mathbb{E}_{\bz \sim \bar{P}_{\btheta}} \sum_{j=1}^{p}(\hat{\theta}_j - \theta_j)^2 \geq \mathbb{E}_{\btheta \sim Q_p} \mathbb{E}_{\bz \sim \bar{P}_{\btheta}} \sum_{j=1}^{p}(\vartheta_{\kappa, a}(z_j) - \theta_j)^2 \geq p\kappa a^2\frac{\alpha}{(1+\beta)^2}. 
\end{equation}

The main task now is to determine $\kappa$ and $a$ such that the required conditions are satisfied. Recall that $0<\alpha<1$ and $c=-\sigma\Phi^{-1}((1-\alpha)/2)$ are chosen such that the equality in (\ref{eq:tuijian2}) holds. We first fix $\alpha$, and consequently $c$, as well as the positive constant $\beta$ in (\ref{eq:xiluo2}). We set $\kappa \triangleq (\log p)^3/p$, which satisfies $0<\kappa<1$ whenever $p\geq \underline{p}$, where $\underline{p}$ is a sufficiently large positive integer. Given the fixed parameters $\alpha$, $c$, and $\beta$, the condition in (\ref{eq:xiluo2}) requires 
\begin{equation*}
  \frac{1}{\sqrt{2\pi\sigma^2/n}}\exp\left[ - \frac{\left( a+ c/\sqrt{n} \right)^2}{2\sigma^2/n} \right] \geq \frac{\kappa}{\beta(1-\kappa)} \frac{1}{\sqrt{2\pi\sigma^2/n}}.
\end{equation*}
Equivalently, after taking logarithms, this condition becomes
\begin{equation*}
\begin{split}
     \frac{\left( a+ c/\sqrt{n} \right)^2}{2\sigma^2/n}& \leq -\log \kappa + \log\beta + \log (1-\kappa)   = \log p -3 \log\log p + \log\beta+ \log \left[ 1 - (\log p)^3/p \right].
\end{split}
\end{equation*}
To maximize the resulting Bayes risk lower bound, we choose $a$ as large as possible while satisfying the above constraint. Thus, we take $a$ to attain equality, which gives
\begin{equation}\label{eq:a_define}
  a\triangleq \sqrt{\frac{2\sigma^2}{n}\big\{ \log p -3 \log\log p + \log\beta+ \log \left[ 1 - (\log p)^3/p \right] \big\}} + \frac{\sigma \Phi^{-1}(\frac{1-\alpha}{2})}{\sqrt{n}}.
\end{equation}
We note that the choice of $a$ in (\ref{eq:a_define}) has the following two important properties, provided that $\alpha \in (0, 1)$ and $\beta>0$ are selected such that $-\sqrt{\log\log p} \lesssim \Phi^{-1}(\frac{1-\alpha}{2}) \lesssim 1$ and $-\log\log p\lesssim \log \beta \leq 3\log\log p$. First, there exists a sufficiently large integer $\underline{p}$ such that, for all $p\geq \underline{p}$, $\sqrt{\sigma^2/n} \leq a \leq \sqrt{2\sigma^2 \log p/n}$. Second, if $p\to\infty$ as $n\to\infty$, then $a = [1+o(1)]\sqrt{2\sigma^2 \log p/n}$.

Before deriving the lower bound for the minimax risk ratio, we analyze the following event under the prior distribution $Q_p$. Define $N(\btheta) \triangleq |\{\theta_j \neq 0, j=1, \ldots, p\}|$ and $\mathcal{A} \triangleq \{N(\btheta) \leq p \kappa+3(p \kappa)^{2 / 3}\}$. Our goal is to upper bound the probability of the complement event $\mathcal{A}^c$. Define the binary random variable $X_j = 1_{\{ \theta_j= {a} \}}$, $j = 1,\ldots,p$. Under the prior distribution $Q_p$, the variables $X_j$ are independent Bernoulli random variables with $\mathbb{E}(X_j)=\kappa$. Therefore, $X_j-\kappa$ are independent, mean-zero random variables satisfying $|X_j - \kappa| \leq 1$. By Bernstein's inequality, we have
  \begin{equation}\label{eq:proaba_A}
    \begin{split}
      \mathbb{P}(\mathcal{A}^c)  & =  \mathbb{P}\left( \sum_{j=1}^{p}X_j - p\kappa > 3(p \kappa)^{2 / 3}\right) \leq \exp\left( - \frac{\frac{9}{2}(p \kappa)^{4 / 3}}{p \kappa(1 - \kappa ) + (p \kappa)^{2 / 3}}  \right)\\
         & \leq \exp\left( - \frac{\frac{9}{2}(p \kappa)^{4 / 3}}{p \kappa + (p \kappa)^{2 / 3}}   \right)  = \exp\left( - \frac{9 (p \kappa)^{4/3}}{2p \kappa + 2(p \kappa)^{2/3}} \right) \leq  \exp\left( - 2.7 (p \kappa)^{1/3} \right)\\
         & = \exp\left( - 2.7 \log p \right) \leq \frac{1}{p^2} = \frac{\kappa}{p(\log p)^3}, 
    \end{split}
  \end{equation}
  where the third inequality follows from $2(p \kappa)^{2/3} \leq 2 \times \frac{2}{3}\times(p \kappa -1) \leq \frac{4p \kappa}{3}$ for $p \kappa \geq 1$, and the last equality uses $\kappa = (\log p)^3/p$.
  
  \subsubsection{Lower bounding the minimax risk by Bayes risk}

We are now in a position to derive a lower bound for the minimax risk ratio. Our approach primarily follows the methods developed in \cite{Donoho1994Ideal} and \cite{Averkamp2003}, while retaining all terms necessary to obtain the lower bound for the finite $p$. 
  
Since for $p\geq\underline{p}$ we have $a\leq \sqrt{2\sigma^2\log p/n}$, the support of the prior distribution $Q_p = \otimes_{j=1}^pF_{\kappa, a}$ is contained in $\Theta^*$. Therefore, the minimax risk ratio in the last term of (\ref{eq:tanhua3}) is lower bounded by 
  \begin{equation}\label{eq:lower_bound_first}
    \begin{split}
      & \min_{\hat{\boldsymbol{\vartheta}}}\max_{\btheta \in \Theta^*} \frac{\mathbb{E}_{\bz \sim \Bar{P}_{\btheta}}\sum_{j=1}^{p} ( \hat{\vartheta}_j -  \theta_j )^2}{\frac{1}{n} + \sum_{j=1}^{p}\frac{\theta_j^2 \sigma^2/n}{\theta_j^2+ \sigma^2/n}} \geq \min_{\hat{\boldsymbol{\vartheta}}} \mathbb{E}_{\btheta \sim Q_p} \left[ \frac{\mathbb{E}_{\bz \sim \Bar{P}_{\btheta}}\sum_{j=1}^{p} ( \hat{\vartheta}_j -  \theta_j )^2}{\frac{1}{n} + \sum_{j=1}^{p}\frac{\theta_j^2 \sigma^2/n}{\theta_j^2+ \sigma^2/n}} \right]\\
       &\geq \min_{\hat{\boldsymbol{\vartheta}}} \mathbb{E}_{\btheta \sim Q_p} \left[ \frac{\mathbb{E}_{\bz \sim \Bar{P}_{\btheta}}\sum_{j=1}^{p} ( \hat{\vartheta}_j -  \theta_j )^2}{\frac{1}{n} + \frac{\sigma^2}{n}N(\btheta)} \right] = \mathbb{E}_{\btheta \sim Q_p} \left[ \frac{\mathbb{E}_{\bz \sim \Bar{P}_{\btheta}}\sum_{j=1}^{p} ( \hat{\vartheta}_j^{\mathrm{r}} -  \theta_j )^2}{\frac{1}{n} + \frac{\sigma^2}{n}N(\btheta)} \right]\\
        & \geq \mathbb{E}_{\btheta \sim Q_p} \left[ \frac{\mathbb{E}_{\bz \sim \Bar{P}_{\btheta}}\sum_{j=1}^{p} ( \hat{\vartheta}_j^{\mathrm{r}} -  \theta_j )^2}{\frac{1}{n} + \frac{\sigma^2}{n}N(\btheta)}1_{\mathcal{A}} \right]  \geq \frac{1}{\frac{1}{n}+\frac{\sigma^2}{n}\left( p \kappa+{3(p \kappa)^{2 / 3}} \right)} \mathbb{E}_{\btheta \sim Q_p} \mathbb{E}_{\bz \sim \bar{P}_{\btheta}}\left[\sum_{j=1}^{p} ( \hat{\vartheta}_j^{\mathrm{r}} -  \theta_j )^2 1_{\mathcal{A}}\right]\\
       & = \frac{1}{\frac{1}{n}+\frac{\sigma^2}{n}\left( p \kappa+{3(p \kappa)^{2 / 3}} \right)}  \left\{\mathbb{E}_{\btheta \sim Q_p} \mathbb{E}_{\bz \sim \Bar{P}_{\btheta}}\sum_{j=1}^{p} ( \hat{\vartheta}_j^{\mathrm{r}} -  \theta_j )^2  - \mathbb{E}_{\btheta \sim Q_p} \mathbb{E}_{\bz \sim \Bar{P}_{\btheta}}\left[\sum_{j=1}^{p} ( \hat{\vartheta}_j^{\mathrm{r}} -  \theta_j )^2 1_{\mathcal{A}^c} \right] \right\}\\
       & \geq \frac{1}{\frac{1}{n}+\frac{\sigma^2}{n}\left( p \kappa+{3(p \kappa)^{2 / 3}} \right)}  \left\{\mathbb{E}_{\btheta \sim Q_p} \mathbb{E}_{\bz \sim \Bar{P}_{\btheta}} \sum_{j=1}^{p}(\vartheta_{\kappa, a}(z_j) - \theta_j)^2  - \mathbb{E}_{\btheta \sim Q_p} \mathbb{E}_{\bz \sim \Bar{P}_{\btheta}}\left[\sum_{j=1}^{p} ( \hat{\vartheta}_j^{\mathrm{r}} -  \theta_j )^2 1_{\mathcal{A}^c} \right]  \right\}, 
    \end{split}
  \end{equation}
  where the second inequality follows from the fact that when $p\geq\underline{p}$, we have $ a \geq \sqrt{\sigma^2/n}$, and hence 
  \begin{equation*}
    \sum_{j=1}^{p}\frac{\theta_j^2 \sigma^2/n}{\theta_j^2+ \sigma^2/n} \leq \sum_{j=1}^{p}\min\left( \theta_j^2, \sigma^2/n \right) \leq \frac{\sigma^2}{n}\sum_{j=1}^{p}1_{\{ \theta_j \neq 0 \}} = \frac{\sigma^2}{n}N(\btheta)
  \end{equation*}
  for any $\btheta \sim Q_p$. The first equality follows from by defining $\hat{\boldsymbol{\vartheta}}^{\mathrm{r}} = (\hat{\vartheta}_1^{\mathrm{r}}, \ldots, \hat{\vartheta}_p^{\mathrm{r}})^{\top}$ as the minimizer of the ratio-type Bayes risk: 
  \begin{equation}\label{eq:risk_ratio_problem}
    \hat{\boldsymbol{\vartheta}}^{\mathrm{r}} \triangleq \argmin_{\hat{\boldsymbol{\vartheta}}} \mathbb{E}_{\btheta \sim Q_p} \mathbb{E}_{\bz \sim \Bar{P}_{\btheta}}  \frac{\sum_{j=1}^{p} ( \hat{\vartheta}_j -  \theta_j )^2}{\frac{1}{n} + \frac{\sigma^2}{n}N(\btheta)} = \argmin_{\hat{\boldsymbol{\vartheta}}} \mathbb{E}_{\btheta \sim Q_p} \mathbb{E}_{\bz \sim \Bar{P}_{\btheta}} \| \hat{\boldsymbol{\vartheta}} - \btheta \|^2 \omega(\btheta),  
  \end{equation}
  where $\omega(\btheta) \triangleq n/(1+\sigma^2N(\btheta))$. The forth inequality in (\ref{eq:lower_bound_first}) follows from the definition of the event $\mathcal{A}$, and the fifth inequality follows from (\ref{eq:bayes_riskss}) that $(\vartheta_{\kappa,a}(z_1),\ldots,\vartheta_{\kappa,a}(z_p))^\top$ is the Bayes estimator of $\btheta$ under the usual squared error loss. 
 
 To proceed with the proof, we need to construct an upper bound for $\mathbb{E}_{\btheta \sim Q_p} \mathbb{E}_{\bz \sim \Bar{P}_{\btheta}}\sum_{j=1}^{p} ( \hat{\vartheta}_j^{\mathrm{r}} -  \theta_j )^2 1_{\mathcal{A}^c}$. Let $\mathbb{E}$ denote the joint distribution of $\btheta$ and $\bz$. Then, we have  
 \begin{equation*}
 \begin{split}
    \mathbb{E}_{\btheta \sim Q_p} \mathbb{E}_{\bz \sim \Bar{P}_{\btheta}} \| \hat{\boldsymbol{\vartheta}} - \btheta \|^2 \omega(\btheta) & = \mathbb{E}_{\bz } \mathbb{E} \left[\| \hat{\boldsymbol{\vartheta}} - \btheta \|^2 \omega(\btheta) | \bz \right]  = \mathbb{E}_{\bz } \sum_{j=1}^{p} \mathbb{E} \left[ \hat{\vartheta}_j^2 \omega(\btheta) -  2\theta_j\hat{\vartheta}_j\omega(\btheta) + \theta_j^2\omega(\btheta) | \bz \right].
 \end{split}
 \end{equation*}
 Minimizing the conditional expectation with respect to $\hat{\vartheta}_j$ yields
$$
0 \leq \hat{\vartheta}_j^{\mathrm{r}}(\bz)=\frac{\mathbb{E}[\theta_j \omega(\btheta ) | \bz ]}{\mathbb{E}[\omega(\btheta ) | \bz ]} \leq \frac{\mathbb{E}[a \omega(\btheta ) | \bz ]}{\mathbb{E}[\omega(\btheta ) | \bz ]} = a.
$$
Therefore, we have the upper bound 
\begin{equation}\label{eq:prob_ine_2}
\begin{split}
  \mathbb{E}_{\btheta \sim Q_p} \mathbb{E}_{\bz \sim \Bar{P}_{\btheta}}\left[\sum_{j=1}^{p} ( \hat{\vartheta}_j^{\mathrm{r}} -  \theta_j )^2 1_{\mathcal{A}^c} \right] \leq p{a}^2\mathbb{P}(\mathcal{A}^c) \leq  \frac{\kappa a^2}{(\log p)^3},
\end{split}
\end{equation}
where the first inequality follows from both $\hat{\vartheta}_j^{\mathrm{r}}$ and $\theta_j$ are between $0$ and $a$, and the second inequality follows from (\ref{eq:proaba_A}). Inserting the inequalities (\ref{eq:bayes_riskss}) and (\ref{eq:prob_ine_2}) into (\ref{eq:lower_bound_first}) gives  
  \begin{equation}\label{eq:aaaasha}
\begin{aligned}
\min_{\hat{\boldsymbol{\vartheta}}}\max_{\btheta \in \Theta^*} \frac{\mathbb{E}_{\bz \sim \Bar{P}_{\btheta}}\sum_{j=1}^{p} ( \hat{\vartheta}_j -  \theta_j )^2}{\frac{1}{n} + \sum_{j=1}^{p}\frac{\theta_j^2 \sigma^2/n}{\theta_j^2+ \sigma^2/n}} & \geq \frac{1}{\frac{1}{n}+\frac{\sigma^2}{n}\left( p \kappa+{3(p \kappa)^{2 / 3}} \right)}  \left[p \kappa {a}^2 \frac{\alpha}{(1+\beta)^2}-\frac{\kappa a^2}{(\log p)^3}\right] \\
& = \frac{p\kappa a^2}{\frac{1}{n}+\frac{\sigma^2}{n}\left( p \kappa+{3(p \kappa)^{2 / 3}} \right)}  \left[ \frac{\alpha}{(1+\beta)^2} - \frac{1}{p(\log p)^3} \right]\\
& = \frac{n}{ \sigma^2}\left[ \frac{\alpha}{(1+\beta)^2} - \frac{1}{p(\log p)^3} \right] \left[ \frac{{a}^2}{1+1/[\sigma^2(\log p)^3]+ 3/\log p} \right],
\end{aligned}
\end{equation}
where the last step in (\ref{eq:aaaasha}) follows from the definition of $\kappa = (\log p)^3/p$. Combining (\ref{eq:aaaasha}) with definition of $a$ in (\ref{eq:a_define}) and the relation (\ref{eq:tanhua3}), we have proved that
\begin{equation}\label{eq:shabini1}
\begin{split}
     &\min_{\hat{\bmu}}\max_{\bmu \in \mathcal{C}(\Theta)} \frac{R_n\left( \hat{\bmu},  \bmu \right)}{n^{-1}+R_n\left(\hat{\bmu}_{\bw^*|\mathcal{M}_{AS}}, \bmu \right)} \geq \left[ \frac{\alpha}{(1+\beta)^2} - \frac{1}{p(\log p)^3} \right] \left[\frac{1}{1+1/[\sigma^2(\log p)^3]+ 3/\log p}  \right]\\
     & \quad\quad\quad\quad\times\left\{\sqrt{2\big\{ \log p -3 \log\log p + \log\beta+ \log \left[ 1 - (\log p)^3/p \right] \big\}} +  \Phi^{-1}\left(\frac{1-\alpha}{2}\right)\right\}^2.
\end{split}
\end{equation}

\subsubsection{Finalizing the proof}

In the asymptotic regime where $p\to\infty$ as $n\to\infty$, we choose $\alpha = 1 - 2\Phi( -\sqrt{2\log\log p} )$ and $\beta = 1/\log p$. When $n \to \infty$, we have $\alpha \to 1$ and $\beta \to 0$. Consequently, the first factor in (\ref{eq:shabini1}) satisfies
\begin{equation*}
  \left[ \frac{\alpha}{(1+\beta)^2} - \frac{1}{p(\log p)^3} \right] \left[\frac{1}{1+1/[\sigma^2(\log p)^3]+ 3/\log p}  \right] \to 1.
\end{equation*}
Moreover, 
\begin{equation*}
  \sqrt{2\big\{ \log p -3 \log\log p + \log\beta+ \log \left[ 1 - (\log p)^3/p \right] \big\}} =[1+o(1)] \sqrt{2\log p}
\end{equation*}
and
\begin{equation*}
  \Phi^{-1}\left(\frac{1-\alpha}{2}\right) = -\sqrt{2\log\log p} = o\left( \sqrt{2\log p} \right).
\end{equation*}
Therefore, in the case $p \to \infty$, we obtain the minimax lower bound
\begin{equation*}
  \min_{\hat{\bmu}}\max_{\bmu \in \mathcal{C}(\Theta)} \frac{R_n\left( \hat{\bmu},  \bmu \right)}{n^{-1}+R_n\left(\hat{\bmu}_{\bw^*|\mathcal{M}_{AS}}, \bmu \right)}\geq[1+o(1)]2\log p.
\end{equation*}

For finite $p$, we can choose $\alpha\in(0,1)$ and $\beta\in(0,1)$ as some fixed constants. By the monotonicity of the lower bound in (\ref{eq:shabini1}) with respect to $p$ when $p$ is large enough, it is straightforward to verify that there exists a sufficiently large positive integer $\underline{p}$ such that, for all $p\geq \underline{p}$, the lower bound is greater than a constant $\underline{c}>1$. 

\subsection{Proof of Theorem~\ref{theo:upper}}\label{sec:a:proof_upper}

\subsubsection{An equivalent expression for the Mallows-type criterion (\ref{eq:MMA_cor})}\label{subsubsec:equivalent}

Recalling that $\hat{\bmu}_{j} = \tilde{\theta}_j \bpsi_j$ and $\tilde{\theta}_j = n^{-1} \by^{\top} \bpsi_j $, we rewrite the criterion in (\ref{eq:MMA_cor}) as follows:
\begin{equation}\label{eq:MMA_cor1}
\begin{split}
     n^{-1}\left\| \by - \sum_{j=1}^{p}w_j \hat{\bmu}_{j} \right\|^2 + 2\lambda^2_2 \sigma^2 \bw^{\top} \boldsymbol{1} &=n^{-1} \left\| \sum_{j=1}^{p}\tilde{\theta}_j \bpsi_j - \sum_{j=1}^{p}w_j \tilde{\theta}_j \bpsi_j + \ba \right\|^2 + 2\lambda^2_2 \sigma^2 \sum_{j=1}^{p}w_j\\
     & = \sum_{j=1}^{p}\left[(1-w_j)^2\tilde{\theta}_j^2 + 2\lambda_2^2 \sigma^2 w_j \right] +n^{-1}\|\ba\|^2 \\
     & = \sum_{j=1}^{p}\left[  \tilde{\theta}_j^2 w_j^2 - (  2 \tilde{\theta}_j^2 - 2\lambda_2^2 \sigma^2  ) w_j \right] + \sum_{j=1}^{p}\tilde{\theta}_j^2 +n^{-1}\|\ba\|^2,
\end{split}
\end{equation}
where $\ba$ is the component of $\by$ that is orthogonal to $\bpsi_1, \ldots, \bpsi_p$ under the inner product $\langle \cdot, \cdot \rangle$. Since the last two terms in (\ref{eq:MMA_cor1}) are independent of $\bw$, the minimizer of the criterion over $[0,1]^p$ is given by
\begin{equation}\label{eq:siwa5}
  \hat{w}_{2j} = \left( 1 - \frac{\lambda_2^2 \sigma^2}{\tilde{\theta}_j^2} \right)_+, \quad j = 1, \ldots, p, 
\end{equation}
where $\hat{w}_{2j}$ is the $j$-th element of $\hat{\bw}_{2}$ defined in (\ref{eq:MMA_cor}). Note that $\hat{w}_{2j}$ depends only on $\tilde{\theta}_j$, where $\tilde{\theta}_j \sim N(\theta_j, \sigma^2/n)$. The risk of the resulting MA estimator (\ref{eq:MMA_est}) is given by
\begin{equation}\label{eq:siwa2}
  \begin{split}
      \mathbb{E} L_n\left(\hat{\bmu}_{\hat{\bw}_2|\mathcal{M}_{U}}, \bmu \right) & = n^{-1}\mathbb{E}\left\| \hat{\bmu}_{\hat{\bw}_2|\mathcal{M}_{U}}  - \bmu \right\|^2 = n^{-1}\mathbb{E}\left\|  \sum_{j=1}^{p} \hat{w}_{2j} \tilde{\theta}_j \bpsi_j - \sum_{j=1}^{p} \theta_j \bpsi_j \right\|^2\\
        &= \sum_{j=1}^{p} \mathbb{E}\left( \hat{w}_{2j} \tilde{\theta}_j - \theta_j \right)^2, \\
  \end{split}
\end{equation}
where the last equality follows from the orthogonality of $\bpsi_1, \ldots, \bpsi_p$ and $n^{-1}\|\bpsi_j\|^2 = 1$ for $j=1,\ldots,p$.

\subsubsection{Univariate risk bound}

To upper bound (\ref{eq:siwa2}), the key step is to upper bound the univariate risk $\mathbb{E}( \hat{w}_{2j} \tilde{\theta}_j - \theta_j )^2$ for $j = 1, \ldots, p$. By (\ref{eq:siwa5}), $\hat{w}_{2j} \tilde{\theta}_j$ can be expressed as
\begin{equation}\label{eq:siwa4}
  \begin{split}
    \hat{w}_{2j} \tilde{\theta}_j&= \left( 1 - \frac{\lambda_2^2 \sigma^2}{\tilde{\theta}_j^2} \right)_+ \tilde{\theta}_j =\left\{\begin{array}{ll}
\tilde{\theta}_j - \frac{\lambda_2^2\sigma^2}{\tilde{\theta}_j} &\quad \tilde{\theta}_j > \lambda_2 \sigma \\
0 &\quad - \lambda_2 \sigma \leq \tilde{\theta}_j \leq \lambda_2 \sigma\\
\tilde{\theta}_j - \frac{\lambda_2^2\sigma^2}{\tilde{\theta}_j} & \quad \tilde{\theta}_j < - \lambda_2 \sigma
\end{array}\right.\\
& = \tilde{\theta}_j - \left\{\begin{array}{ll}
 \frac{\lambda_2^2\sigma^2}{\tilde{\theta}_j} &\quad \tilde{\theta}_j > \lambda_2 \sigma \\
 \tilde{\theta}_j &\quad - \lambda_2 \sigma \leq \tilde{\theta}_j \leq \lambda_2 \sigma\\
 \frac{\lambda_2^2\sigma^2}{\tilde{\theta}_j} & \quad \tilde{\theta}_j < - \lambda_2 \sigma.
\end{array}\right. = \tilde{\theta}_j - S(\tilde{\theta}_j),\\
  \end{split}
\end{equation}
where 
\begin{equation*}
  S(\theta)= \begin{cases} \frac{\lambda_2^2 \sigma^2 }{ \theta}, & |\theta|>\lambda_2 \sigma \\ \theta, & |\theta| \leq \lambda_2 \sigma .\end{cases}
\end{equation*}

In the following, we derive three upper bounds for $\mathbb{E}( \hat{w}_{2j} \tilde{\theta}_j - \theta_j )^2$, each of which is suited to a different regime of the magnitude of $|\theta_j|$. The first bound is given by
\begin{equation}\label{eq:bound1}
\begin{split}
     \mathbb{E}( \hat{w}_{2j} \tilde{\theta}_j - \theta_j )^2& = \mathbb{E} \left[  \tilde{\theta}_j - \theta_j - S(\tilde{\theta}_j) \right]^2 \leq \mathbb{E}(\tilde{\theta}_j - \theta_j)^2 + \mathbb{E}\big[S^2(\tilde{\theta}_j)\big] + 2 \sqrt{\mathbb{E}(\tilde{\theta}_j - \theta_j)^2} \sqrt{\mathbb{E}\big[S^2(\tilde{\theta}_j)\big]} \\
     & = \left[\sqrt{\mathbb{E}(\tilde{\theta}_j - \theta_j)^2} + \sqrt{\mathbb{E}\big[S^2(\tilde{\theta}_j)\big]} \right]^2 \leq \left(\frac{\sigma}{\sqrt{n}} + \lambda_2 \sigma \right)^2 \\
     & = \left( \lambda_2 + \frac{1}{\sqrt{n}} \right)^2 \sigma^2, 
\end{split}
\end{equation}
where the first inequality follows from the Cauchy–Schwarz inequality, and the second inequality follows from $\Var(\tilde{\theta}_j) = \sigma^2/n$ and $|S(\theta)| \leq \lambda_2 \sigma$ for all $\theta$. Note that the upper bound (\ref{eq:bound1}) is valid without the Gaussian assumption on $\tilde{\theta}_j$.

The second upper bound is
\begin{equation}\label{eq:bound2}
  \begin{split}
     \mathbb{E}( \hat{w}_{2j} \tilde{\theta}_j - \theta_j )^2  & \leq 2 \mathbb{E} (\hat{w}_{2 j} \tilde{\theta}_j)^2+2 \theta_j^2 \leq 2 \mathbb{E}(\tilde{\theta}_j^2) +2 \theta_j^2  =  4\theta_j^2 + \frac{2\sigma^2}{n}, 
  \end{split}
\end{equation}
where the second inequality follows from the fact that $|\hat{w}_{2j}| \leq 1$ as shown in (\ref{eq:siwa5}), and the last equality follows from $\mathbb{E}(\tilde{\theta}_j^2) = (\mathbb{E}\tilde{\theta}_j)^2+  \Var(\tilde{\theta}_j) = \theta_j^2 + \sigma^2/n$. This bound also does not require Gaussianity of $\tilde{\theta}_j$. 

The third bound is derived by exploiting the Gaussian assumption $\tilde{\theta}_j \sim N(\theta_j, \sigma^2/n)$. Applying Stein's identity \citep{Stein1981}, we have 
\begin{equation}\label{eq:bound3_step1}
  \begin{split}
     \mathbb{E}( \hat{w}_{2j} \tilde{\theta}_j - \theta_j )^2  &= \mathbb{E} \left[ \left( \tilde{\theta}_j^2 - \frac{\sigma^2}{n} \right) 1_{\{ |\tilde{\theta}_j | \leq \lambda_2 \sigma \}} \right] + \mathbb{E} \left[ \left( \frac{\lambda_2^4\sigma^4+2\lambda_2^2\sigma^2 \frac{\sigma^2}{n}}{\tilde{\theta}_j^2}+ \frac{\sigma^2}{n}\right) 1_{\{ |\tilde{\theta}_j | > \lambda_2 \sigma \}}\right] \\
       & =\mathbb{E} \left[ \left( \tilde{\theta}_j^2 - \frac{\sigma^2}{n} \right)  \right] + \mathbb{E} \left[ \left( \frac{\lambda_2^4\sigma^4+2\lambda_2^2\sigma^2 \frac{\sigma^2}{n}}{\tilde{\theta}_j^2}-\tilde{\theta}_j^2+ \frac{2\sigma^2}{n}\right) 1_{\{ |\tilde{\theta}_j | > \lambda_2 \sigma \}}\right]\\
       & \leq \theta_j^2 + \left(\frac{\lambda_2^4\sigma^4+2\lambda_2^2\sigma^2 \frac{\sigma^2}{n}}{\lambda_2^2 \sigma^2}-\lambda_2^2 \sigma^2+ \frac{2\sigma^2}{n}\right)\mathbb{P}\left( |\tilde{\theta}_j | > \lambda_2 \sigma \right)\\
       & = \theta_j^2 + \frac{4\sigma^2}{n}\mathbb{P}\left( |\tilde{\theta}_j | > \lambda_2 \sigma \right).\\
  \end{split}
\end{equation}
To bound the tail probability, we use the Taylor expansion argument in (A1$\cdot$3) of \cite{Donoho1994Ideal} (see also the inequality on line 18 of page 486 of \citet{Gao1998}), which yields 
\begin{equation*}
\begin{split}
     &   \mathbb{P}\left( |\tilde{\theta}_j | > \lambda_2 \sigma \right) = \mathbb{P}\left( \left|\frac{\tilde{\theta}_j}{\sigma/\sqrt{n}} \right| > \sqrt{n}\lambda_2 \right)  \leq \frac{2 \phi(\sqrt{n}\lambda_2)}{\sqrt{n}\lambda_2} + \frac{n\theta_j^2}{4\sigma^2}, 
\end{split}
\end{equation*}
where $\phi$ is the density function of the standard normal distribution. Substituting this bound into (\ref{eq:bound3_step1}) gives the third univariate risk bound: 
\begin{equation}\label{eq:bound3}
  \begin{split}
       &\mathbb{E}( \hat{w}_{2j} \tilde{\theta}_j - \theta_j )^2 \leq \theta_j^2 + \frac{4\sigma^2}{n}\left( \frac{2 \phi(\sqrt{n}\lambda_2)}{\sqrt{n}\lambda_2} + \frac{n\theta_j^2}{4\sigma^2} \right) = 2\theta_j^2 + \frac{8\sigma^2\phi(\sqrt{n}\lambda_2)}{n\sqrt{n}\lambda_2}.
  \end{split}
\end{equation}

\subsubsection{Finalizing the proof}

To complete the proof, we follow the approach in \cite{Donoho1994Ideal} by separately upper bounding the univariate risk $\mathbb{E}( \hat{w}_{2j} \tilde{\theta}_j - \theta_j )^2$ in three different regimes of $|\theta_j|$. 

The first regime is $\theta_j^2 \geq \frac{\sigma^2\log p}{8n}$. Recall that $\lambda_2 = \sqrt{ 2\log p/n} $. Using the first bound (\ref{eq:bound1}), the univariate risk of the Adap estimator is upper bounded by
\begin{equation*}
  \mathbb{E}( \hat{w}_{2j} \tilde{\theta}_j - \theta_j )^2 \leq \left( \frac{\sqrt{2\log p} + 1}{\sqrt{n}} \right)^2\sigma^2 = \left( \frac{2\log p + 2\sqrt{2\log p} + 1}{n}  \right) \sigma^2.
\end{equation*}
Meanwhile, the $j$-th term in the optimal all-subset MA risk (\ref{eq:ideal_MA_risk}) is lower bounded by
\begin{equation}\label{eq:woaisss1}
  \frac{\theta_j^2 \sigma^2/n}{\theta_j^2 + \sigma^2/n} = \frac{ \sigma^2/n}{1 + \frac{\sigma^2/n}{\theta_j^2}} = \frac{\sigma^2}{n} \frac{ 1}{1 + \frac{\sigma^2/n}{\theta_j^2}} \geq \frac{\sigma^2}{n} \frac{\log p}{\log p +8}.
\end{equation}
Thus, the univariate risk ratio in this regime is upper bounded by 
\begin{equation}\label{eq:regime_1}
  \frac{\mathbb{E}( \hat{w}_{2j} \tilde{\theta}_j - \theta_j )^2}{\frac{\theta_j^2 \sigma^2/n}{\theta_j^2 + \sigma^2/n}} \leq \frac{(2\log p + 2\sqrt{2\log p} + 1)(\log p +8)}{\log p} \leq \begin{cases} \bar{C}, & p \text{\: is fixed} \\ [1+o(1)]2\log p , &  p \to \infty.\end{cases} 
\end{equation}

The second regime is $\frac{ 2\sigma^2}{n\log p} \leq \theta_j^2 < \frac{\sigma^2\log p}{8n}$. Applying the second upper bound (\ref{eq:bound2}), we have
\begin{equation}\label{eq:regime_2}
\begin{split}
    \frac{\mathbb{E}( \hat{w}_{2j} \tilde{\theta}_j - \theta_j )^2}{\frac{\theta_j^2 \sigma^2/n}{\theta_j^2 + \sigma^2/n}} & \leq  \frac{ 4\theta_j^2+2\sigma^2/n}{\frac{\theta_j^2 \sigma^2/n}{\theta_j^2 + \sigma^2/n}} = \frac{\left( 4\theta_j^2+2\sigma^2/n \right)\left( \theta_j^2 + \sigma^2/n \right)}{\theta_j^2 \sigma^2/n} \leq \frac{\left( 4\theta_j^2+2\sigma^2/n \right)^2}{\theta_j^2 \sigma^2/n} \\
     & =  \left( \frac{4|\theta_j|}{\sigma/\sqrt{n}} + \frac{2\sigma}{\sqrt{n}|\theta_j|} \right)^2 = \left( 4t+ \frac{2}{t} \right)^2\\
     & \leq \left( \sqrt{2 \log p} + \frac{8}{\sqrt{2 \log p}} \right)^2 \leq \begin{cases} \bar{C}, & p \text{\: is bounded} \\ [1+o(1)]2\log p , &  p \to \infty.\end{cases} 
\end{split}
\end{equation}
where the third inequality follows from $\sqrt{\frac{2}{\log p}} \leq  t \triangleq \frac{|\theta_j|}{\sigma/\sqrt{n}} \leq \frac{\sqrt{2\log p}}{4}$. It is worth noting that the upper bounds in \eqref{eq:regime_1} and \eqref{eq:regime_2} do not require the Gaussian assumption.

The last case is $0 \leq \theta_j^2 < \frac{ 2\sigma^2}{n\log p}$. By (\ref{eq:bound3}), the risk ratio satisfies
\begin{equation}\label{eq:mmpliu_1}
 \frac{\mathbb{E}( \hat{w}_{2j} \tilde{\theta}_j - \theta_j )^2}{\frac{1}{np}+\frac{\theta_j^2 \sigma^2/n}{\theta_j^2 + \sigma^2/n}} \leq\frac{\frac{8\sigma^2\phi(\sqrt{n}\lambda_2)}{n\sqrt{n}\lambda_2}+ 2\theta_j^2}{\frac{1}{np}+\frac{\theta_j^2 \sigma^2/n}{\theta_j^2 + \sigma^2/n}} \leq \frac{\frac{8\sigma^2\phi(\sqrt{n}\lambda_2)}{n\sqrt{n}\lambda_2}}{\frac{1}{np}}+ \frac{2\theta_j^2}{\frac{\theta_j^2 \sigma^2/n}{\theta_j^2 + \sigma^2/n}}.
\end{equation}
Since $\theta_j^2 < \frac{2\sigma^2}{n\log p}$, the second term on the right-hand side of \eqref{eq:mmpliu_1} is upper bounded by 
\begin{equation}\label{eq:mmpliu_2}
  \frac{2\theta_j^2}{\frac{\theta_j^2 \sigma^2/n}{\theta_j^2 + \sigma^2/n}} = \frac{2\left(\theta_j^2 + \sigma^2/n\right)}{\sigma^2/n} \leq 2+ \frac{4}{\log p}.
\end{equation}
Moreover, the first term on the right-hand side of \eqref{eq:mmpliu_1} admits the following expression:  
\begin{equation}\label{eq:mmpliu_3}
\begin{split}
   \frac{\frac{8\sigma^2\phi(\sqrt{n}\lambda_2)}{n\sqrt{n}\lambda_2}}{\frac{1}{np}} & = \frac{\frac{8\sigma^2\phi(\sqrt{2\log p})}{n\sqrt{2\log p}}}{\frac{1}{np}} = \frac{8p\sigma^2\phi(\sqrt{2\log p})}{\sqrt{2\log p}} \\
     & = \frac{8p\sigma^2 \frac{1}{\sqrt{2\pi}}\exp\left( - \frac{2\log p}{2} \right)}{\sqrt{2\log p}} = \frac{4\sigma^2}{\sqrt{\pi \log p}}.
\end{split}
\end{equation}
Substituting (\ref{eq:mmpliu_2}) and (\ref{eq:mmpliu_3}) into (\ref{eq:mmpliu_1}) yields
\begin{equation}\label{eq:regime_3}
  \frac{\mathbb{E}( \hat{w}_{2j} \tilde{\theta}_j - \theta_j )^2}{\frac{1}{np}+\frac{\theta_j^2 \sigma^2/n}{\theta_j^2 + \sigma^2/n}} \leq 2+ \frac{4}{\log p}  + \frac{4\sigma^2}{\sqrt{\pi \log p}} \leq \bar{C}.
\end{equation}

Combining the results in (\ref{eq:regime_1}), (\ref{eq:regime_2}), and (\ref{eq:regime_3}), we conclude that the univariate risk ratio is upper bounded by
\begin{equation*}
  \frac{\mathbb{E}( \hat{w}_{2j} \tilde{\theta}_j - \theta_j )^2}{\frac{1}{np}+\frac{\theta_j^2 \sigma^2/n}{\theta_j^2 + \sigma^2/n}} \leq \begin{cases} \bar{C}, & p \text{\: is fixed} \\ [1+o(1)]2\log p , &  p \to \infty.\end{cases} 
\end{equation*}
Multiplying both sides of the above display by $\frac{1}{np}+\frac{\theta_j^2 \sigma^2/n}{\theta_j^2 + \sigma^2/n}$ and summing over all $j$, the desired result follows. 

\subsection{Proof of Corollary~\ref{coro:unknown-variance}} 

Under the i.i.d. Gaussian assumption on $\epsilon_1,\ldots,\epsilon_n$ and Assumption~\ref{ass:complete}, the variance estimator satisfies $\hat{\sigma}^2=\frac{\sigma^2}{n-p} \chi_{n-p}^2$, and is independent of $\tilde{\theta}_j = n^{-1}\bpsi_j^{\top}\by$, $j = 1, \ldots, p$. Consequently, $\tilde{\theta}_j | \hat{\sigma}^2 \sim N(\theta_j, \sigma^2/n)$. Using properties of the chi‑square distribution, we obtain  
\begin{equation}\label{eq:chi_xx}
  \mathbb{E}[\hat{\sigma}^2]=\sigma^2, \quad \mathbb{E}[\hat{\sigma}]=\sigma \sqrt{\frac{2}{n-p}} \frac{\Gamma\left(\frac{n-p+1}{2}\right)}{\Gamma\left(\frac{n-p}{2}\right)}=\sigma+O\left[(n-p)^{-1}\right] = \sigma+O\left(n^{-1}\right),
\end{equation}
where the second equality follows from Stirling's formula and the condition $n-p\asymp n$. Moreover, $\operatorname{Var}(\hat{\sigma})=O(n^{-1})$ and $\mathbb{E}[(\hat{\sigma}-\sigma)^2]=O(n^{-1})$. Define $\hat{w}_{2j}(s) \triangleq ( 1 - \lambda_2^2 s/\tilde{\theta}_j^2 )_+$ for $s>0$, $j = 1, \ldots, p$. In analogy with (\ref{eq:MMA_cor1})--(\ref{eq:siwa5}), the $j$-th element of $\tilde{\bw}_2$ can be expressed as $\tilde{w}_{2j} = \hat{w}_{2j}(\hat{\sigma}^2)$. Since $\hat{\sigma}^2$ is independent of $\tilde{\theta}_j$, the proof of Corollary~\ref{coro:unknown-variance} follows the same arguments as those in Section~\ref{sec:a:proof_upper}. 

For the regime $\theta_j^2 \geq \frac{\sigma^2 \log p}{8 n}$, conditioning on $\hat{\sigma}^2$ and applying inequality (\ref{eq:bound1}) yields  
\begin{equation*}
\begin{split}
   \mathbb{E}\left[\left(\tilde{w}_{2 j} \tilde{\theta}_j-\theta_j\right)^2\right]  & = \mathbb{E} \left\{ \mathbb{E}\left[\left(\hat{w}_{2 j}(\hat{\sigma}^2) \tilde{\theta}_j-\theta_j\right)^2 \mid \hat{\sigma}^2\right] \right\} \leq \mathbb{E}\left(\frac{\sigma}{\sqrt{n}}+\lambda_2 \hat{\sigma}\right)^2 \\
     & =  \frac{\sigma^2}{n}+2 \lambda_2 \frac{\sigma}{\sqrt{n}} \mathbb{E}[\hat{\sigma}]+\lambda_2^2 \mathbb{E}\left[\hat{\sigma}^2\right] =\sigma^2\left(\lambda_2+\frac{1}{\sqrt{n}}\right)^2+2 \lambda_2 \frac{\sigma}{\sqrt{n}}(\mathbb{E}[\hat{\sigma}]-\sigma)\\
     & \leq  \left( \frac{\sqrt{2\log p} + 1}{\sqrt{n}} \right)^2\sigma^2  + O\left(\frac{\sqrt{\log p}}{n^2}\right), 
\end{split}
\end{equation*}
where the last line uses (\ref{eq:chi_xx}) and $\lambda_2= \sqrt{ 2\log p/n}$. Combining this result with (\ref{eq:woaisss1}), we conclude that the additional term $O(\sqrt{\log p}/n^2)$ does not affect the risk ratio. Hence, (\ref{eq:regime_1}) continues to hold provided that $n - p\asymp n$. 

For the regime $\frac{2 \sigma^2}{n \log p} \leq \theta_j^2<\frac{\sigma^2 \log p}{8 n}$, we have the upper bound  
$$
\mathbb{E}\left[\left(\hat{w}_{2 j}(s) \tilde{\theta}_j-\theta_j\right)^2\right] \leq 4 \theta_j^2+\frac{2 \sigma^2}{n}, \quad \forall s \geq 0,
$$
which is completely independent of $s$. Therefore, the upper bounds in (\ref{eq:regime_2}) also remains valid provided that $n - p\asymp n$.

For the regime $0 \leq \theta_j^2<\frac{2 \sigma^2}{n \log p}$, given any $\hat{\sigma}^2 = s> 0$, based on the similar analysis in (\ref{eq:bound3_step1}) and (\ref{eq:bound3}), we have 
\begin{equation}\label{eq:zoulle1}
  \begin{split}
     \mathbb{E}\left[\left(\hat{w}_{2 j}(s) \tilde{\theta}_j-\theta_j\right)^2\right] & \leq \theta_j^2 + \frac{4\sigma^2}{n}\mathbb{P}\left( |\tilde{\theta}_j | > \lambda_2 \sqrt{s} \right) \leq \theta_j^2 + \frac{4\sigma^2}{n}\left[  \frac{2 \phi(\sqrt{ns/\sigma^2}\lambda_2)}{\sqrt{ns/\sigma^2}\lambda_2} + \frac{n\theta_j^2}{4\sigma^2} \right] \\
     & = 2 \theta_j^2 + \frac{8\sigma^2 \phi(\sqrt{ns/\sigma^2}\lambda_2)}{n\sqrt{ns/\sigma^2}\lambda_2}. 
  \end{split}
\end{equation}
Define the event $\mathcal{A}=\{|\hat{\sigma}^2-\sigma^2| \leq \sigma^2/\log p\}$. Since $\operatorname{Var}(\hat{\sigma}^2) = 2 \sigma^4 / (n - p)$, Chebyshev's inequality yields  $\mathbb{P}(\mathcal{A}^c) \leq \frac{2 \sigma^4 / (n - p)}{\sigma^4 /(\log p)^2}\leq C_0 (\log p)^2/n $. On $\mathcal{A}$, we have $\hat{\sigma} / \sigma \in[\sqrt{1-1 / \log p}, \sqrt{1+1 / \log p}]$. Combining with (\ref{eq:zoulle1}), we have the upper bound 
\begin{equation*}
\begin{split}
   \mathbb{E}\left[\left(\tilde{w}_{2 j} \tilde{\theta}_j-\theta_j\right)^2\right] & = \mathbb{E} \left\{ \mathbb{E}\left[\left(\hat{w}_{2 j}(\hat{\sigma}^2) \tilde{\theta}_j-\theta_j\right)^2 \mid \hat{\sigma}^2\right] 1_{\mathcal{A}} \right\} + \mathbb{E} \left\{ \mathbb{E}\left[\left(\hat{w}_{2 j}(\hat{\sigma}^2) \tilde{\theta}_j-\theta_j\right)^2 \mid \hat{\sigma}^2\right] 1_{\mathcal{A}^c} \right\}\\
     & \leq \mathbb{E} \left\{ \left [2 \theta_j^2 + \frac{8\sigma^2 \phi(\sqrt{ns/\sigma^2}\lambda_2)}{n\sqrt{ns/\sigma^2}\lambda_2} \right] 1_{\mathcal{A}} \right\} + \left( \theta_j^2 + \frac{4\sigma^2}{n} \right) \mathbb{P}(\mathcal{A}^c)\\
     & \leq 2 \theta_j^2 + \frac{4\sigma^2}{n \sqrt{\pi}} \frac{\exp(- \log p + 1)}{\sqrt{\log p -1}} + C_0 \left( \theta_j^2 + \frac{4\sigma^2}{n} \right) \frac{(\log p)^2}{n}\\
     & \leq C \left[ \theta_j^2 + \frac{1}{np \sqrt{\log p}} + \left(\frac{\log p}{n} \right)^2 \right]. 
\end{split}
\end{equation*}
Then, using the similar proof technique as that in (\ref{eq:mmpliu_1})--(\ref{eq:regime_3}), we see that the risk ratio in this regime is upper bounded by a constant provided $p(\log p)^2 = O(n)$.

\subsection{Proof of Corollary~\ref{coro:non-Gaussian}}

Note that the proofs of the univariate bounds in \eqref{eq:regime_1} and \eqref{eq:regime_2} only rely on the identities $ \mathbb{E}(\tilde{\theta}_j) = \theta_j$ and $\Var(\tilde{\theta}_j) + \sigma^2/n$, which hold under Assumption~\ref{ass:4_moment} for any orthogonal basis satisfying $\| \bpsi_j \|^2 = n$. 

\subsubsection{Proof of (\ref{eq:universal_case})}

The proof of \eqref{eq:universal_case} follows a similar argument to that used in Section~\ref{sec:a:proof_upper}. Similar to \eqref{eq:siwa5}, the $j$-th element of the estimated weight vector $\hat{\bw}_c$ is given by $\hat{w}_{cj} = ( 1 - \lambda_c^2 \sigma^2/\tilde{\theta}_j^2 )_+$, where $\lambda_c = \sqrt{ (c\log p)/n }$. Note that the univariate risk bounds in \eqref{eq:bound1} and \eqref{eq:bound2} remain valid under Assumption~\ref{ass:4_moment} with $\lambda_2$ replaced by $\lambda_c$. 

First, consider the regime $\theta_j^2 \geq \frac{c \sigma^2\log p}{16 n}$. Following the similar argument as in \eqref{eq:regime_1}, we have 
\begin{equation}\label{eq:regime_u1}
  \frac{\mathbb{E}( \hat{w}_{cj} \tilde{\theta}_j - \theta_j )^2}{\frac{\theta_j^2 \sigma^2/n}{\theta_j^2 + \sigma^2/n}} \leq \frac{(c\log p + 2\sqrt{c\log p} + 1)(c\log p +16)}{c\log p} \leq \begin{cases} \bar{C}, & p \text{\: is fixed} \\ [1+o(1)]c\log p , &  p \to \infty.\end{cases} 
\end{equation}
In the regime $\frac{4\sigma^2}{c n \log p} \leq \theta_j^2 \leq \frac{c \sigma^2\log p}{16 n}$, following the similar argument as in \eqref{eq:regime_2}, we obtain 
\begin{equation}\label{eq:regime_u2}
\begin{split}
    \frac{\mathbb{E}( \hat{w}_{cj} \tilde{\theta}_j - \theta_j )^2}{\frac{\theta_j^2 \sigma^2/n}{\theta_j^2 + \sigma^2/n}} & \leq \left( \sqrt{c \log p} + \frac{8}{\sqrt{c \log p}} \right)^2 \leq \begin{cases} \bar{C}, & p \text{\: is fixed} \\ [1+o(1)]c\log p , &  p \to \infty.\end{cases} 
\end{split}
\end{equation}
Finally, consider the regime $0 \leq \theta_j^2 < \frac{4\sigma^2}{c n \log p}$. We first establish the following univariate risk bound: 
\begin{equation}\label{eq:hess3}
  \begin{split}
     \mathbb{E}( \hat{w}_{cj} \tilde{\theta}_j - \theta_j )^2  &  = \mathbb{E} \left[ \theta_j^2 1_{\{ |\tilde{\theta}_j| \leq \lambda_c \sigma \}} + \left( \tilde{\theta}_j - \theta_j - \frac{\lambda_c^2 \sigma^2}{\tilde{\theta}_j} \right)^2 1_{\{ |\tilde{\theta}_j| > \lambda_c \sigma \}}  \right] \\
       &  \leq \theta_j^2 + 2 \mathbb{E} \left[ ( \tilde{\theta}_j - \theta_j )^2 1_{\{ |\tilde{\theta}_j| > \lambda_c \sigma \}} \right] + 2 \lambda_c^2 \sigma^2 \mathbb{P}(|\tilde{\theta}_j| > \lambda_c \sigma)\\
       & \leq \theta_j^2 + 2 \sqrt{\mathbb{E}( \tilde{\theta}_j - \theta_j )^4} \sqrt{\mathbb{P}(|\tilde{\theta}_j| > \lambda_c \sigma)} + 2 \lambda_c^2 \sigma^2 \mathbb{P}(|\tilde{\theta}_j| > \lambda_c \sigma), 
  \end{split}
\end{equation}
where the first inequality follows from $(a+b)^2\leq 2a^2+2b^2$ and the fact that $\tilde{\theta}_j^2 > \lambda_c^2 \sigma^2$ when $1_{\{ |\tilde{\theta}_j| > \lambda_c \sigma \}} = 1$, and last inequality follows from the Cauchy–Schwarz inequality. Recall that $\tilde{\theta}_j - \theta_j = n^{-1} \bpsi_j^{\top}\beps$ under Assumption~\ref{ass:complete}. Since $\| \bpsi_j \|^2 = n$, $\max_{j,i}|\psi_{ji}|  \leq C_0$, and $\epsilon_i$'s are independent zero mean random errors with $\mathbb{E} \epsilon_i^4 \leq \nu$, a standard fourth-moment calculation gives 
\begin{equation}\label{eq:hess2}
  \mathbb{E}( \tilde{\theta}_j - \theta_j )^4 \leq \frac{C}{n^2}, 
\end{equation}
where $C>0$ is a constant. Next, we bound the tail probability $\mathbb{P}(|\tilde{\theta}_j| > \lambda_c \sigma)$. Since $\tilde{\theta}_j - \theta_j = n^{-1} \bpsi_j^{\top}\beps = n^{-1} \sum_{i=1}^{n}\psi_{ji}\epsilon_i$, the Berry--Esseen theorem gives  
\begin{equation}\label{eq:BS_bound}
\begin{split}
  \sup_{t}\left| \mathbb{P}\left( \frac{\sqrt{n}(\tilde{\theta}_j-\theta_j) }{ \sigma} \leq t \right)  - \Phi(t)   \right| & =\sup_{t}\left| \mathbb{P}\left( \frac{n^{-1} \bpsi_j^{\top}\beps}{n^{-1} \| \bpsi_j \| \sigma} \leq t \right)  - \Phi(t)   \right| \\
     &  \leq  \frac{C \max_{1 \leq i \leq n}|\psi_{ji}| }{\sqrt{n}} \leq \frac{C}{\sqrt{n}}. 
\end{split}
\end{equation}
Therefore, when $c \geq \frac{4}{\log p}$ and $|\tilde{\theta}_j|>\lambda_c \sigma$, we have $|\sqrt{n}\tilde{\theta}_j/\sigma|>\sqrt{c \log p} \geq 2$ and $|\sqrt{n}\theta_j/\sigma| < 2/\sqrt{c \log p} \leq 1$ in the current regime. Consequently, $|\sqrt{n}\tilde{\theta}_j/\sigma - \sqrt{n}\theta_j/\sigma| \geq \frac{\sqrt{c \log p}}{2}$. Hence,  
\begin{equation}\label{eq:hess1}
\begin{aligned}
\mathbb{P}(|\tilde{\theta}_j|>\lambda_c \sigma) & \leq \mathbb{P}\left( \left| \frac{\sqrt{n}(\tilde{\theta}_j-\theta_j) }{ \sigma} \right| \geq  \frac{\sqrt{c \log p}}{2} \right) \leq  2\left[1-\Phi\left(\frac{\sqrt{c \log p}}{2}\right)\right]+\frac{C}{\sqrt{n}} \\
& \leq \frac{4}{\sqrt{2 \pi c \log p}} \exp \left(-\frac{c \log p}{8}\right)+\frac{C}{\sqrt{n}} \\
& =\frac{C}{\sqrt{c \log p}} p^{-c / 8}+\frac{C}{\sqrt{n}}, 
\end{aligned}
\end{equation}
where the second inequality follows from (\ref{eq:BS_bound}), and the third inequality follows from $1-\Phi(x) \leq \frac{\phi(x)}{x}$. Inserting (\ref{eq:hess2}) and (\ref{eq:hess1}) into (\ref{eq:hess3}), we obtain 
\begin{equation*}
  \begin{split}
      \mathbb{E}( \hat{w}_{cj} \tilde{\theta}_j - \theta_j )^2 & \leq \theta_j^2+\frac{C}{n} \sqrt{\frac{p^{-c / 8}}{\sqrt{\log p}}+\frac{1}{\sqrt{n}}}+\frac{C \log p}{n}\left(\frac{p^{-c / 8}}{\sqrt{\log p}}+\frac{1}{\sqrt{n}}\right).\\
  \end{split}
\end{equation*}
Using the similar strategies in (\ref{eq:mmpliu_1}), we have the risk ratio in this regime is upper bounded by 
\begin{equation}\label{eq:chuzz1}
  \begin{split}
     \frac{\mathbb{E}( \hat{w}_{cj} \tilde{\theta}_j - \theta_j )^2}{\frac{1}{np}+\frac{\theta_j^2 \sigma^2/n}{\theta_j^2 + \sigma^2/n}} & \leq \frac{\theta_j^2}{\frac{\theta_j^2 \sigma^2/n}{\theta_j^2 + \sigma^2/n}} + C p \sqrt{\frac{p^{-c / 8}}{\sqrt{\log p}}+\frac{1}{\sqrt{n}}} + C p \log p \left(\frac{p^{-c / 8}}{\sqrt{\log p}}+\frac{1}{\sqrt{n}}\right)\\
       &  \leq 1 + \frac{4}{c \log p} + \frac{C p^{1 - \frac{c}{16}}}{(\log p)^{\frac{1}{4}}} + \frac{C p}{n^{\frac{1}{4}}} + C p^{1 - \frac{c}{8} } \sqrt{\log p} + \frac{C p \log p}{\sqrt{n}}\\
       & \leq \bar{C},
  \end{split}
\end{equation}
provided that $c \geq 16$ and $p=O(n^{1/4})$. Summing over all $j$ yields the desired results. 

\subsubsection{Proof of (\ref{eq:universal_case_2})}

The upper bounds in \eqref{eq:regime_u1} and \eqref{eq:regime_u2} also hold under the Bernstein condition on $\epsilon_i$. Therefore, to prove \eqref{eq:universal_case_2}, it remains to establish an upper bound for the regime $0 \leq \theta_j^2 < \frac{4\sigma^2}{c n \log p}$. More specifically, it suffices to bound the tail probability $\mathbb{P}(|\tilde{\theta}_j| > \lambda_c \sigma)$ appearing in the last line of \eqref{eq:hess3}. 

We write $\tilde{\theta}_j - \theta_j = n^{-1} \sum_{i=1}^{n}\psi_{ji}\epsilon_i = \sum_{i=1}^{n} \omega_i $, where $\omega_i \triangleq n^{-1} \psi_{ji}\epsilon_i$. Note that $\mathbb{E}\omega_i=0$ and $\mathbb{E}(\omega_i^2)=\frac{\psi_{j i}^2}{n^2} \sigma^2$. For all $k \geq 2$, we have 
\begin{equation*}
  \mathbb{E}\left|\omega_i\right|^k=\frac{\left|\psi_{j i}\right|^k}{n^k} \mathbb{E}\left|\epsilon_i\right|^k \leq \frac{1}{2} k!\frac{\left|\psi_{j i}\right|^k}{n^k} B^{k-2} \sigma^2 \leq \frac{1}{2} k!\frac{\psi_{j i}^2 \sigma^2}{n^2}\left(\frac{C_0 B}{n}\right)^{k-2}. 
\end{equation*}
Then based on the Bernstein inequality, we have 
\begin{equation*}
  \begin{split}
     \mathbb{P}\left( \left|\tilde{\theta}_j - \theta_j \right| \geq t \right) & = \mathbb{P}\left( \left| \sum_{i=1}^{n}\omega_i \right| \geq t \right) \leq 2 \exp \left\{-\frac{t^2}{2 \sigma^2 / n+2\left(C_0 B / n\right) t}\right\} = 2 \exp \left(-\frac{n t^2}{2 \sigma^2+2 C_0 B t}\right).
  \end{split}
\end{equation*}
When $c \geq \frac{4}{\log p}$, we see that $\{|\tilde{\theta}_j| > \lambda_c \sigma \}$ implies $\{ |\tilde{\theta}_j - \theta_j | \geq \frac{\lambda_c \sigma}{2} \}$. Therefore, when $n$ is sufficiently large, there exists a sufficiently small $\delta>0$ such that  
\begin{equation*}
  \mathbb{P} \left( |\tilde{\theta}_j| > \lambda_c \sigma \right) \leq \mathbb{P}\left( \left|\tilde{\theta}_j - \theta_j \right| \geq \frac{\lambda_c \sigma}{2} \right) \leq 2 \exp\left( - \frac{c\log p}{4(2+\delta)} \right) = 2p^{-\frac{c}{4(2+\delta)}}. 
\end{equation*}
Using the inequality (\ref{eq:hess3}) again, we have 
\begin{equation*}
  \mathbb{E}( \hat{w}_{cj} \tilde{\theta}_j - \theta_j )^2 \leq \theta_j^2 + \frac{C}{n}p^{-\frac{c}{8(2+\delta)}} + \frac{C}{n} (\log p) p^{-\frac{c}{4(2+\delta)}}. 
\end{equation*}
Similar to (\ref{eq:chuzz1}), the risk ratio is then upper bounded by 
\begin{equation*}
  \begin{split}
     \frac{\mathbb{E}( \hat{w}_{cj} \tilde{\theta}_j - \theta_j )^2}{\frac{1}{np}+\frac{\theta_j^2 \sigma^2/n}{\theta_j^2 + \sigma^2/n}} & \leq 1 + \frac{4}{c \log p} + Cp^{1-\frac{c}{8(2+\delta)}} + C (\log p) p^{1-\frac{c}{4(2+\delta)}} \leq \bar{C}, 
  \end{split}
\end{equation*}
provided that $c > 16$. Summing over all $j$ yields the desired results.

\newpage
\bibliographystyle{apalike}
\bibliography{MAbibfile}

\end{document}